%% file: BurgerIozziWienhard.tex
\documentclass[12pt]{amsart}
\usepackage{amscd,amssymb}
\usepackage[mathscr]{eucal}
\usepackage[english]{babel}
\input pstricks
\input pst-node

\usepackage{pstricks, pst-node}
\input xy
\usepackage[all]{xy}

\usepackage{amsmath}

\usepackage[latin1]{inputenc}%

\newlength{\itemlaenge}

\begin{document}

%%%%%%%%%% Definitionen von Umgebungen %%%%%%%%%%%%%%%%%%%%%%%%

\newtheoremstyle{mytheorem}% name
  {}%      Space above, empty = `usual value'
  {}%      Space below
  {\slshape}% Body font
  {}%         Indent amount (empty = no indent, \parindent = para indent)
  {\scshape}% Thm head font
  {.}%        Punctuation after head
  { }%     Space after thm head: " " = normal interword space;
        %       \newline = linebreak
  {}% Thm head spec

\newtheoremstyle{mydefinition}% name
  {}%      Space above, empty = `usual value'
  {}%      Space below
  {\upshape}% Body font
  {}%         Indent amount (empty = no indent, \parindent = para indent)
  {\scshape}% Thm head font
  {.}%        Punctuation after thm head
  { }%     Space after thm head: " " = normal interword space;
        %       \newline = linebreak
  {}% Thm head spec

\theoremstyle{mytheorem}
\newtheorem{lemma}{Lemma}[section]
\newtheorem{prop}[lemma]{Proposition}
\newtheorem*{prop*}{Proposition}
\newtheorem{prop_intro}{Proposition}
\newtheorem{cor}[lemma]{Corollary}
\newtheorem{cor_intro}[prop_intro]{Corollary}
\newtheorem{thm}[lemma]{Theorem}
\newtheorem{thm_intro}[prop_intro]{Theorem}
\newtheorem*{thm*}{Theorem}
\theoremstyle{mydefinition}
\newtheorem{rem}[lemma]{Remark}
\newtheorem*{rem*}{Remark}
\newtheorem{rem_intro}[prop_intro]{Remark}
\newtheorem{rems_intro}[prop_intro]{Remarks}
\newtheorem*{notation*}{Notation}
\newtheorem*{warning*}{Warning}
\newtheorem{rems}[lemma]{Remarks}
\newtheorem{defi}[lemma]{Definition}
\newtheorem*{defi*}{Definition}
\newtheorem{defi_intro}[prop_intro]{Definition}
\newtheorem{defis}[lemma]{Definitions}
\newtheorem{exo}[lemma]{Example}
\newtheorem{exos}[lemma]{Examples}
\newtheorem{exo_intro}[prop_intro]{Example}
\newtheorem{exos_intro}[prop_intro]{Examples}

\numberwithin{equation}{section}

\newcommand{\bqn}{\begin{eqnarray*}}
\newcommand{\eqn}{\end{eqnarray*}}
\newcommand{\bq}{\begin{eqnarray}}
\newcommand{\eq}{\end{eqnarray}}
\newcommand{\be}{\begin{enumerate}}
\newcommand{\ee}{\end{enumerate}}
\newcommand{\ba}{\begin{aligned}}
\newcommand{\ea}{\end{aligned}}

\newcommand{\bibURL}[1]{{\unskip\nobreak\hfil\penalty50{\tt#1}}}

\def\ti{-\allowhyphens}
\newcommand{\thismonth}{\ifcase\month % case 0 --- impossible!
  \or January\or February\or March\or April\or May\or June%
  \or July\or August\or September\or October\or November%
  \or December\fi}
\newcommand{\thismonthyear}{{\thismonth} {\number\year}}
\newcommand{\thisdaymonthyear}{{\number\day} {\thismonth} {\number\year}}
%
%
%
%%%%%%%%%%%%%%%%ù Definitionen von Symbolen %%%%%%%%%%%%%%%%%%%%%%ù

\renewcommand{\AA}{{\mathbb A}}
\newcommand{\BB}{{\mathbb B}}
\newcommand{\CC}{{\mathbb C}}
\newcommand{\DD}{{\mathbb D}}
\newcommand{\FF}{{\mathbb F}}
\newcommand{\HH}{{\mathbb H}}
\newcommand{\GG}{{\mathbb G}}
\newcommand{\KK}{{\mathbb K}}
\newcommand{\NN}{{\mathbb N}}
\newcommand{\PP}{{\mathbb P}}
\newcommand{\QQ}{{\mathbb Q}}
\newcommand{\RR}{{\mathbb R}}
\newcommand{\TT}{{\mathbb T}}
\newcommand{\ZZ}{{\mathbb Z}}

\newcommand{\Aa}{{\mathcal A}}
\newcommand{\Bb}{{\mathcal B}}
\newcommand{\Cc}{{\mathcal C}}
\newcommand{\Dd}{{\mathcal D}}
\newcommand{\Ee}{{\mathcal E}}
\newcommand{\Ff}{{\mathcal F}}
\newcommand{\Hh}{{\mathcal H}}
\newcommand{\Ii}{{\mathcal I}}
\newcommand{\Jj}{{\mathcal J}}
\newcommand{\Kk}{{\mathcal K}}
\newcommand{\Ll}{{\mathcal L}}
\newcommand{\Mm}{{\mathcal M}}
\newcommand{\Nn}{{\mathcal N}}
\newcommand{\Oo}{{\mathcal O}}
\newcommand{\Qq}{{\mathcal Q}}
\newcommand{\Pp}{{\mathcal P}}
\newcommand{\Rr}{{\mathcal R}}
\newcommand{\Ss}{{\mathcal S}}
\newcommand{\Tt}{{\mathcal T}}
\newcommand{\Vv}{{\mathcal V}}
\newcommand{\Xx}{{\mathcal X}}
\newcommand{\Yy}{{\mathcal Y}}
\newcommand{\Zz}{{\mathcal Z}}

\newcommand{\fraka}{{\mathfrak a}}
\newcommand{\frakc}{{\mathfrak c}}
\newcommand{\frake}{{\mathfrak e}}
\newcommand{\frakg}{{\mathfrak g}}
\newcommand{\frakh}{{\mathfrak h}}
\newcommand{\frakk}{{\mathfrak k}}
\newcommand{\frakl}{{\mathfrak l}}
\newcommand{\frakm}{{\mathfrak m}}
\newcommand{\frakn}{{\mathfrak n}}
\newcommand{\frako}{{\mathfrak o}}
\newcommand{\frakp}{{\mathfrak p}}
\newcommand{\frakq}{{\mathfrak q}}
\newcommand{\frakr}{{\mathfrak r}}
\newcommand{\fraks}{{\mathfrak s}}
\newcommand{\frakt}{{\mathfrak t}}
\newcommand{\fraku}{{\mathfrak u}}

\newcommand{\frakpp}{{\mathfrak p_+}}
\newcommand{\frakpm}{{\mathfrak p_-}}
\newcommand{\frakkc}{{{\mathfrak k}_\CC}}
\newcommand{\frakgc}{{{\mathfrak g}_\CC}}
\newcommand{\frakpc}{{{\mathfrak p}_\CC}}
\newcommand{\frakB}{{\mathfrak B}}

\newcommand{\tk}{{\tt k}}

\newcommand{\Stab}{\operatorname{Stab}}
\newcommand{\sign}{\operatorname{sign}}
\newcommand{\aut}{\operatorname{Aut}}
\newcommand{\Int}{\operatorname{Int}}
\newcommand{\End}{\operatorname{End}}
\newcommand{\Is}{\operatorname{Is}}
\newcommand{\SU}{\operatorname{SU}}
\newcommand{\SL}{\operatorname{SL}}
\newcommand{\PU}{\operatorname{PU}}
\newcommand{\mo}{\operatorname{mod}}
\newcommand{\ko}{\operatorname{k}}
\newcommand{\Image}{\operatorname{Image}}
\newcommand{\Lie}{\operatorname{Lie}}
\newcommand{\codim}{\operatorname{codim}}
\newcommand{\esssup}{\operatorname{ess\,sup}}
\newcommand{\supp}{\operatorname{supp}}
\newcommand{\length}{\operatorname{length}}
\newcommand{\rad}{\operatorname{Rad}}
\newcommand{\essgr}{\operatorname{Ess\,Gr}}
\newcommand{\essim}{\operatorname{Ess\,Im}}
\newcommand{\QHB}{{\operatorname{QH}_\Bb}}
\newcommand{\QHch}{{\operatorname{QH}_\mathrm{c}^\mathrm{h}}}

\renewcommand{\a}{\alpha}
\newcommand{\e}{\epsilon}
\newcommand{\eps}{\epsilon}
\renewcommand{\b}{\beta}
\newcommand{\g}{\gamma}
\newcommand{\G}{\Gamma}
\renewcommand{\L}{\Lambda}
\renewcommand{\l}{\lambda}
\newcommand{\T}{{\rm T}}

\newcommand{\dD}{{\mathbf D}}
\newcommand{\gG}{{\mathbf G}}
\newcommand{\hH}{{\mathbf H}}
\newcommand{\pP}{{\mathbf P}}
\newcommand{\lL}{{\mathbf L}}
\newcommand{\qQ}{{\mathbf Q}}
\newcommand{\nN}{{\mathbf N}}
\newcommand{\uU}{{\mathbf U}}
\newcommand{\vV}{{\mathbf V}}
\newcommand{\wW}{{\mathbf W}}

\newcommand{\<}{\langle}
\renewcommand{\>}{\rangle}

\def\ol{\overline}

\def\h{{\rm H}}
\def\hb{{\rm H}_{\rm b}}
\def\ehb{{\rm EH}_{\rm b}}
\def\ha{{\rm H}_{(G,K)}}
\def\hc{{\rm H}_{\rm c}}
\def\hh{{\rm {\widehat H}}}
\def\hhc{{\rm {\widehat H}}_{\rm c}}
\def\hcb{{\rm H}_{\rm cb}}
\def\hhcb{{\rm {\widehat H}}_{\rm cb}}
\def\hhb{{\rm {\widehat H}}_{\rm b}}
\def\ehbc{{\rm EH}_{\rm cb}}
\def\linfty{L^\infty}
\def\linftyw{L^\infty_{\rm w*}}
\def\linftya{L^\infty_{\mathrm{w*,alt}}}
\def\cba{\mathrm{C}_{\mathrm{b,alt}}}
\def\balt{\mathcal{B}^\infty_{\mathrm{alt}}}
\def\la{L^\infty_{\mathrm{alt}}}
\def\cb{{\rm C}_{\rm b}}
\def\tb{\mathrm{t}_\mathrm{b}}
\def\Tb{\mathrm{T}_\mathrm{b}}
\def\binfty{\mathcal B^\infty_{\mathrm alt}}

\def\one{\mathbf{1\kern-1.6mm 1}}
\def\sous#1#2{{\raisebox{-1.5mm}{$#1$}\backslash \raisebox{.5mm}{$#2$}}}
\def\rest#1{{\raisebox{-.95mm}{$\big|$}\raisebox{-2mm}{$#1$}}}
\def\homeo#1{{\sl H\!omeo}^+\!\left(#1\right)}
\def\thomeo#1{\widetilde{{\sl \!H}\!omeo}^+\!\left(#1\right)}
\def\bu{\bullet}
\def\weak{weak-* }
\def\property{\textbf{\rm\textbf A}}
\def\cont{\mathcal{C}}
\def\id{{\it I\! d}}
\def\opposite{^{\rm op}}
\def\oddex#1#2{\left\{#1\right\}_{o}^{#2}}
\def\comp#1{{\rm C}^{(#1)}}
\def\ro{\varrho}
\def\ti{-\allowhyphens}
\def\lra{\longrightarrow}

\def\Homeo{\operatorname{Homeo}}
\def\fix{{\operatorname{Fix}}}
\def\Mat{{\operatorname{Mat}}}
\def\zent{{\operatorname{Zent}}}
\def\sym{{\operatorname{Sym}}}
\def\stab{{\operatorname{Stab}}}
\def\arg{{\operatorname{arg}}}
\def\HTP{{\operatorname{HTP}}}
\def\h2{{\operatorname{H_2}}}
\def\h1{{\operatorname{H_1}}}
\def\pr{{\operatorname{pr}}}
\def\rk{{\operatorname{r}}}
\def\tr{{\operatorname{tr}}}
\def\Tr{{\operatorname{Tr}}}
\def\codim{{\operatorname{codim}}}
\def\nt{{\operatorname{nt}}}
\def\d{{\operatorname{d}}}
\def\Gr{{\operatorname{Gr}}}
\def\id{{\operatorname{Id}}}
\def\ker{{\operatorname{ker}}}
\def\im{{\operatorname{im}}}
\def\rot{{\operatorname{Rot}}}
\def\rotk{{\operatorname{Rot}_\kappa}}
\def\rep{{\operatorname{Rep}}}
\def\Sp{{\operatorname{Sp}}}
\def\PU{{\operatorname{PU}}}
\def\PSL{{\operatorname{PSL}}}
\def\ad{\operatorname{ad}}
\def\Ad{\operatorname{Ad}}
\def\adg{\operatorname{ad}_\frakg}
\def\adp{\operatorname{ad}_\frakp}
\def\det{{\operatorname{det}}}
\def\deta{{\operatorname{det}_A}}
\def\diag{{\operatorname{diag}}}
\def\hom{\operatorname{Hom}}
\def\homc{\operatorname{Hom}_\mathrm{c}}
\def\I{\operatorname{I}}
\def\isp{\operatorname{Is}_{\<\cdot,\cdot\>}}
\def\isptwo{\operatorname{Is}_{\<\cdot,\cdot\>}^{(2)}}
\def\ispth{\operatorname{Is}_{\<\cdot,\cdot\>}^{(3)}}
\def\isf{\operatorname{Is}_F}
\def\isfi{\operatorname{Is}_{F_i}}
\def\isft{\operatorname{Is}_F^{(3)}}
\def\isfit{\operatorname{Is}_{F_i}^{(3)}}
\def\isftwo{\operatorname{Is}_F^{(2)}}
\def\lin{\operatorname{Lin}(L_+,L_-)}
\def\r{\operatorname{r}}
\def\Iso{\operatorname{Iso}}
\def\val{\operatorname{Val}}
\def\Sm{\operatorname{Sm}}
\def\diam{\operatorname{diam}}

\def\bg{B_\frakg}
\def\creg{C_{\rm reg}}
\def\c{\operatorname{c}}
\def\C{{\mathrm{C}}}
\def\cs{{\check S}}
\def\cst{{\check S}^{(3)}}
\def\gmodp{\gG(\RR)/\pP(\RR)}
\def\gmodq{\gG(\RR)/\qQ(\RR)}
\def\gr{{\rm g}_\RR^*}
\def\grs{{\rm g}_\RR^{(\Sigma,\star)}}
\def\gz{{\rm g}_\ZZ^*}
\def\gzs{{\rm g}_\ZZ^{(\Sigma,\star)}}
\def\ll{{\Ll_1,\Ll_2}}
\def\kahler{ K\"ahler }
\def\k{\kappa}
\def\kd{\kappa_\DD}
\def\kg{\kappa_G}
\def\kgb{\kappa_G^\mathrm{b}}
\def\kbr{\kappa_\RR^b}
\def\kbz{\kappa_\ZZ^b}
\def\kx{\kappa_\Xx}
\def\ku{\kappa_u}
\def\kr{\kappa_\RR}
\def\krlb{\kappa_{\RR,L}^\mathrm{b}}
\def\klb{\k_L^\mathrm{b}}
\def\kz{\kappa_\ZZ}
\def\kub{\kappa^{\rm b}_u}
\def\kgB{\kappa_{G,\mathrm{B}}}
\def\kgBb{\kappa_{G,\mathrm{B}}^{\rm b}}
\def\kgjb{\kappa_{G,J}^{\rm b}}
\def\kljb{\kappa_{L,J}^{\rm b}}
\def\khib{{\kappa_{H,i}^\mathrm{b}}}
\def\kxb{\kappa_\Xx^{\rm b}}
\def\kyb{\kappa_\Yy^{\rm b}}
\def\kdb{\kappa_\DD^{\rm b}}
\def\kgz{\kappa_{G,\ZZ}}
\def\kgbz{\kappa_{G,\ZZ}^{\rm b}}
\def\kuz{\kappa_{u,\ZZ}}
\def\kubz{\kappa_{u,\ZZ}^{\rm b}}
\def\kur{\kappa_{u,\RR}}
\def\kubr{\kappa_{u,\RR}^{\rm b}}
\def\kib{\kappa_i^{\rm b}}
\def\kibt{\tilde\kappa_i^{\rm b}}
\def\kszb{\kappa_{\Sigma,\ZZ}^\mathrm{b}}
\def\ksb{\kappa_{\Sigma}^\mathrm{b}}
\def\kxi{\kappa_{\Xx,i}}
\def\kxib{\kappa_{\Xx,i}^\mathrm{b}}
\def\kgib{\kappa_{G,i}^\mathrm{b}}
\def\oc{\overline c}
\def\om{\overline m}
\def\oG{\overline G}
\def\pii{\pi_i}
\def\piit{\tilde\pi_i}
\def\psupq{{\rm PSU}(p,q)}
\def\psuvi{{\rm PSU}\big(V,\<\cdot,\cdot\>_i\big)}
\def\rg{r_G}
\def\rhoba{\rho_\mathrm{b}^\ast}
\def\rast{\rho^\ast}
\def\PSU{\mathrm{PSU}}
\def\slv{{\rm SL}(V)}
\def\supq{{\rm SU}(p,q)}
\def\suv{{\rm SU}\big(V,\<\cdot,\cdot\>\big)}
\def\suvi{{\rm SU}\big(V,\<\cdot,\cdot\>_i\big)}
\def\jm{J^{(m)}}
\def\vc{V_\CC}
\def\vconj{\overline v}
\def\wconj{\overline w}

\def\binfty{\mathcal B^\infty_{\mathrm {alt}}}
\def\hcb{{\rm H}_{\rm cb}}
\def\to{\rightarrow}
\def\bu{\bullet}
\def\la{L^\infty_{\mathrm{alt}}}
\def\hb{{\rm H}_{\rm b}}
\def\hc{{\rm H}_{\rm c}}
\def\h{{\rm H}}

%%%%%%%%%%%%%%%% Definitions onlt for this paper

\def\bs{{\partial\Sigma}}
\def\gx{{G_\Xx}}
\def\gxi{{G_{\Xx_i}}}
\def\hcbps{\hcb^2\big(\pi_1(\Sigma),\RR\big)}
\def\hbps{\hb^2\big(\pi_1(\Sigma),\RR\big)}
\def\hpg{\operatorname{Hom}\big(\pi_1(\Sigma),G\big)}
\def\hmaxpg{\operatorname{Hom}_\mathrm{max}\big(\pi_1(\Sigma),G\big)}
\def\ps{{\pi_1(\Sigma)}}
\def\tsr{\mathrm{T}(\Sigma,\rho)}
\def\tksr{\mathrm{T}_\kappa(\Sigma,\rho)}

\def\sk{\mathsf k}
\renewcommand{\phi}{\varphi}

\def\eI{{^{I}\!\mathrm{E}}}
\def\eII{{^{I\!I}\!\mathrm{E}}}
\def\dI{{^{I}\!d}}
\def\dII{{^{I\!I}\!d}}
\def\fI{{^{I}\!F}}
\def\fII{{^{I\!I}\!F}}

\def\ind{\mathbf{i}}
\def\Ind{\mathrm{Ind}}
\def\No{N\raise4pt\hbox{\tiny o}\kern+.2em}
\def\no{n\raise4pt\hbox{\tiny o}\kern+.2em}
\def\bsl{\backslash}
%
%
%
%
%%%%%%%%%%%%%% Definitionen von Matrizen, etc. %%%%%%%%%%%%%%%%%%%%%%%

\newcommand{\mat}[4]{
\left(
\begin{array}{cc}
{\scriptstyle #1} & {\scriptstyle #2} \\ {\scriptstyle #3} & {\scriptstyle
#4}
\end{array}\right)}

\newcommand{\matt}[6]{
\left(
\begin{array}{cc}
{\scriptstyle #1} & {\scriptstyle #2} \\ {\scriptstyle #3} & {\scriptstyle
#4}\\{\scriptstyle #5} & {\scriptstyle #6}
\end{array}\right)}

\newcommand{\db}[4]{0=d\beta(y_{#1},y_{#2},y_{#3},y_{#4})
= \beta(y_{#1},y_{#2},y_{#3})+\beta(y_{#1},y_{#3},y_{#4})-\beta(y_{#1},y_{#2},y_{#4})
-\beta(y_{#2},y_{#3},y_{#4})}
%
%
%
%
%
%%%%%%%%%%%%%%%%%%% Titel, etc. %%%%%%%%%%%%%%%%%%ù

\title[Maximal Representations]{Surface Group Representations\\ with Maximal Toledo Invariant}
\author[M.~Burger]{Marc Burger}
\email{burger@math.ethz.ch}
\address{FIM, ETH Zentrum, R\"amistrasse 101, CH-8092 Z\"urich, Switzerland}
\author[A.~Iozzi]{Alessandra Iozzi}
\email{iozzi@math.ethz.ch}
\address{D-Math, ETH Zentrum, R\"amistrasse 101, CH-8092 Z\"urich, Switzerland}
\author[A.~Wienhard]{Anna Wienhard}
\email{wienhard@math.princeton.edu}
\address{Department of Mathematics, Princeton University, Fine Hall - Washington Road, Princeton, NJ 08540, USA}
\thanks{A.I. and A.W. were partially supported by FNS grant
  PP002-102765. A.W. was partially supported by the National Science
  Foundation under agreement No. DMS-0111298 and No. DMS-0604665}
%\keywords{ }
%\subjclass{ }

\date{\today}

\maketitle

\centerline{\it Dedicated to Domingo Toledo on his 60th birthday}
\vskip2cm

\begin{abstract} 
We develop the theory of maximal representations of the 
fundamental group $\pi_1(\Sigma)$ of a compact connected oriented surface 
$\Sigma$ (possibly with boundary) into Lie groups $G$ of Hermitian type.
For any homomorphism $\rho:\pi_1(\Sigma)\to G$, 
we define the {\it Toledo invariant} $\T(\Sigma,\rho)$, 
a numerical invariant which has both topological and analytical interpretations. 
We establish important properties of $\T(\Sigma,\rho)$, among which continuity, uniform boundedness on the 
representation variety, additivity
under connected sum of surfaces and congruence relations mod $\ZZ$.
We thus obtain information about the representation
variety as well as striking geometric properties of {\it maximal} representations, 
that is representations whose Toledo invariant achieves the maximum value. 

Moreover we establish properties of boundary maps associated to maximal 
representations which generalize naturally
monotonicity properties of semiconjugations of the circle.

We define a rotation number function for general locally compact groups and study it 
in detail for groups of Hermitian type.
Properties of the rotation number together with the existence of boundary maps 
lead to additional invariants for maximal representations and show 
that the subset of maximal representations is always real semialgebraic.  

In the case of surfaces without boundary some of the results were announced in \cite{Burger_Iozzi_Wienhard_ann}.
\end{abstract}

\maketitle
%
%
%
%
%
%%%%%%%%%%%%%%%%%%%%%%%%%%%%%%%%% Aufbau, Dateien %%%%%%%%%%%%%%%%%%%%%%%%%%%

%\input{titel}

\tableofcontents

%\vskip2cm
%\input{surface}
\input{intro}
\vskip1cm
\input{prelim}
\vskip1cm
\input{tolno}
\vskip1cm
\input{structure1}

\vskip1cm
\input{regularity}
\vskip1cm
\input{structure2}
\vskip1cm
\input{rotation}
\vskip1cm
\input{formula}
\vskip1cm
\input{examples}
\vskip1cm
\input{notation}
\vskip1cm
\bibliographystyle{amsplain}
\bibliography{refs}
%\vskip1cm
%\input{abstract}
%\vskip1cm
%\input{leben}
\providecommand{\bysame}{\leavevmode\hbox to3em{\hrulefill}\thinspace}
\providecommand{\MR}{\relax\ifhmode\unskip\space\fi MR }
\providecommand{\MRhref}[2]{%
  \href{http://www.ams.org/mathscinet-getitem?mr=#1}{#2}
}
\providecommand{\href}[2]{#2}
\vskip1cm
\end{document}

%% file: intro.tex
\section{Introduction}\label{sec:intro}
Let $\Sigma$ be an oriented compact surface with boundary $\partial\Sigma$ and 
let $G$ be a connected semisimple Lie group with finite center.
The problem of understanding the representation variety $\hpg$
has received considerable interest.  
A major theme is the problem of singling out special components of this
representation variety which should generalize Teichm\"uller space and 
then studying the geometric significance of the representations belonging to such components.

If $G$ is a split real group and $\bs=\emptyset$, 
there is a component whose global properties were studied by Hitchin
\cite{Hitchin};
the geometric significance of these representations was recently 
brought into the open on the one hand by the work of Labourie \cite{Labourie_anosov}
relating them to Anosov structures, 
and on the other hand by Fock and Goncharov \cite{Fock_Goncharov_local, Fock_Goncharov_convex}
studying them via the notion of positivity introduced by Lusztig \cite{Lusztig_positivity}.

When $G$ is of Hermitian type and $\bs=\emptyset$, 
one can define the {\it Toledo invariant} of a representation and 
hence the notion of {\it maximal} representation: 
these form a union of connected components of the representation variety.
The global properties of these components were investigated by Garc\'{\i}a-Prada,
Bradlow and Gothen using Higgs bundles 
\cite{Gothen, Bradlow_GarciaPrada_Gothen, Bradlow_GarciaPrada_Gothen_cacanerveuse} 
and the geometric properties of maximal representations were investigated 
by the authors in \cite{Burger_Iozzi_Wienhard_ann, Burger_Iozzi_Labourie_Wienhard}.

The purpose of this paper is to introduce and study the notion of Toledo invariant
when $\bs\neq\emptyset$ and 
investigate the structure of the corresponding maximal representations.  
The treatment includes the case in which $\partial\Sigma=\emptyset$
on which it sheds new light.

The main results are the structure theorem (Theorem~\ref{thm_intro:thm3}),
the regularity properties of boundary maps (Theorem~\ref{thm_intro:thm5})
and the formula for the Toledo invariant in terms of rotations numbers
(Theorem~\ref{thm_intro:thm8}).  For more background on the study of maximal 
representations we refer to 
\cite{Bradlow_GarciaPrada_Gothen, Bradlow_GarciaPrada_Gothen_cacanerveuse, Burger_Iozzi_difff, Burger_Iozzi_Labourie_Wienhard, Burger_Iozzi_Wienhard_ann, Goldman_thesis, Goldman_88, Goldman_torus, Hernandez, Toledo_89, Labourie_energy, Wienhard_mapping}.

\subsection{The Toledo invariant}\label{subsec:tol-inv}
Let $G$ be a group of Hermitian type (see \S~\ref{subsubsec:2.1.1}), 
so that in particular the associated  symmetric space $\Xx$ is Hermitian 
of noncompact type;  then $\Xx$ carries a unique Hermitian (normalized) metric 
of minimal holomorphic sectional curvature $-1$. 
The K\"ahler form $\omega_\Xx$ of this metric gives rise 
in the familiar way to a continuous class $\kg\in\hc^2(G,\RR)$ and,
owing to the isomorphism between bounded continuous and continuous cohomology
in degree two, to the {\it bounded K\"ahler class} $\kgb\in\hcb^2(G,\RR)$ 
(see \S~\ref{subsec:2.1}).  The bounded K\"ahler class is the source
of new invariants for representations and has been considered in 
\cite{Burger_Iozzi_supq, Burger_Iozzi_Wienhard_ann, Burger_Iozzi_Labourie_Wienhard, Burger_Iozzi_Wienhard_kahler, Burger_Iozzi_Wienhard_tight}.

Let $\Sigma$ be a connected oriented compact surface with boundary $\bs$,
and $\rho:\ps\to G$ a representation.
When $\bs=\emptyset$, the Toledo invariant is given by 
the evaluation of $\rho^\ast(\kg)$ on the fundamental class $[\Sigma]$,
\bqn
\tsr=\big\langle\rho^\ast(\kg),[\Sigma]\big\rangle\,.
\eqn
In the general case we obtain by pullback in bounded cohomology a bounded class
\bqn
\rho^\ast(\kgb)\in\hb^2\big(\ps,\RR\big)\cong\hb^2(\Sigma,\RR)\,.
\eqn
The canonical map $j_{\partial \Sigma}:\hb^2(\Sigma,\bs,\RR)\to\hb^2(\Sigma,\RR)$
from singular bounded cohomology relative to $\bs$ 
to singular bounded cohomology is an isomorphism (see (\ref{sec:prelim}.d) in \S~\ref{subsec:2.2}),
and we define
\bqn
\tsr=\big\langle j_{\partial \Sigma}^{-1}\rho^\ast(\kgb),[\Sigma,\bs]\big\rangle\,,
\eqn
where now $j_{\partial \Sigma}^{-1}\rho^\ast(\kgb)$ is considered as an ordinary relative class and $[\Sigma,\bs]$
is the relative fundamental class.  
The above construction applies to any class $\kappa\in\hcb^2(G,\RR)$ and 
we denote by $\tksr$ the resulting invariant.
This generalization will be useful when we consider integral classes.
This construction circumvents the fact that $\h^2(\Sigma,\RR)=0$ when $\bs\neq\emptyset$;
indeed in all cases $\hb^2(\Sigma,\RR)$  is infinite dimensional, 
provided $\chi(\Sigma)\leq-1$.

The basic properties of the Toledo invariant are summarized in the following

\begin{thm_intro}\label{thm_intro:thm1} Let $G$ be a group of Hermitian type and 
$\rho:\ps\to G$ a representation.  Then 
\be
\item $|\tsr|\leq|\chi(\Sigma)|\r_\Xx$, where $\r_\Xx$ is the rank of $\Xx$.
\item The map $\mathrm{T}(\Sigma,\,\cdot\,)$ is continuous on $\hpg$;  if $\bs=\emptyset$,
its range is finite, while if $\bs\neq\emptyset$ its range is the interval
\bqn
\big[-|\chi(\Sigma)|\r_\Xx,|\chi(\Sigma)|\r_\Xx\big]\,.
\eqn
\item If $\Sigma$ is the connected sum of two (connected)
surfaces $\Sigma_i$ along a separating loop, then 
\bqn
\tsr=\mathrm{T}(\Sigma_1,\rho_1)+\mathrm{T}(\Sigma_2,\rho_2)\,,
\eqn
\ee
where $\rho_i$ is the restriction of $\rho$ to $\pi_1(\Sigma_i)$.
\end{thm_intro}

\noindent
{\small{[Theorem~\ref{thm_intro:thm1} follows from Corollary~\ref{cor:3.3}, Proposition~\ref{prop:3.8}, Corollary~\ref{cor:3.10} and Proposition~\ref{prop:3.1}.]}}

Here and in the sequel, an essential role is played by Theorem~\ref{thm:3.2},
where we identify the Toledo invariant with an invariant defined in analytic terms,
introduced and studied in \cite{Burger_Iozzi_difff}. 
In view of Theorem~\ref{thm_intro:thm1}, we set the following

\begin{defi_intro}\label{def:maximal} A representation $\rho:\ps\to G$
is {\it maximal} if 
\bqn
\tsr=|\chi(\Sigma)|\r_\Xx\,.
\eqn
\end{defi_intro}

We denote by $\hmaxpg$ the subspace of the representation variety
consisting of maximal representations.
Notice that when $\bs=\emptyset$  
this subspace is a union of components of the representation variety,
while when $\bs\neq\emptyset$ the whole representation variety is connected
(see however Corollary~\ref{cor_intro:cor10}).

\subsection{Geometric properties of maximal representations}\label{subsec:geomprop}
Before treating the case of a general group of Hermitian type,
we state the structure theorem for maximal representation into $\mathrm{PU}(1,1)$
which generalizes to the case of surfaces with boundary 
Goldman's characterization of maximal representations 
when $\bs=\emptyset$ \cite{Goldman_thesis, Goldman_88}.

\begin{thm_intro}\label{thm_intro:thm2}  Let $\Sigma$ be a connected oriented surface
such that $\chi(\Sigma)\leq-1$.  A representation $\rho:\ps\to\mathrm{PU}(1,1)$
is maximal if and only if it is the holonomy representation of a complete
hyperbolic metric on the interior $\Sigma^\circ$ of $\Sigma$.
\end{thm_intro}

\begin{rem_intro} Theorem~\ref{thm_intro:thm2} is proved in \S~\ref{subsec:3.2};
the ``only if'' part is a consequence of \cite{Burger_Iozzi_difff}
together with the formula in Theorem~\ref{thm:3.2}; the ``if'' part,
while being Gauss--Bonnet's theorem in the boundaryless case,
requires a quite different argument when $\bs\neq\emptyset$.

We observe that when $\bs\neq\emptyset$, the invariant $\tsr$ depends not only on 
$\ps$ but also on $\Sigma$: in fact, if $\Sigma_1,\Sigma_2$ are
nondiffeomorphic surfaces with isomorphic fundamental groups,
$i:\pi_1(\Sigma_1)\to\pi_1(\Sigma_2)$ is an isomorphism 
and $\rho:\pi_1(\Sigma_2)\to\mathrm{PU}(1,1)$ is a maximal representation,
then it follows from Theorem~\ref{thm_intro:thm2} that $\rho\circ i$ is not maximal.
\end{rem_intro}

The first result beyond the case $\mathrm{PU}(1,1)$ was obtained by Toledo \cite{Toledo_89}
who, in the boundaryless case, showed that a maximal representation 
into $\mathrm{PU}(1,m)$ stabilizes a complex geodesic.  
It turns out that the appropriate generalization of complex geodesic is,
in this context, the notion of {\it maximal tube type subdomain}.
Roughly speaking, {\it tube type domains} are bounded symmetric domains 
which admit a model which corresponds to the upper half space model in the case
of the Poincar\'e disk, and their significance for rigidity questions of isometric group actions
already appeared in \cite{Burger_Iozzi_supq, Burger_Iozzi_Wienhard_kahler}.  
For a general Hermitian symmetric space $\Xx$,
maximal tube type subdomains exist, are of rank equal to the rank of $\Xx$
and are $G$-conjugate.  As alluded to above, 
the maximal tube type subdomains in complex hyperbolic $n$-space are the complex geodesics.

The main structure theorem for maximal representations is:

\begin{thm_intro}\label{thm_intro:thm3} Let $\gG$ be a connected
semisimple algebraic group defined over $\RR$ 
such that $G=\gG(\RR)^\circ$ is of Hermitian type.
Let $\Sigma$ be a compact connected oriented surface 
with (possibly empty) boundary and $\chi(\Sigma)\leq-1$.
If $\rho:\ps\to G$ is a maximal representation, then 
\be
\item $\rho$ is injective with discrete image;
\item the Zariski closure $\hH<\gG$ of
the image of $\rho$ is reductive;
\item the reductive Lie group $H:=\hH(\RR)^\circ$ has compact centralizer in $G$
and the symmetric space $\Yy$ associated to $H$ is Hermitian of tube type;
\item $\rho\big(\ps\big)$ stabilizes a maximal tube type subdomain $\Tt\subset\Xx$.
\ee
\end{thm_intro}

\noindent
{\small{[Theorem~\ref{thm_intro:thm3} is proved in \S~\ref{sec:structure1} 
when $\rho$ has Zariski dense image, 
while the general case is treated in \S~\ref{sec:structure2}.]}}

\begin{rem_intro}\label{rem_intro:5}
\be
\item In the case in which $\bs=\emptyset$,
Theorem~\ref{thm_intro:thm3} was announced in \cite{Burger_Iozzi_Wienhard_ann}.
\item When $\partial\Sigma=\emptyset$, Theorem~\ref{thm_intro:thm3}(4)
was obtained by Hern\'andez for $G=\PU(2,m)$ \cite{Hernandez}.
Assuming that the representation is reductive, 
Bradlow, Garc\'{\i}a-Prada and Gothen also obtained Theorem~\ref{thm_intro:thm3}(4)
for $G=\SU(n,m)$ \cite{Bradlow_GarciaPrada_Gothen}
and for $\mathrm{SO}^\ast(2n)$ \cite{Bradlow_GarciaPrada_Gothen_cacanerveuse}.
In each of these works the Toledo invariant appears as the first Chern class of an appropriate
complex line bundle over $\Sigma$.
\item When $\partial\Sigma\neq\emptyset$ and $G=\PU(1,m)$,
Koziarz and Maubon introduced \cite{Koziarz_Maubon} an invariant lying
in the de Rham cohomology of $\Sigma$ with compact support, whose evaluation
on $[\Sigma,\bs]$ can be shown to be equal to our notion of Toledo invariant;
in this context, they obtained in \cite{Koziarz_Maubon} Theorem~\ref{thm_intro:thm3}(4)
as well as Theorem~\ref{thm_intro:thm2}.
\ee
\end{rem_intro}

The symmetric space $\Yy$ in Theorem~\ref{thm_intro:thm3} is the variety
of maximal compact subgroups of $H$; since $H$ has compact centralizer in $G$,
there is a unique totally geodesic embedding $i:\Yy\to\Xx$
which is not necessarily holomorphic but is {\it tight}.
This latter notion, which is analytic in nature, stems from our approach via 
bounded cohomology, see \cite{Burger_Iozzi_Wienhard_tight}.

A special case of Theorem~\ref{thm_intro:thm3} is when the homomorphism $\rho$
has Zariski dense image.  Then $\Yy=\Xx$ and hence
$\Xx$ is of tube type.  This result is optimal in the sense 
that every tube type domain admits a maximal representation with Zariski dense image.
More precisely, let $d:\DD\to\Xx$ be a {\it diagonal disk}, 
also called {\it tight holomorphic disk} in 
\cite{Clerc_Orsted_2, Burger_Iozzi_Wienhard_tight}
(see (\ref{sec:prelim}.b)\S~\ref{subsubsec:2.1.2} for the definition), 
and $\Delta:\mathrm{SU}(1,1)\to\gG(\RR)^\circ$
a homomorphism associated to $d$.

\begin{thm_intro}\label{thm_intro:thm4} Assume that $\Xx$ is of tube type, 
that $\chi(\Sigma)\leq-2$, and let 
\bqn
h:\ps\to\mathrm{SU}(1,1)
\eqn
be a complete hyperbolization of $\Sigma^\circ$.  
If the surface is of type $(g,n)=(1,2)$ or $(0,4)$, we assume that $h$ sends one,
respectively two, boundary components of $\bs$ to hyperbolic elements.
Then $\rho_0:=\Delta\circ h:\ps\to G$ admits a deformation $(\rho_t)_{t\geq0}$ 
such that:
\be
\item $\rho_t$ is maximal for all $t\geq0$, and
\item $\rho_t$ has Zariski dense image for all $t>0$.
\ee
\end{thm_intro}

\noindent
{\small{[This theorem is proved in \S~\ref{sec:examples}.]}}

\subsection{Boundary maps}\label{subsec:bdrymaps}
Maximal representations give rise to boundary maps with special regularity properties
which in turn play an important role in the study of the set of maximal representations
and in the construction of new invariants thereof.
{\it Monotonicity (or positivity)} is one of these properties 
and in order to express it we need the notion of maximal triples of points 
in the Shilov boundary of a symmetric domain: those are the vertices of ideal
geodesic triangles of maximal K\"ahler area in a sense made precise
in \cite{Clerc_Orsted_2} (see \S~\2.1.3 for the definition).

\begin{thm_intro}\label{thm_intro:thm5} Let $h:\ps\to\mathrm{PU}(1,1)$
be a complete hyperbolization of $\Sigma^\circ$ of finite area and $\rho:\ps\to G$ 
a representation into a group of Hermitian type.  
Then $\rho$ is maximal if and only if there exists a left continuous map
\bqn
\varphi:\partial\DD\to\cs
\eqn
with values in the Shilov boundary $\cs$ of the bounded symmetric domain
associated to $G$ such that 
\be
\item $\varphi$ is strictly $\rho\circ h^{-1}$-equivariant, and
\item $\varphi$ is monotone, that is it maps positively oriented triples on $\partial\DD$ 
to maximal triples on $\cs$.
\ee
\end{thm_intro}

\noindent
{\small{[Theorem~\ref{thm_intro:thm5} is proved in \S~\ref{sec:regularity} in the case in which
$\rho$ has Zariski dense image and in \S~\ref{sec:structure2} in the general case.]}}

\begin{rem_intro} The theorem holds true also if ``left continuous'' is replaced
by ``right continuous''.
\end{rem_intro}

The characterization in Theorem~\ref{thm_intro:thm5} clarifies the
relation between maximal representations and the Hitchin or positive
representations into split real Lie groups which were 
recently studied by Labourie \cite{Labourie_anosov}, Guichard
\cite{Guichard_hyperconvex}, and Fock and Goncharov
\cite{Fock_Goncharov_local}. Indeed in the latter the authors
established a similar characterization in terms of equivariant maps
from $\partial \DD$ into (full) flag varieties which send positively
oriented triple in $\partial\DD$ to positive triples of flags in the
sense of Lusztig.

In the only case when $G$ is of Hermitian type as well as real split,
namely when $G$ is locally isomorphic to a symplectic group $\Sp(V)$,
the Shilov boundary can be identified with the space of Lagrangian
subspaces in $V$ and in this case the notion of maximality of triples
in $\cs$ coincides with the notion of positivity of triples in partial
flag varieties -- such as the space of Lagrangians -- defined by
Lusztig in \cite{Lusztig}.  Thus Theorem~\ref{thm_intro:thm5} implies
that the space of positive representations into ${\rm PSp}(V)$
defined by Fock and Goncharov
\cite[Definition~1.10]{Fock_Goncharov_local}, is a proper subset of
the space of maximal representations. For the Hitchin component (when
$\partial \Sigma = \emptyset$) this was observed in
\cite{Burger_Iozzi_Labourie_Wienhard}.

\bigskip
The issue of continuity of the boundary map $\varphi$ presents itself naturally. 
When $\bs=\emptyset$, the continuity of $\varphi$ was established in 
\cite{Burger_Iozzi_Labourie_Wienhard} in the case in which $G=\Sp(V)$ is a symplectic group,
as a byproduct of the construction of an Anosov system;
the case of a general group of Hermitian type will be treated in a forthcoming paper.
When $\bs\neq\emptyset$, then already in the case $G=\mathrm{PU}(1,1)$
the map $\varphi$ will not be in general continuous as 
the case in which $\rho$ is an infinite area hyperbolization indicates.
In fact, if $G=\mathrm{PU}(1,1)$, the map $\varphi$ is a semiconjugacy 
in the sense of Ghys \cite{Ghys_87} and this will be used in \S~\ref{subsec:8.2}
to define a canonical integral bounded class
\bqn
\kszb\in\hb^2\big(\ps,\ZZ\big)
\eqn
which, when $\bs=\emptyset$, corresponds to the fundamental class under the comparison map.
Letting $\ksb\in\hb^2\big(\ps,\RR\big)$ denote the corresponding real class, 
we will establish in \S~\ref{subsec:8.2} (see Corollary~\ref{cor:8.6}) the following

\begin{cor_intro}\label{cor_intro:cor6} For any homomorphism 
$\rho:\ps\to G$, the following are equivalent:
\be
\item $\rho$ is maximal, and
\item $\rho^\ast(\kappa)=\lambda_G(\kappa)\ksb$, for all $\kappa\in\hcb^2(G,\RR)$,
where $\lambda_G$ is a certain explicit linear form on $\hcb^2(G,\RR)$.
\ee
\end{cor_intro}

The extent to which Corollary~\ref{cor_intro:cor6} does not hold for integral classes
will be the source of new invariants of maximal representations 
which will be given an explicit form once rotation numbers are introduced
in the next section.

\subsection{Toledo invariant and rotation numbers}\label{subsec:tolinvrot}
In order to define a notion of integral class in continuous bounded cohomology
we consider, for $G$ a locally compact second countable group 
and $A=\ZZ$ or $\RR$, the cohomology $\hhcb^\bullet(G,A)$ 
of the complex of bounded Borel cochains on $G$ which turns out to coincide
with bounded continuous cohomology if $A=\RR$ (see \S~\ref{subsec:2.3}).
Given $\kappa\in\hhcb^2(G,\ZZ)$, we introduce in \S~\ref{sec:rotation} the rotation number 
\bqn
\rot_\kappa:G\to\RR/\ZZ\,,
\eqn
which is a class function whose restriction to any amenable closed subgroup is a homomorphism
and we show its continuity (Corollary~\ref{cor:7.5}).  

The rotation number $\rot_\kappa$ generalizes the classical rotation
number of an orientation preserving homeomorphism of the circle as
well as the symplectic rotation number introduced by Barge and Ghys in
\cite{Barge_Ghys} and the construction of Clerc and Koufany in
\cite{Clerc_Koufany}.  The exact relations are discussed in
\S~\ref{sec:rotation}.

When $G$ if of Hermitian type and $K<G$ is a maximal compact subgroup, 
the basic properties of the rotation number $\rotk$ are summarized in the following

\begin{thm_intro}\label{thm_intro:thm7}  
\be
\item The map 
\bqn
\ba
\hhcb^2(G,\ZZ)&\to\homc(K,\RR/\ZZ)\\
\kappa\quad\,\,&\mapsto\quad\rot_\kappa|_K
\ea
\eqn
is an isomorphism.
\item The change of coefficients $\hhcb^2(G,\ZZ)\to\hcb^2(G,\RR)$ is injective 
with image a lattice.
\item For every $g\in G$,
\bqn
\rot_\kappa(g)=\rot_\kappa(k)\,,
\eqn
where $k\in K$ is conjugate to the elliptic component $g_e$ in the
refined Jordan decomposition $g=g_eg_hg_u$ of $g$.
\item The unique continuous lift
\bqn
\widetilde\rot_\kappa:\widetilde G\to\RR
\eqn
vanishing at $e$ is a homogeneous quasimorphism.
\ee
\end{thm_intro}

\noindent
{\small{[Theorem~\ref{thm_intro:thm7} is proved in Propositions~\ref{prop:7.6}, \ref{prop:7.7} 
and Theorem~\ref{thm:7.8}; for the refined Jordan decomposition see \cite[\S~2]{Borel_comm}.]}}

We turn now to the formula of the Toledo invariant $\tksr$ 
when $\kappa$ is an integral bounded class.
For this we assume that $\bs\neq\emptyset$ and let 
\bq\label{eq:ps}
 \quad\quad\pi_1(\Sigma)
=\bigg\<a_1,b_1,\dots,a_g,b_g,c_1,\dots,c_n:\,
\prod_{i=1}^g[a_i,b_i]\prod_{j=1}^nc_j=e\bigg>
\eq
be a presentation where the elements $c_i$ represent loops which are freely homotopic
to the corresponding boundary components of $\bs$ with positive orientation.
Given a homomorphism $\rho:\ps\to G$, let $\widetilde\rho:\ps\to \widetilde G$
be a lift of $\rho$ to the universal covering $\widetilde G$, 
taking into account that $\ps$ is free. 

\begin{thm_intro}\label{thm_intro:thm8}  Let $\kappa\in\hhcb^2(G,\ZZ)$.
Then
\bqn
 \mathrm{T}_\k(\Sigma,\rho)
=-\sum_{j=1}^n\widetilde\rotk\big(\widetilde\rho(c_j)\big)\,.
\eqn
\end{thm_intro}

\noindent
{\small{[Theorem~\ref{thm_intro:thm8} is proved in \S~\ref{subsec:8.1}.]}}

When the boundary of $\Sigma$ is empty, a formula for $\mathrm{T}_\kappa$ (see Theorem~\ref{thm:8.4}) can be obtained 
by cutting $\Sigma$ along a separating loop and using Theorem~\ref{thm_intro:thm8}
together with the additivity property of the Toledo invariant in Theorem~\ref{thm_intro:thm1}(3).

\bigskip
In conjunction with Corollary~\ref{cor_intro:cor6}, rotation  numbers give rise
to nontrivial invariants of maximal representations;
recalling that the Shilov boundary $\cs$ of the bounded symmetric domain $\Dd$
associated to $G$ is a homogeneous space with typical stabilizer $Q$ 
and letting $e_G$ denote the exponent of the finite group $Q/Q^\circ$, we have:

\begin{thm_intro}\label{thm_intro:thm9}  Let $\kappa\in\hhcb^2(G,\ZZ)$ and $\rho_0: \ps \to G$ a maximal representation.
\be
\item For every maximal representation $\rho:\ps\to G$ the map
\bqn
\ba
\mathrm{R}^{\rho_0}_\kappa(\rho):\ps&\longrightarrow\qquad\qquad\RR/\ZZ\\
\gamma\quad&\mapsto\rot_\kappa\big(\rho(\gamma)\big) - \rot_\kappa\big(\rho_0(\gamma)\big)
\ea
\eqn
is a homomorphism.
\item If $\Dd$ is of tube type, then $\mathrm{R}^{\rho_0}_\kappa(\rho)$ takes
values in $e_G^{-1}\ZZ/\ZZ$ and
\bqn
\hmaxpg\to\hom\big(\ps,\RR/\ZZ\big)
\eqn
is constant on connected components.
\ee
\end{thm_intro}

\noindent
{\small{[Theorem~\ref{thm_intro:thm9} is proved in \S~\ref{subsec:8.2}.]}}

In Example~\ref{exo:8.8} we describe for $G=\Sp(V)$, $\dim(V)=4m$ and $\bs=\emptyset$ that already 
$\rot_\kappa(\rho):\ps \to \ZZ/2\ZZ$ is a homomorphism and discuss  
its relationship with the first
Stiefel--Withney class of a certain real vector bundle constructed
using the boundary map from Theorem~\ref{thm_intro:thm5}.

\bigskip
We turn now to our final application to representation varieties.  
For this we assume again that $\bs\neq\emptyset$ and use the familiar presentation 
of $\ps$ given in \eqref{eq:ps}.  Then

\bqn
 \hom^{\cs}\big(\ps,G\big)
:=\big\{\rho\in\hom\big(\ps,G\big):\,
  \rho(c_i)\hbox{ has at least}\\ \hbox{ one fixed point in }\cs,\,1\leq i\leq n\big\}\,,
\eqn
is a semialgebraic set if $G$ is real algebraic and we have as a consequence 
of Theorem~\ref{thm_intro:thm5}, that 
\bq\label{eq:max<cs}
 \hom_\mathrm{max}\big(\pi_1(\Sigma),G\big)
\subset\hom^{\cs}\big(\pi_1(\Sigma),G\big)\,.
\eq

\begin{cor_intro}\label{cor_intro:cor10}  Let $\kappa\in\hhcb^2(G,\ZZ)$
and assume that $\Dd$ is of tube type.  Then:
\be
\item $\tksr\in e_G^{-1}\ZZ$ for every $\rho\in \hom^{\cs}\big(\ps,G\big)$, and
\item $\hmaxpg$ is a union of connected components of the set $\hom^{\cs}\big(\ps,G\big)$.
\ee
\end{cor_intro}

\noindent
{\small{[Corollary~\ref{cor_intro:cor10} is proved in \S~\ref{subsec:8.3}.]}}

An alternative boundary condition might be imposed by fixing instead a
set $\Cc=\{\Cc_1,\dots,\Cc_n\}$ of conjugacy classes in $G$ and
defining \bqn \hom^\Cc\big(\pi_1(\Sigma),G\big)
:=\big\{\rho\in\hom\big(\pi_1(\Sigma),G\big):\,
\rho(c_i)\in\Cc_i,\,1\leq i\leq n\big\}\,. \eqn Then
$\hom^\Cc\big(\pi_1(\Sigma),G\big)$ is also a semialgebraic set and it
follows immediately from Theorem~\ref{thm_intro:thm8} that $T_\kappa$
is constant on its connected components.

Notice however that Corollary~\ref{cor_intro:cor10}
implies that for many choices of conjugacy classes the intersection
$\hom_\mathrm{max}\big(\pi_1(\Sigma),G\big)\cap
\hom^\Cc\big(\pi_1(\Sigma),G\big)$ considered above
is actually empty.  For example in the case
when $\Sigma$ has precisely one boundary component
Theorem~\ref{thm_intro:thm8} readily implies that for any maximal
representation the rotation number of the conjugacy class
$\Cc=\{\Cc_1\}$ has to be zero.  From a different point of view,
fixing a conjugacy class $\Cc$ with nonzero rotation number, gives a
modified Milnor--Wood type inequality as in
Theorem~\ref{thm_intro:thm1}(1) for the Toledo invariant restricted to
$\hom^{\Cc}\big(\pi_1(\Sigma),G\big)$.  Goldman showed in
\cite{Goldman_torus} that for the one-punctured torus, representations
into ${\rm PSL}(2,\RR)$ maximizing the Toledo invariant with respect
to this modified Milnor--Wood type inequality correspond to singular
hyperbolic structures on the torus with cone type singularities in the
puncture.

\bigskip
\bigskip
\noindent
{\it Acknowledgments:} The authors are grateful to D.~Toledo for his
crucial comments in the beginning of our study of maximal
representations and for his continuing interest and support.  The
authors thank N.~A'Campo for asking the question which triggered the
study of maximal representations of surface groups with boundary,
Y.~Benoist for suggesting which formula computes the Toledo invariant
(see Theorem~\ref{thm_intro:thm8}) and F.~Labourie for various
enjoyable and useful conversations.  Our thanks go also to S.~Bradlow,
O.~Garc\'{\i}a-Prada, P.~Gothen and I.~Mundet i Riera for interesting
discussions concerning the relation between our work and their
approach through Higgs bundles. The third author wants to thank
O.~Guichard for his questions and remarks, which lead to a correction
in Theorem~\ref{thm_intro:thm9}.
%and, more specifically, about \S~\ref{subsubsec:8.2.1}.

%%% Local Variables: 
%%% mode: latex
%%% TeX-master: "toledo"
%%% End: 

%% file: prelim.tex
\section{Preliminaries}\label{sec:prelim}

\subsection{Hermitian symmetric spaces, bounded continuous cohomology}\label{subsec:2.1}
\subsubsection{\ }\label{subsubsec:2.1.1} A Lie group $G$ is of {\it Hermitian type} 
if it is connected, semisimple with finite center and no compact factors,
and if the associated symmetric space is Hermitian. A Lie group $G$ is of type (RH)
if it is connected reductive with compact center and the quotient $G/G_c$ by the largest
connected compact normal subgroup $G_c$ is of Hermitian type.

If $G$ is a locally compact group, $\hc^\bullet(G,\RR)$ denotes the continuous cohomology
with $\RR$-trivial coefficients, while $\hcb^\bullet(G,\RR)$ is the bounded continuous cohomology;
for the general theory concerning the latter and its relation to the former, we
refer to \cite{Monod_book, Burger_Monod_GAFA, Burger_Iozzi_app, Burger_Iozzi_formula}.

When $G$ is of type (RH) and $\Xx$ is its associated symmetric space, we have isomorphisms
\bqn
\xymatrix@1{
 \Omega^2(\Xx)^G\ar[r]
&\hc^2(G,\RR)
&\hcb^2(G,\RR)\ar[l]\,,
}
\eqn
where the first is the Van Est isomorphism between the complex $\Omega^\bullet(\Xx)^G$ 
of $G$-invariant differential forms on $\Xx$ and $\hc^\bullet(G,\RR)$ \cite{Van_Est},
while the second is the comparison map which in degree two is an isomorphism 
\cite{Burger_Monod_GAFA}.
Given $\omega\in\Omega^2(\Xx)^G$ and $x\in\Xx$ a basepoint,
the function
\bqn
c_\omega(g_0,g_1,g_2):=\frac{1}{2\pi}\int_{\Delta(g_0x,g_1x,g_2x)}\omega\,,
\eqn
where $\Delta(g_0x,g_1x,g_2x)$ denotes a smooth triangle with geodesic sides,
defines a homogeneous $G$-invariant cocycle which is moreover bounded;
when $\omega=\omega_\Xx$ is the K\"ahler form for the unique $G$-invariant Hermitian metric 
of minimal holomorphic sectional curvature $-1$ ({\em normalized metric}),
we let $\kg\in\hc^2(G,\RR)$ and $\kgb\in\hcb^2(G,\RR)$ 
denote the corresponding classes and refer to $\kg$ (respectively $\kgb$)
as the K\"ahler class (respectively bounded K\"ahler class).
For the Gromov norm of $\kgb$, we have that 
\bq\label{eq:gromov-norm}
\|\kgb\|=\frac{\rk_\Xx}{2}\,,
\eq
where $\rk_\Xx$ is the rank of $\Xx$.

\subsubsection{\ }\label{subsubsec:2.1.2}  Let $G$ be of Hermitian type and 
$\Xx$ be the associated symmetric space.  Then:
\be
\item[(\ref{sec:prelim}.a)] A {\it maximal polydisk} in $\Xx$ is the image
  of a totally geodesic and holomorphic embedding
  $t:\DD^{\rk_\Xx}\to\Xx$ of a product of $\rk_\Xx$ Poincar\'e disks.
  
\item[(\ref{sec:prelim}.b)] A {\it diagonal disk} (or {\it tight
    holomorphic disk} in $\Xx$ is the image of the diagonal
  $\DD\subset\DD^{\rk_\Xx}$ under $t$; we will denote by $d:\DD\to\Xx$
  the resulting totally geodesic and holomorphic embedding.  \ee

To the above objects are associated a connected finite covering $L$ of $\mathrm{PU}(1,1)$
and homomorphisms $\tau:L^{\rk_\Xx}\to G$ and $\Delta:L\to G$ with respect to which 
$t$ and $d$ are equivariant.
We have moreover that
\bq\label{eq:disks}
\tau^\ast(\kgb)= \kappa_{L^{\rk_\Xx}}^\mathrm{b}\qquad\text{ and } \Delta^\ast(\kgb)=\rk_\Xx\klb\,.
\eq

\subsubsection{\ }\label{subsubsec:2.1.3}
Let $G$ be of type (RH), $\Xx$ the associated symmetric space, 
$\Dd$ the bounded domain realization and $\cs$ its Shilov boundary.  
Then $\cs$ is a homogeneous $G$-space of the form $G/Q$,
where $Q$ is a specific parabolic subgroup (which is maximal if $\Xx$ is irreducible). Two points $x,y\in\cs$
are {\it transversal} if $(x,y)$ lies in the open $G$-orbit in $\cs^2$.
Let $\cst$ denote the set of triples of pairwise transversal points.
It was shown by Clerc and \O rsted \cite{Clerc_Orsted_TG, Clerc_Orsted_2} that the map
\bqn
\ba
\Dd^3\,\,&\longrightarrow\qquad\RR\\
(x,y,z)&\mapsto\frac{1}{2\pi}\int_{\Delta(x,y,z)}\omega_\Dd\,, 
\ea
\eqn
where $\omega_\Dd$ is the K\"ahler form for the normalized metric 
on $\Dd$, extends continuously to $\cst$;
Clerc showed then that by taking appropriate tangential limits
one obtains a well defined $G$-invariant Borel cocycle
\bqn
\beta_{\cs}:\cs^3\to\RR
\eqn
extending the previous one, which satisfies
\bqn
|\beta_{\cs}(x,y,z)|\leq\frac{\rk_\Dd}{2}\,,
\eqn
where $\rk_\Dd = \rk_\Xx$ (see \cite[Theorem~5.3]{Clerc_maslov}).
In the following $\b_\cs$ will be referred to as the {\em generalized Maslov cocycle}
and a triple $x,y,z\in\cs$ for which $\b_\cs(x,y,z)=\frac{\rk_\Dd}{2}$ will
be called {\it maximal}.

Let now $\big(\balt(\cs^\bullet)\big)$ denote the complex of bounded alternating
Borel cocycles on $\cs$.  Then, under the canonical map
\bq\label{eq:can}
\h^\bullet\big(\balt(\cs^\bullet)^G\big)\to\hcb^\bullet(G,\RR)
\eq
(see \cite[\S~4.2.]{Burger_Iozzi_Wienhard_kahler}, \cite[Corollary~2.2.]{Burger_Iozzi_app})
the class defined by $\beta_{\cs}$ corresponds to $\kgb$.

\subsection{Bounded singular and bounded group cohomology}\label{subsec:2.2}
A ``space'' will always refer to a countable CW-complex 
and $A$ will be one of the coefficients $\ZZ,\RR,\RR/\ZZ$.  
For a pair of spaces $Y\subset X$,
$\h^\bullet(X,Y,A)$ and $\hb^\bullet(X,Y,A)$ denote respectively 
the singular relative cohomology with coefficients in $A$ 
and its bounded counterpart;
observe that $\h^\bullet(X,Y,\RR/\ZZ)=\hb^\bullet(X,Y,\RR/\ZZ)$.
Also, $\h^\bullet\big(\pi_1(X),A\big)$ and $\hb^\bullet\big(\pi_1(X),A\big)$ denote
respectively the group cohomology and the bounded group cohomology 
of $\pi_1(X)$ with $A$-coefficients and 
$\h^\bullet\big(\pi_1(X),\RR/\ZZ\big)=\hb^\bullet\big(\pi_1(X),\RR/\ZZ\big)$.
These cohomology theories come with the following natural comparison maps
\bqn
\ba
&\hb^\bullet(X,Y,A)\to\h^\bullet(X,Y,A)\\
&\hb^\bullet\big(\pi_1(X),A\big)\to\h^\bullet\big(\pi_1(X),A\big)\\
&\h^\bullet\big(\pi_1(X),A\big)\to\h^\bullet\big(X,A\big)\\
&\hb^\bullet\big(\pi_1(X),A\big)\to\hb^\bullet\big(X,A\big)\,,
\ea
\eqn
where the last two are induced by the classifying map $X\to B\pi_1(X)$.

We recall the following facts:
\be
\item[(\ref{sec:prelim}.c)] the short exact sequence
\bq\label{eq:coeff}
\xymatrix@1{
 0\ar[r]
&\ZZ\ar[r]
&\RR\ar[r]
&\RR/\ZZ\ar[r]
&0
}
\eq
gives rise to long exact sequences in each of the four cohomology theories;
these sequences are natural with respect to the four comparison maps;
\item[(\ref{sec:prelim}.d)] the inclusion of spaces $Z_1\subset Z_2\subset X$
induces a long exact sequence in singular relative and bounded singular
relative cohomology which fit into the long exact sequences coming from 
the coefficient sequence in (\ref{sec:prelim}.c);
\item[(\ref{sec:prelim}.e)] in general the comparison map
\bqn
g_X:\hb^\bullet\big(\pi_1(X),\RR\big)\to\hb^\bullet(X,\RR)
\eqn
is an isomorphism \cite{Gromov_82, Brooks, Ivanov}, referred to as {\em Gromov isomorphism};
as a consequence, if each connected component of $Z_1$ and $Z_2$ 
has amenable fundamental group the map
\bqn
j_{Z_1,Z_2}:\hb^\bullet(X,Z_2,\RR)\to\hb^\bullet(X,Z_1,\RR)
\eqn
is an isomorphism; when $Z_1=\emptyset$ we set $j_{\emptyset,Z_2}=:j_{Z_2}$ for ease of notation.  
We will only need the above isomorphism when all spaces involved are $K(\pi,1)$'s, 
in which case we have for all coefficients $A$ a commutative diagram
\bqn
\xymatrix{
 \hb^\bullet\big(\pi_1(X),A\big)\ar[r]^-{g_X}\ar[d]
&\hb^\bullet(X,A)\ar[d]\\
 \h^\bullet\big(\pi_1(X),A\big)\ar[r]^-{g_X}
&\h^\bullet(X,A)\,,
}
\eqn
where the horizontal maps are isomorphisms.
\ee

\subsection{(Bounded) Borel cohomology versus (bounded) continuous cohomology}\label{subsec:2.3}
Given a locally compact group $G$ and $A=\ZZ,\RR,\RR/\ZZ$, we have the complexes
$\big(\mathrm{C}(G^\bullet,A)\big)$, $\big(\cb(G^\bullet,A)\big)$, 
$\big(\mathcal{B}(G^\bullet,A)\big)$ and $\big(\mathcal{B}_\mathrm{b}(G^\bullet,A)\big)$, 
of $A$-valued continuous, bounded continuous, Borel and bounded Borel cochains on $G$
which lead, by taking the cohomology of the $G$-invariants,
to the $A$-valued continuous $\hc^\bullet(G,A)$, bounded continuous $\hcb^\bullet(G,A)$,
Borel $\hhc^\bullet(G,A)$ and bounded Borel $\hhcb^\bullet(G,A)$ cohomology.
Of course when $A=\ZZ$ the first two cohomology theories are not of much use
and their Borel version is a natural substitute.
We have at any rate comparison maps coming from the obvious inclusions of complexes
\bqn
\xymatrix{
 \hcb^\bullet(G,A)\ar[r]\ar[d]
&\hc^\bullet(G,A)\ar[d]\\
 \hhcb^\bullet(G,A)\ar[r]
&\hhc^\bullet(G,A)\,.
}
\eqn
For us the following facts will be of importance:
\be
\item[(\ref{sec:prelim}.f)] the short exact sequence 
\bqn
\xymatrix@1{
 0\ar[r]
&\ZZ\ar[r]
&\RR\ar[r]
&\RR/\ZZ\ar[r]
&0
}
\eqn
gives rise to long exact sequences in Borel and bounded Borel cohomology 
which are compatible with respect to the comparison map;
\item[(\ref{sec:prelim}.g)] $\hcb^\bullet(G,\RR/\ZZ)=\hc^\bullet(G,\RR/\ZZ)$
and $\hhcb^\bullet(G,\RR/\ZZ)=\hhc^\bullet(G,\RR/\ZZ)$;
\item[(\ref{sec:prelim}.h)] if $A=\RR,\ZZ$, then $\hhc^1(G,A)=\homc(G,A)$
and $\hcb^1(G,A)=0$;
\item[(\ref{sec:prelim}.i)] $\hcb^\bullet(G,\RR)=\hhcb^\bullet(G,\RR)$, which
can be checked by using the regularization operators defined in \cite[\S~4]{Blanc};
\item[(\ref{sec:prelim}.l)] if $G$ is a Lie group, the comparison map
$\hc^\bullet(G,\RR)\to\hhc^\bullet(G,\RR)$ is an isomorphism \cite[Theorem~3]{Wigner_73}.
\ee

%%% Local Variables: 
%%% mode: latex
%%% TeX-master: "toledo"
%%% End: 

%% file: tolno.tex
\section{Toledo numbers, basic properties and first consequences}\label{sec:tolno}

\subsection{Definitions and basic properties}\label{subsec:3.1}
Let $\Sigma$ be a compact oriented surface with (possibly empty) boundary $\partial\Sigma$,
$G$ a locally compact group and $\rho:\ps\to G$ a homomorphism.
Using the diagram
\bqn
\xymatrix{
 \hcb^2(G,\RR)\ar[r]^{\rho^\ast}
&\hbps\ar[r]\ar[r]^{g_\Sigma}
&\hb^2(\Sigma,\RR)\\
& &\hb^2(\Sigma,\bs,\RR)\ar[u]_{j_{\partial \Sigma}}
}
\eqn
where $\rho^\ast$ is the pullback in bounded cohomology,
$g_\Sigma$ the Gromov isomorphism in (\ref{sec:prelim}.e) and
$j_\bs$ is the isomorphism in (\ref{sec:prelim}.e) in bounded singular cohomology
induced by the inclusion $(\Sigma,\emptyset)\to(\Sigma,\bs)$,
we make the following 
\begin{defi}
%define the {\it Toledo number of $\rho$ relative to a class $\kappa\in\hcb^2(G,\RR)$} by
The {\it Toledo number of $\rho$ relative to a class $\kappa\in\hcb^2(G,\RR)$} is
\bqn
\tksr:=\big\langle(j_\bs)^{-1}g_\Sigma\rho^\ast(\kappa),[\Sigma,\bs]\big\rangle\,.
\eqn
Here $[\Sigma,\bs]\in\h^2(\Sigma,\bs,\RR)$ denotes the relative fundamental class
and $(j_\bs)^{-1}g_\Sigma\rho^\ast(\kappa)$ is considered as an ordinary relative cohomology class.
\end{defi}

When $G$ is of type (RH) the {\it Toledo number $\tsr$ of $\rho$} is defined 
as $\tksr$ where $\kappa=\kgb$ is the bounded K\"ahler class (see \S~\ref{subsubsec:2.1.1}).
The following two properties are immediate:
\be
\item[--] if $\rho_1$ and $\rho_2$ are $G$-conjugate, then 
$\mathrm{T}_\kappa(\Sigma,\rho_1)=\mathrm{T}_\kappa(\Sigma,\rho_2)$;
\item[--] if $f:\Sigma_1\to\Sigma_2$ is a continuous map of degree $d\geq1$,
$\rho_i:\pi_1(\Sigma_i)\to G$ are homomorphisms related by $\rho_1=\rho_2 f_\ast$, 
where $f_\ast$ is the morphism induced on the fundamental groups, then 
$\mathrm{T}_\kappa(\Sigma_1,\rho_1)=d\cdot\mathrm{T}_\kappa(\Sigma_2,\rho_2)$.
\ee
The next results describe the behavior of the Toledo numbers
under natural topological operations on surfaces.

\begin{prop}\label{prop:3.1} Let $\Sigma$ be a surface and $\rho:\ps\to G$ a homomorphism.
\be
\item (Additivity) If $\Sigma=\Sigma_1\cup_C\Sigma_2$ is the connected sum of two 
sub surfaces $\Sigma_i$ along a separating loop $C$, then
\bqn
\mathrm{T}_\kappa(\Sigma,\rho)=\mathrm{T}_\kappa(\Sigma_1,\rho_1)+\mathrm{T}_\kappa(\Sigma_2,\rho_2)\,,
\eqn
where $\rho_i$ is the restriction of $\rho$ to $\pi_1(\Sigma_i)$.
\item (Invariance under gluing)  If $\Sigma'$ is the surface obtained by cutting $\Sigma$ along
a nonseparating loop $C$ and $i:\Sigma'\to\Sigma$ is the canonical map, then
\bqn
\mathrm{T}_\kappa(\Sigma',\rho i_\ast)=\mathrm{T}_\kappa(\Sigma,\rho)\,.
\eqn
\ee
\end{prop}

\begin{proof} Here we prove the additivity property, 
the proof of the invariance under gluing proceeds along similar lines.
Let $\alpha\in\hb^2(\Sigma,\bs)$ and let 
\bqn
j:=j_{\bs\cup C,\bs}:\hb^2(\Sigma,\bs\cup C)\to\hb^2(\Sigma,\bs)
\eqn
be the morphism given by the inclusion $(\Sigma,\bs)\to(\Sigma,\bs\cup C)$;
notice that $j$ is an isomorphism since every connected component of $\bs\cup C$ has amenable fundamental group.
Then:
\bqn
\ba
 \big\langle\alpha,[\Sigma,\bs]\big\rangle
&=\big\langle j^{-1}(\alpha),[\Sigma_1,\bs_1]+[\Sigma_2,\bs_2]\big\rangle\\
&=\big\langle j^{-1}(\alpha)|_{\Sigma_1},[\Sigma_1,\bs_1]\big\rangle
 +\big\langle j^{-1}(\alpha)|_{\Sigma_2},[\Sigma_2,\bs_2]\big\rangle\,.
\ea
\eqn
Using that $j^{-1}(\alpha)|_{\Sigma_i}=(j_{i})^{-1}(\alpha|_{\Sigma_i})$, where $j_i:= j_{\bs_i-C, \bs_i}$ we get
\bqn
 \big\langle\alpha,[\Sigma,\bs]\big\rangle
=\big\langle (j_{1})^{-1}(\alpha|_{\Sigma_1}),[\Sigma_1,\bs_1]\big\rangle
 +\big\langle (j_{2})^{-1}(\alpha|_{\Sigma_2}),[\Sigma_2,\bs_2]\big\rangle\,.
\eqn
Specializing to $j_\bs(\alpha)=g_\Sigma\rho^\ast(\kappa)$ and 
observing that 
\bqn
(j_{i})^{-1}(\alpha|_{\Sigma_i})=(j_{i})^{-1}(\alpha)|_{\Sigma_i}=(j_{\bs_i})^{-1}g_{\Sigma_i}\rho_i^\ast(\kappa)
\eqn
concludes the proof.
%Using that $j^{-1}(\alpha)|_{\Sigma_i}=(j_{\bs_i})^{-1}(\alpha|_{\Sigma_i})$ we get
%\bqn
% \big\langle\alpha,[\Sigma,\bs]\big\rangle
%=\big\langle (j_{\bs_1})^{-1}(\alpha|_{\Sigma_1}),[\Sigma_1,\bs_1]\big\rangle
% +\big\langle (j_{\bs_2})^{-1}(\alpha|_{\Sigma_2}),[\Sigma_2,\bs_2]\big\rangle\,.
%\eqn
%Specializing to $j_\bs(\alpha)=g_\Sigma\rho^\ast(\kappa)$ and 
%observing that 
%\bqn
%(j_{\bs_i})^{-1}(\alpha|_{\Sigma_i})=(j_{\bs_i})^{-1}(\alpha)|_{\Sigma_i}=g_{\Sigma_i}\rho_i^\ast(\kappa)
%\eqn
%concludes the proof.
\end{proof}

\subsection{The analytic formula}\label{subsec:3.2}
In this section we relate the Toledo numbers introduced in 
\S~\ref{subsec:3.1} to invariants introduced and studied in \cite{Burger_Iozzi_difff}. 
Let $G$ be a locally compact group. 
Let $L$ be a finite connected covering of $\mathrm{PU}(1,1)$, 
$\Gamma<L$ a lattice and $\rho:\Gamma\to G$ a homomorphism. 
Composing the transfer map
\bqn
\mathrm{T}_\mathrm{b}:\hb^2(\Gamma,\RR)\to\hcb^2(L,\RR)
\eqn
with the pullback $\rho^\ast$, we obtain the bounded Toledo map
\bqn
\mathrm{T}_\mathrm{b}(\rho):\hcb^2(G,\RR)\to\hcb^2(L,\RR)=\RR\klb
\eqn
(see \cite{Burger_Iozzi_difff}) which leads to an invariant $\tb(\rho,\kappa)\in\RR$
given by 
\bqn
\mathrm{T}_\mathrm{b}(\rho)(\kappa)=\tb(\rho,\kappa)\klb
\eqn
for $\kappa\in\hcb^2(G,\RR)$.  When $G$ is of type (RH) we set, 
in analogy with \S~\ref{subsec:3.1}, 
\bqn
\tb(\rho)=\tb(\rho,\kgb)\,.
\eqn

\begin{thm}\label{thm:3.2} Let $h:\ps\to\Gamma$ be an isomorphism
whose composition with the projection to $\mathrm{PU}(1,1)$ is the developing homomorphism
of a complete hyperbolic structure on $\Sigma^\circ$ with finite area.  Then
\bqn
\T_\kappa(\Sigma,\rho)=|\chi(\Sigma)|\,\tb(\rho\circ h^{-1},\kappa)
\eqn
for any homomorphism $\rho:\ps\to G$ and $\kappa\in\hcb^2(G,\RR)$.
\end{thm}

We defer the proof of this theorem until \S~\ref{subsec:3.3} and collect here a few important consequences.

\begin{cor}\label{cor:3.3}
We have the following Milnor--Wood type bounds: 
\be
\item $|\T_\kappa(\Sigma,\rho)|\leq2|\chi(\Sigma)|\,\|\kappa\|$;
\item if $G$ is of type (RH), then $|\T(\Sigma,\rho)|\leq|\chi(\Sigma)|\,\rk_\Xx$.
\ee
\end{cor}

\begin{proof} The first assertion follows from Theorem~\ref{thm:3.2} and the fact that the transfer
$\T_\mathrm{b}$ and the pullback are both norm decreasing.  The second assertion
follows from the first one and the equality $\|\kgb\|=\frac{\rk_\Xx}{2}$ (see \S~\ref{subsec:2.1}).
\end{proof}

The following are then the two main concepts of this paper:

\begin{defi}\label{def:3.4} Let $G$ be a group of type (RH).
\be
\item A homomorphism $\rho:\ps\to G$ of a surface group $\ps$ is {\it maximal} if $\T(\Sigma,\rho)=|\chi(\Sigma)|\rk_\Xx$.
\item A homomorphism $\rho:\Gamma\to G$ of a lattice $\Gamma<L$ is {\it maximal} if
$\tb(\rho)=\rk_\Xx$.
\ee
\end{defi}

Observe that the first definition generalizes the concept of maximal representation given in the introduction
and puts it in the context of groups of type (RH) which will turn out to be the right one for the proofs.
The second concept of maximality is equivalent to the one introduced in 
\cite{Burger_Iozzi_difff}, the equivalence being given by \cite[Lemma~5.3]{Burger_Iozzi_difff}.
The relationship between the above definitions is given by Theorem~\ref{thm:3.2}.

\begin{proof}[{\it Proof of Theorem~\ref{thm_intro:thm2}}] Let $h:\ps\to\mathrm{PU}(1,1)$
be a hyperbolization of $\Sigma^\circ$ with finite area and image $\Gamma$;
in particular $h$ is induced by a diffeomorphism $f:\Sigma^\circ\to\Gamma\backslash\DD$.
Let $\rho:\ps\to\mathrm{PU}(1,1)$ be a homomorphism.  
If $\rho$ is maximal, then, by Theorem~\ref{thm:3.2}, $\rho\circ h^{-1}:\Gamma\to G$ is maximal 
as a representation of the lattice $\Gamma$ and it follows then from
\cite[Lemma~5.2 and Corollary~11]{Burger_Iozzi_difff} that $\rho\circ h^{-1}$
is induced by a diffeomorphism 
\bqn
f_\rho:\Gamma\backslash\DD\to\rho\big(\ps\big)\backslash\DD
\eqn
which implies that $\rho$ itself is induced by the diffeomorphism 
\bqn
f_\rho\circ f:\Sigma^\circ\to\rho\big(\ps\big)\backslash\DD\,.
\eqn

Conversely, if $\rho$ is induced by a complete hyperbolic metric on $\Sigma^\circ$,
there exists a semiconjugation $F:\partial\DD\to\partial\DD$ in the sense of Ghys
\cite{Ghys_87} with $\rho(\gamma)F=F h(\gamma)$.
Since $\kappa_{\mathrm{PU}(1,1)}^\mathrm{b}$ is the bounded real Euler class 
\cite{Ghys_87}, we have that 
\bqn
\rho^\ast(\kappa_{\mathrm{PU}(1,1)}^\mathrm{b})=h^\ast(\kappa_{\mathrm{PU}(1,1)}^\mathrm{b})
\eqn
and hence
\bqn
(\rho\circ h^{-1})^\ast(\kappa_{\mathrm{PU}(1,1)}^\mathrm{b})=\kappa_{\mathrm{PU}(1,1)}^\mathrm{b}|_\Gamma
\eqn
which, applying the transfer map, implies that $\tb(\rho\circ h^{-1})=1$ 
and thus $\rho$ is maximal by Theorem~\ref{thm:3.2}.
\end{proof}

\subsection{Proof of Theorem~\ref{thm:3.2}}\label{subsec:3.3}
The statement of Theorem~\ref{thm:3.2} can be reformulated as follows.
Let $\Gamma<L$ be a torsionfree lattice; we consider the finite area surface
$S=\Gamma\backslash\DD$ as interior of a compact surface $\overline S$ with boundary $\partial\ol S$,
which is a union of circles.  Given a homomorphism $\rho:\Gamma\to G$ 
and identifying $\Gamma$ with $\pi_1(S)=\pi_1(\ol S)$,
the assertion is that 
\bqn
\T_\kappa(\ol S,\rho)=|\chi(S)|\tb(\kappa,\rho)\,.
\eqn
For $T\geq0$ large enough, let $S_{\geq T}$ denote the union of the convex cusp neighborhoods
bounded by horocycles of length $1/T$.  It is easy to verify that 
if $\beta\in\hb^2(\Gamma,\RR)$ we have 
\bq\label{eq:verif}
 \big\langle(j_{\partial\ol S})^{-1}g_{\ol S}(\beta),[\ol S,\partial\ol S]\big\rangle
=\big\langle(j_T)^{-1}g_S(\beta),[S,S_{\geq T}]\big\rangle
\eq
where, for ease of notation, $j_T$ refers to the canonical isomorphism 
\bqn
\hb^2(S,S_{\geq T},\RR)\to\hb^2(S,\RR)\,.
\eqn
Introducing the notation
\bqn
\Tb(\beta)=\tau(\beta)\klb\,,
\eqn
where $\Tb$ is the transfer operator, the theorem will then follow from \eqref{eq:verif}
and the proposition below applied to $\beta=\rho^\ast(\kappa)$.

\begin{prop}\label{prop:3.5} With the above notation
\bqn
\big\langle(j_T)^{-1}g_S(\beta),[S,S_{\geq T}]\big\rangle=\tau(\beta)|\chi(S)|\,.
\eqn
\end{prop}

The rest of this subsection is devoted to the proof of Proposition~\ref{prop:3.5}
for which we will need the three lemmas below.
We fix the following notation, if $Y$ is any topological space, let $S_m(Y)$ denote the set of singular $m$-simplices 
and $F_\mathrm{b}(Y,\RR)$ the space of bounded $m$-cochains.

\begin{lemma}[{Loeh-Strohm, \cite[Theorem~2.37]{Strohm}}]\label{lem:3.6} 
Let $U\subset\DD$ be a convex subset and 
$\Lambda< L$ be a discrete torsionfree subgroup preserving $U$. 
The canonical isomorphism 
\bqn
\xymatrix@1{
 \hb^m(\Lambda,\RR)\ar[r]^-\cong
&\hb^m(\Lambda\backslash U,\RR)
}
\eqn
can be implemented by the map
\bqn
\ba
\cba(U^{m+1},\RR)^\Lambda&\to F_\mathrm{b}^m(\Lambda\backslash U,\RR)\\
f\quad\quad&\longmapsto\qquad\quad\overline f\,,
\ea
\eqn
defined by $\overline f(\sigma):=f(\tilde\sigma_0,\dots,\tilde\sigma_m)$,
where $\sigma:\Delta^m\to\Lambda\backslash U$ is an $m$-simplex
and $\tilde\sigma:\Delta^m\to U$ is a lift with vertices
$\tilde\sigma_0,\dots,\tilde\sigma_m$.
\end{lemma}

The next lemma follows from standard properties of the transfer map and will also be useful later on.
Let $A(x,y,z)$ denote the area of a geodesic triangle in $\DD$ with vertices $x,y,z$,
and let $\mu$ be the $L$-invariant probability measure on $\Gamma\backslash L$. 

\begin{lemma}\label{lem:3.7} Let $a\in\cba(\DD^3,\RR)^\Gamma$ be a representative of the class
$\beta\in\hb^2(\Gamma,\RR)$.  Then
\bqn
\int_{\Gamma\backslash L}a(gx,gy,gz)d\mu(g)=\frac{\tau(\beta)}{2\pi} A(x,y,z)\,.
\eqn
\end{lemma}

We will also need to represent the relative cycle $[S,S_{\geq T}]$ using smearing 
in the context of relative measure homology.
For $\e\in S_2(\DD)$ we consider as usual the continuous map
\bqn
\ba
m_\e:\Gamma\backslash L&\to S_2(S)\\
\Gamma g&\mapsto p(g\e)
\ea
\eqn
where $p:\DD\to S$ is the canonical projection and define 
\bqn
S m_T(\e)=(m_\e)_\ast\big(\mu|_{(\Gamma\backslash L)_T}\big)\,,
\eqn
where $(\Gamma\backslash L)_T=\{\Gamma g\in\Gamma\backslash L: p(g0)\in S_{\leq T}\}$
and $S_{\leq T}$ is the complement of $S_{\geq T}$.
The following is then a verification proceeding along standard arguments.

\begin{lemma}\label{lem:3.8} Let $\sigma:\Delta^2\to\DD$ be a geodesic simplex
and $\sigma'$ its reflection alone one side.  
Then there is $C>1$ such that the boundary of the measured chain
\bqn
\mu_{CT}:=S m_{CT}(\sigma)-S m_{CT}(\sigma')
\eqn
has it support in $S_1(S_{\geq T})$ and $\mu_{CT}$ represents the relative cycle
\bqn
\frac{2A(\sigma)}{A(S)}[S,S_{\geq T}]
\eqn
where $A$ refers to the hyperbolic area.
\end{lemma}

\begin{proof}[{\it Proof of Proposition~\ref{prop:3.5}}] In the notation of Lemma~\ref{lem:3.6},
let $a\in\cba(D^3,\RR)^\Gamma$ be such that $\ol a$ is a representative of $g_S(\beta)$.
Then for $T_0$ large enough $\ol a$ restricted to $S_{\geq T_0}$ is trivial in bounded cohomology
and using Lemma~\ref{lem:3.6} applied to appropriate cusp neighborhoods we get a continuous bounded
function
\bqn
\ol f:S_1(S_{\geq T_0})\to\RR
\eqn
with $\ol a|_{S_2(S_{\geq T_0})}=d\ol f$.

For $T\geq T_0$ define $f_T:S_1(S)\to\RR$ as being equal to $\ol f$ on simplices in $S_{\geq T}$
and zero otherwise, and let $a_T:=\ol a-d f_T$.
Then $a_T$ is a bounded Borel function on $S_2(S)$ and $\|a_T\|_\infty\leq C$ for some
constant $C>0$.  Moreover $a_T$ is a representative of
$(j_T)^{-1}g_S(\beta)$.  For $T_2\geq T_1\geq T_0$ we clearly have
\bq\label{eq:3.5}
 \big\langle(j_{T_1})^{-1}g_S(\beta),[S,S_{\geq T_1}]\big\rangle
=\big\langle(j_{T_1})^{-1}g_S(\beta),[S,S_{\geq T_2}]\big\rangle\,.
\eq
According to Lemma~\ref{lem:3.8} the right hand side equals
\bqn
\left\langle a_{T_1},\frac{A(S)}{2A(\sigma)}\mu_{CT_2}\right\rangle
\eqn
which, letting $T_2\to\infty$, gives
\bqn
\frac{A(S)}{2A(\sigma)}\int_{\Gamma\backslash L}\big(a_{T_1}(pg(
\sigma))-a_{T_1}(pg(\sigma'))\big)d\mu(g)\,.
\eqn
Since however the left hand side of \eqref{eq:3.5} is independent of $T_1$, 
we let $T_1\to\infty$ and, using the dominated convergence theorem, obtain 
\bqn
\ba
  &\big\langle(j_{T})^{-1}g_S(\beta),[S,S_{\geq T}]\big\rangle\\
=&\frac{A(S)}{2A(\sigma)}\int_{\Gamma\backslash L}
  \big(a(g\sigma_0,g\sigma_1,g\sigma_2)-a(g\sigma'_0,g\sigma'_1,g\sigma'_2)\big)d\mu(g)
\ea
\eqn
which, together with Lemma~\ref{lem:3.7} proves the proposition.
\end{proof}

\subsection{Continuity}\label{subsec:3.4}
We will now use Theorem~\ref{thm:3.2} to show the following

\begin{prop}\label{prop:3.8} Let $G$ be a group of Hermitian type and 
let $\kappa\in\hcb^2(G,\RR)$.  Then the map
\begin{equation}\label{eq:map}
\ba
\T_\kappa(\Sigma,\,\cdot\,):\hom\big(\ps,G\big)&\longrightarrow\quad\RR\\
\rho\qquad&\mapsto\tksr
\ea
\end{equation}
is continuous.
\end{prop}

Together with the following basic example of maximal representation, the continuity of
$\rho\mapsto\tksr$ allows us to determine the range of the map in (\ref{eq:map}) 
when $\kappa = \kgb$ and $\bs\neq\emptyset$.

\begin{exo}\label{exo:3.9} Let $G$ be of Hermitian type with associated symmetric space $\Xx$,
$d:\DD\to\Xx$ a diagonal disk (see (\ref{sec:prelim}.b)) 
and $\Delta:L\to G$ the corresponding homomorphism,
where $L$ is an appropriate finite covering of $\mathrm{PU}(1,1)$.  
Then if $h:\ps\to\Gamma<L$ is a finite area hyperbolization of $\Sigma^\circ$,
the homomorphism $\rho:=\Delta\circ h$ is maximal.  Indeed, we have that
\bqn
\Delta^\ast(\kgb)=\rk_\Xx\klb
\eqn
and hence (Theorem~\ref{thm:3.2})
\bqn
 \T(\Sigma,\Delta\circ h) 
=|\chi(\Sigma)|\tb(\Delta|_\Gamma)
=|\chi(\Sigma)|\rk_\Xx\,.
\eqn
Observe that if $d'$ is the composition of $d$ with an antiholomorphic isometry of $\DD$
and $\Delta':L\to G$ is the corresponding homomorphism, then
$\T(\Sigma,\Delta'\circ h)=-|\chi(\Sigma)|\rk_\Xx$.
\end{exo}

\begin{cor}\label{cor:3.10} Assume that $\partial\Sigma\neq\emptyset$.  Then the range
of the map $T(\Sigma, \,\cdot\,)$ is the interval $\big[-\rk_\Xx|\chi(\Sigma)|,\rk_\Xx|\chi(\Sigma)|\big]$.
\end{cor}

\begin{proof} By Corollary~\ref{cor:3.3} the range is contained in the above interval and,
by Example~\ref{exo:3.9} it contains the endpoints.  Since $\bs\neq\emptyset$,
$\pi_1(\Sigma)$ is a free group and hence $\hpg$ is connected.
The corollary then follows from Proposition~\ref{prop:3.8}.
\end{proof}

Turning now to the proof of Proposition~\ref{prop:3.8} we will need the following:

\begin{lemma}\label{lem:3.11} Let $G$ be a connected semisimple
Lie group with associated symmetric space $\Xx$, 
let $\C(\DD,\Xx)$ be the space of continuous maps 
from the Poincar\'e disk $\DD$ into $\Xx$, with the topology
of uniform convergence on compact sets,
and let $\Gamma$ be a torsionfree lattice in $\PU(1,1)$.  
Then there is
a continuous map
\bqn
\ba
\hom(\Gamma,G)&\to\mathrm{C}(\DD,\Xx)\\
\rho\qquad&\mapsto \quad F_\rho
\ea
\eqn
such that $F_\rho$ is equivariant with respect to $\rho:\Gamma\to G$.
\end{lemma}

\begin{proof} Let $\Kk$ be a simplicial complex such that $|\Kk|$ 
is homeomorphic to $\Gamma\backslash\DD$.  
Let $\rho:\Gamma\to G$ be a homomorphism and $F:\widetilde \Kk^{(0)}\to\Xx$
a $\rho$-equivariant map defined on the $0$-skeleton of the universal covering of $\Kk$.
Using barycentric coordinates on the simplices of $\widetilde \Kk$
and the center of mass in $\Xx$, one obtains a canonical continuous
extension $F^{ext}:\widetilde \Kk\to\Xx$ which is thus $\rho$-equivariant
and depends continuously on $F$.  Fix $Y\subset\widetilde \Kk^{(0)}$
a complete set of representatives of $\Gamma$-orbits in $\widetilde \Kk^{(0)}$
and fix any map $f:Y\to\Xx$.  Then given $\rho:\Gamma\to G$,
we define $f_\rho:\widetilde \Kk^{(0)}\to\Xx$ 
as the unique $\rho$-equivariant extension of $f$ and 
$F_\rho:=(f_\rho)^{ext}$\,.  The assertion that $\rho\mapsto F_\rho$
is continuous follows from the continuity of the center of mass construction in $\Xx$.
\end{proof}

\begin{proof}[Proof of Proposition~\ref{prop:3.8}] We realize, as we may,
the bounded continuous cohomology of $G$ on the complex
$\big(\cba(\Xx^\bullet,\RR)\big)$ of bounded continuous alternating cochains on $\Xx$
and similarly for $L$ and $\Gamma$ on $\big(\cba(\DD^\bullet,\RR)\big)$.
Given a homomorphism $\rho:\G\to G$, the continuous $\rho$-equivariant
map $F_\rho\in\mathrm{C}(\DD,\Xx)$ in the previous lemma
induces by precomposition a map of complexes
\bqn
 \big(\cba(\Xx^\bullet,\RR)\big)^G\to
\big(\cba(\DD^\bullet,\RR)\big)^\Gamma
\eqn
which, according to \cite{Burger_Iozzi_app}, 
represents the pullback $\rho^\ast:\hcb^2(G,\RR)\to\hb^2(\Gamma,\RR)$.
In particular, if $c:\Xx^3\to\RR$ is a bounded continuous $G$-invariant 
alternating cocycle representing $\kappa\in\hcb^2(G,\RR)$, 
then the cocycle
\bqn
(z_1,z_2,z_3)\mapsto c\big(F_\rho(z_1),F_\rho(z_2),F_\rho(z_3)\big)
\eqn
represents $\rast(\k)\in\hb^2(\Gamma,\RR)$ \cite{Burger_Iozzi_app}. 
We deduce then from Lemma~\ref{lem:3.7} that 
\bqn
 \int_{\Gamma\backslash L}c\big(F_\rho(gz_1),F_\rho(gz_2),F_\rho(gz_3)\big)d\mu(g)
=\frac1{2\pi}\tb(\rho,\kappa)A(z_1,z_2,z_3)\,.
\eqn

If now $\rho_n\to\rho$, then according to Lemma~\ref{lem:3.11}
$F_{\rho_n}\to F_\rho$ uniformly on compact sets and hence
\bqn
c\big(F_{\rho_n}(gz_1),F_{\rho_n}(gz_2),F_{\rho_n}(gz_3)\big)\to
c\big(F_\rho(gz_1),F_\rho(gz_2),F_\rho(gz_3)\big)
\eqn
pointwise. Since
\bqn
    \big|c\big(F_{\rho_n}(gz_1),F_{\rho_n}(gz_2),F_{\rho_n}(gz_3)\big)\big|
\leq\|c\|_\infty\,,
\eqn
the dominated convergence theorem implies that 
$\tb(\rho_n,\kappa)\to\tb(\rho,\kappa)$.
\end{proof}

%%% Local Variables: 
%%% mode: latex
%%% TeX-master: "toledo"
%%% End: 

%% file: structure1.tex
\section{Structure of maximal representations: the Zariski dense case}\label{sec:structure1}
In this section we will investigate the structure of maximal 
homomorphisms $\rho:\Gamma\to G$, where, as before, 
$\Gamma<L$ is a lattice in a finite connected covering $L$ of $\PU(1,1)$, 
$G=\Iso(\Xx)^\circ$ is the connected component of the group
of isometries of an irreducible Hermitian symmetric space $\Xx$,
and we assume now that the image of $\rho$ is Zariski dense.  
More precisely, if $\gG$ is the connected adjoint $\RR$-group 
associated to the complexification of the Lie algebra of $G$,
we will prove the following:

\begin{thm}\label{thm:4.1}  If $\rho:\Gamma\to\gG(\RR)^\circ$ is a
maximal representation with Zariski dense image, then:
\be
\item the Hermitian symmetric space $\Xx$ is of tube type;
\item the image of $\rho$ is discrete;
\item the representation $\rho$ is injective, modulo possibly the center
$\Zz(\Gamma)$ of $\Gamma$.
\ee
\end{thm}

\subsection{The formula}\label{subsec:4.1}
Here we will use heavily the results of
\cite{Burger_Iozzi_Wienhard_kahler}.  In particular, 
let $\Dd$ be the bounded domain realization of $\Xx$ 
and $\cs$ its Shilov boundary.  Recall that $\cs=G/Q$
where $Q$ is a specific maximal parabolic subgroup of $G$; 
and denote by $\cs^{(2)}$ the set of pairs of transverse
points in $\cs$.

The lattice $\Gamma$ acts on the boundary of the Poincar\'e disk
$\partial\DD$ and, as is well known,
the space $(\partial\DD,\lambda)$ where $\lambda$ is the round
measure on $\partial\DD$ is a Poisson boundary for $\Gamma$;
moreover, the $\Gamma$-action on $\partial\DD\times\partial\DD$
is ergodic.  With this we can apply 
\cite[Proposition~7.2]{Burger_Iozzi_supq} and 
\cite[Theorem~4.7]{Burger_Iozzi_Wienhard_kahler} to conclude:

\begin{thm}\label{thm:4.2} Assume that $\rho:\Gamma\to\gG(\RR)^\circ=G$ is a homomorphism
with Zariski dense image.  Then there exists a $\rho$-equivariant
measurable map $\varphi:\partial\DD\to\cs$ such that 
$\big(\varphi(x_1),\varphi(x_2)\big)\in\cs^{(2)}$ for almost
all $(x_1,x_2)\in(\partial\DD)^2$.
\end{thm}

Using the boundary map $\varphi$, we now give an explicit cocycle 
on $(\partial\DD)^3$ representing the pullback 
$\rho^\ast(\kgb)\in\hb^2(\Gamma,\RR)$.
To this purpose, if $\beta_{\cs}:\cs^3\to\RR$ is the generalized Maslov cocycle (see \S~\ref{subsubsec:2.1.3}), 
we have:

\begin{cor}[{\cite[Proposition~4.6]{Burger_Iozzi_Wienhard_kahler}}]\label{cor:4.3} 
Under the canonical isomorphism
\bqn
\hb^2(\Gamma,\RR)\cong\Zz\linftya\big((\partial\DD)^3,\RR\big)^\Gamma
\eqn
the class $\rhoba(\kgb)$ corresponds to the cocycle
\bqn
\ba
(\partial\DD)^3\,&\longrightarrow\qquad\quad\RR\\
(x,y,z)&\mapsto\beta_{\cs}\big(\varphi(x),\varphi(y),\varphi(z)\big)\,.
\ea
\eqn
\end{cor}

We then conclude:

\begin{cor}\label{cor:4.4} Let $\rho:\Gamma\to G=\gG(\RR)^\circ$ be a homomorphism
with Zariski dense image and $\varphi:\partial\DD\to\cs$ a measurable
$\rho$-equivariant boundary map.  Then 
if $\mu$ is the $L$-invariant probability measure on $\Gamma\backslash L$, 
we have that 
\bq\label{eq:formula}
\int_{\Gamma\backslash L}
    \beta_{\cs}\big(\varphi(gx),\varphi(gy),\varphi(gz)\big)d\mu(g)=
\tb(\rho)\beta_{\partial\DD}(x,y,z)
\eq
for almost every $(x,y,z)\in(\partial\DD)^3$.
In particular, if $\rho$ is maximal,
\bq\label{eq:bb}
\beta_{\cs}\big(\varphi(x),\varphi(y),\varphi(z)\big)=
\rk_\Xx\,\beta_{\partial\DD}(x,y,z)
\eq
for almost every $(x,y,z)\in\partial\DD$
and thus
\bq\label{eq:4.3}
\rho^\ast(\kgb)=\r_\Xx\klb\,.
\eq
\end{cor}

\begin{proof} We have that 
\bqn
\mathrm{T}_\mathrm{b}\big(\rho^\ast(\kgb)\big)=\tb(\rho)\klb\,,
\eqn
so that the formula follows from Corollary~\ref{cor:4.3}
and the functoriality of the transfer operator in \cite[III.8]{Monod_book}.

Assume now that $\rho$ is maximal, that is $\tb(\rho)=\r_\Xx$.
Fix $(x_0,y_0,z_0)\in(\partial\DD)^3$ such that 
$\beta_{\partial\DD}(x_0,y_0,z_0)=\frac12$ and (\ref{eq:formula}) holds.
Since for every $a,b,c \in\cs$
\bqn
|\beta_{\cs}(a,b,c)|\leq\frac{\rk_\Xx}{2}\,,
\eqn
(\cite{Clerc_Orsted_2} -- see also \cite[Theorem~4.2]{Burger_Iozzi_Wienhard_kahler})
and since $\mu$ is a probability measure,
we deduce from (\ref{eq:formula}) that for almost every $g\in L$
\bq\label{eq:eq}
 \beta_{\cs}\big(\varphi(g x_0),\varphi(g y_0),\varphi(g z_0)\big)
=\rk_\Xx\beta_{\partial\DD}(x_0,y_0,z_0)\,.
\eq
Similarly, if $\beta_{\partial\DD}(x_0,y_0,z_0)= -\frac12$,
by the same argument we deduce that (\ref{eq:eq}) holds for almost every $g\in L$.
Thus the function on $(\partial\DD)^3$
\bqn
(x,y,z)\mapsto\beta_{\cs}\big(\varphi(x),\varphi(y),\varphi(z)\big)
\eqn
is essentially $L$-invariant and (\ref{eq:eq}) implies then 
that this function coincides almost everywhere 
with $\rk_\Xx\beta_{\partial\DD}$.
\end{proof}

\subsection{$\Xx$ is of tube type}\label{subsec:4.2}
In this section we prove the first assertion of Theorem~\ref{thm:4.1}.
For this we will use our characterization of tube type domains
obtained in \cite{Burger_Iozzi_Wienhard_kahler}.  
In particular we defined in \cite[\S~2.4]{Burger_Iozzi_Wienhard_kahler} 
the {\it Hermitian triple product}, a $G$-invariant map
\bqn
\<\<\,\cdot\,,\,\cdot\,,\,\cdot\,\>\>:\cs^{(3)}\to\RR^\times\backslash\CC^\times
\eqn
on the set $\cs^{(3)}$ of triples of points in $\cs$ 
which are pairwise transverse, and which is related to the generalized Maslov cocycle by
\bqn
\<\<x,y,z\>\>\equiv e^{i\pi p_\Xx\beta_{\cs}(x,y,z)}\mod\RR^\times
\eqn
for all $(x,y,z)\in\cs^{(3)}$, where $p_\Xx$ is an integer 
defined in terms of the root system associated to $G$.  Let $\cs=G/Q$
and $\qQ$ be the $\RR$-parabolic subgroup of $\gG$ with $\qQ(\RR)=Q$.

Then if $A^\times$ is the $\RR$-algebraic group $\CC^\times\times\CC^\times$
(with real structure $(\lambda,\mu)\mapsto(\overline\mu,\overline\lambda)$)
and $\CC^\times\one=\{(\lambda,\lambda)\in A:\,\lambda\in\CC^\times\}$,
we constructed a rational $\gG$-invariant map defined over $\RR$
\bqn
\<\<\,\cdot\,,\,\cdot\,,\,\cdot\,\>\>_\CC:(\gG/\qQ)^3\to\CC^\times\one\backslash A^\times\,,
\eqn
which we called the {\it complex Hermitian triple product},
and which is related to the Hermitian triple product
by the commutative diagram
\bqn
\xymatrix{
 (\gG/\qQ)^3\ar[rr]^{\<\<\,\cdot\,,\,\cdot\,,\,\cdot\,\>\>_\CC} &
&\CC^\times\one\backslash A^\times\\
 \cs^{(3)}\ar[rr]^{\<\<\,\cdot\,,\,\cdot\,,\,\cdot\,\>\>}
                   \ar[u]^{\imath^3} &
&\RR^\times\backslash\CC^\times\ar[u]_\Delta
}
\eqn
where $\imath$ is given by the $G$-map $\imath:\cs\to\gG/\qQ$ 
sending $\cs$ to $(\gG/\qQ)(\RR)$, and 
$\Delta\big([\lambda]\big)=\big[(\lambda,\overline\lambda)\big]$,
(see \cite[Corollary~2.11]{Burger_Iozzi_Wienhard_kahler});
the statement includes the fact that 
the domain of definition of $\<\<\,\cdot\,,\,\cdot\,,\,\cdot\,\>\>_\CC$ 
contains $\cs^{(3)}$.

Given now $(a,b)\in\cs^{(2)}$, let as in 
\cite[\S~5.1]{Burger_Iozzi_Wienhard_kahler}
\bqn
\Oo_{a,b}\subset\gG/\qQ 
\eqn
be the Zariski open subset on which the map
\bqn
\ba
p_{a,b}:\Oo_{a,b}&\to\,\CC^\times\one\backslash A^\times\\
x\,\,\,&\mapsto\<\<a,b,x\>\>_\CC
\ea
\eqn
is defined.  We have then (see \cite[Lemma~5.1]{Burger_Iozzi_Wienhard_kahler})
that if for some $m\in\ZZ\setminus\{0\}$ the map
\bqn
\ba
\Oo_{a,b}&\to\CC^\times\one\backslash A^\times\\
x\,\,\,&\mapsto \,\,p_{a,b}(x)^m
\ea
\eqn
is constant, then $\Xx$ is of tube type.  
Now we apply (\ref{eq:bb}) in Corollary~\ref{cor:4.4}
to get that if $\rho:\Gamma\to G=\gG(\RR)^\circ$ is maximal, 
then 
\bqn
\beta_{\cs}\big(\varphi(x),\varphi(y),\varphi(z)\big)=\pm\frac{\rk_\Xx}{2}
\eqn
for almost every $(x,y,z)$, which implies that 
\bq\label{eq:chtp=1}
\big<\big\<\varphi(x),\varphi(y),\varphi(z)\big\>\big\>^2\equiv1\mod\RR^\times\,.
\eq
In particular, fix $x,y$ such that $\big(\varphi(x),\varphi(y)\big)\in\cs^{(2)}$
and such that (\ref{eq:chtp=1}) holds for almost all $z\in\partial\DD$:
letting $E\subset \partial \DD$ be this set of full measure, 
we may assume that $E$ is $\Gamma$-invariant and $\varphi(E)\subset\Oo_{a,b}$
where $a=\varphi(x)$ and $b=\varphi(y)$.  But $\varphi(E)$ being
$\rho(\Gamma)$-invariant is Zariski dense in $\gG/\qQ$ and 
hence Zariski dense in the open set $\Oo_{a,b}$;
since the map $x\mapsto p_{a,b}(x)^2$ is constant on $\varphi(E)$
it is so on $\Oo_{a,b}$, which implies that $\Xx$ is of tube type.

\subsection{The image of $\rho$ is discrete}\label{subsec:discrete}
Under the hypothesis of Theorem~\ref{thm:4.1}, we know now that
$\Xx$ is of tube type.  Then the generalized Maslov cocycle $\beta_{\cs}$
takes on $\cs^{(3)}$ exactly $\rk_\Xx+1$ values, namely
\bq\label{eq:values}
\left\{-\frac{\rk_\Xx}{2},-\frac{\rk_\Xx}{2}+1,\dots,\frac{\rk_\Xx}{2}-1, \frac{\rk_\Xx}{2}\right\}
\eq
so that 
\bqn
\cs^{(3)}=\cup_{i=0}^{\rk_\Xx}\Oo_{-\rk_\Xx+2i}\,,
\eqn
where $\Oo_{-\rk_\Xx+2i}$ is the preimage via $\beta_{\cs}$ of $\frac{\rk_\Xx}{2}+i$,
which incidentally is open since $\beta_{\cs}$ is continuous on $\cs^{(3)}$
(see \cite[Corollary~3.7]{Burger_Iozzi_Wienhard_kahler}).
With the above notations, it follows from (\ref{eq:bb}) in Corollary~\ref{cor:4.4} that
\bq\label{eq:r-r}
\big(\varphi(x),\varphi(y),\varphi(z)\big)\in\Oo_{-\rk_\Xx}\cup\Oo_{\rk_\Xx}
\eq
for almost all $(x,y,z)\in(\partial\DD)^3$.
Let us now denote by $\essim\varphi\subset\cs$ the essential image of $\varphi$,
that is the support of the pushforward $\varphi_\ast(\lambda)$
of the round measure $\lambda$ on $\partial\DD$.  
Then $\essim\varphi$ is closed and $\rho(\Gamma)$-invariant.
It follows then from (\ref{eq:r-r}) that 
\bqn
(\essim\varphi)^3\subset\overline{\Oo_{-\rk_\Xx}}\cup\overline{\Oo_{\rk_\Xx}}\,,
\eqn
where the closure on the right hand side is taken in $\cs^3$.
There are now two cases. Either $\rk_\Xx=1$, $\Xx=\DD$, $G=\PU(1,1)$,
which is the case treated in \cite{Burger_Iozzi_difff};
or $\rk_\Xx\geq2$, and then $\overline{\Oo_{-\rk_\Xx}}\cup\overline{\Oo_{\rk_\Xx}}$
is not the whole of $\cs^{(3)}$, since its complement contains 
at least $\Oo_{-\rk_\Xx+2}$; thus $(\essim\varphi)^3\neq\cs^3$
and, since $(\essim\varphi)^3$ is $\rho(\Gamma)^3$-invariant
closed and $\cs$ is $G$-homogeneous, this implies that 
$\rho(\Gamma)$ is not dense in $G$.  Since a Zariski dense subgroup of $G$ is either discrete or dense,  
we have that $\rho(\Gamma)$ is discrete.

\subsection{The representation $\rho$ is injective}\label{subsec:4.4}
Assume that 
\bqn
\ker(\rho)\not<\Zz(\Gamma)\,.
\eqn
Then it is easy to see that there is $\gamma\in\ker\rho$ of infinite order,
so that we may choose a nonempty open interval $I\subset\partial\DD$ 
such that $I,\gamma I,\gamma^2 I$ are pairwise disjoint and positively
oriented, that is $\beta_{\DD}(x,y,z)=\frac12$ for all 
$(x,y,z)\in I\times\gamma I\times\gamma^2 I$.
Now choose three open nonvoid intervals $I_1,I_2,I_3$ in $I$
which are pairwise disjoint and positively oriented.
Then it follows from (\ref{eq:bb}) in Corollary~\ref{cor:4.4} that 
\bqn
E_1:=\left\{(x,y,z)\in I_3\times\gamma I_2\times\gamma^2 I_1:\,
\beta_{\cs}\big(\varphi(x),\varphi(y),\varphi(z)\big)=\frac{\rk_\Xx}{2}\right\}
\eqn
is of full measure in $I_3\times\gamma I_2\times\gamma^2 I_1$, while
\bqn
E_2:=\left\{(x,y,z)\in I_3\times I_2\times I_1:\,
\beta_{\cs}\big(\varphi(x),\varphi(y),\varphi(z)\big)=-\frac{\rk_\Xx}{2}\right\}
\eqn
is of full measure in $I_3\times I_2\times I_1$.
Thus 
\bqn
E_1':=\big\{(x,y,z)\in E_1:\,(x,\gamma^{-1}y,\gamma^{-2}z)\in E_2\big\}
\eqn
is of full measure; using the almost everywhere equivariance of $\varphi$
and the assumption that $\rho(\gamma)=\id$, 
we conclude that for almost every $(x,y,z)\in E_1'$
\bqn
  \frac{\rk_\Xx}{2}
=\beta_{\cs}\big(\varphi(x),\varphi(y),\varphi(z)\big)
=\beta_{\cs}\big(\varphi(x),\varphi(\gamma^{-1}y),\varphi(\gamma^{-2}z)\big)
=-\frac{\rk_\Xx}{2}\,,
\eqn
which is a contradiction.
This shows that $\ker\rho\subset\Zz(\Gamma)$ and completes the proof of
Theorem~\ref{thm:4.1}.

%%% Local Variables: 
%%% mode: latex
%%% TeX-master: "toledo"
%%% End: 

%% file: regularity.tex
\section{Regularity properties of the boundary map: the Zariski dense case}\label{sec:regularity}
In this section we generalize a technique from \cite{Burger_Iozzi_Labourie_Wienhard}
to study the boundary map $\varphi:\partial\DD\to\cs$
associated to a maximal representation $\rho:\Gamma\to G=\gG(\RR)^\circ$
with Zariski dense image and establish the existence of strictly
equivariant maps with additional regularity properties.  
This relies in an essential way on the fact established in Theorem~\ref{thm:4.1}
asserting that, in the situation described above, the symmetric space $\Xx$
associated to $G$ is irreducible and of tube type. 

\begin{thm}\label{thm:5.1} Let $\Gamma$ be a lattice in 
a finite connected covering of $\PU(1,1)$ and let
$\rho:\Gamma\to G$ a maximal representation with Zariski dense image.  
Then there are two Borel maps
\bqn
\varphi_\pm:\partial\DD\to\cs
\eqn
with the following properties:
\be
\item $\varphi_+$ and $\varphi_-$ are strictly $\rho$-equivariant;
\item $\varphi_-$ is left continuous and $\varphi_+$ is right continuous;
\item for every $x\neq y$, $\varphi_\epsilon(x)$ is transverse to $\varphi_\delta(y)$
for all $\epsilon,\delta\in\{+,-\}$;
\item for all $x,y,z\in\partial\DD$,
\bqn
 \beta_{\cs}\big(\varphi_\epsilon(x),\varphi_\delta(y),\varphi_\eta(z)\big)
=\rk_\Xx\,\beta_{\partial\DD}(x,y,z)\,,
\eqn
for all $\epsilon,\delta,\eta\in\{+,-\}$.
\ee
Moreover $\varphi_+$ and $\varphi_-$ are the unique maps satisfying (1) and (2).
\end{thm}

\subsection{General properties of boundary maps}\label{subsec:general}
Let $\Gamma<L$ be a lattice in a connected finite covering $L$ of $\PU(1,1)$ 
as above, $\gG$ a connected semisimple group, 
$\pP$ a parabolic subgroup, both defined over $\RR$, 
$G=\gG(\RR)$, $P=\pP(\RR)$, $\rho:\Gamma\to G$ a homomorphism and
$\lambda$ the round measure on $\partial\DD$.  
If $\varphi:\partial\DD\to G/P$ is a $\rho$-equivariant measurable map
and $\rho$ has Zariski dense image, one sees immediately that the image
of $\varphi$ cannot be contained in a proper algebraic subset.
The following proposition is a strengthening of this statement, 
showing that from the point of view of the round measure,
the essential image of $\varphi$ meets any proper algebraic subset 
in a set of measure zero. Namely:

\begin{prop}\label{prop:5.2} 
If $\vV\subset\gG/\pP$ is any proper Zariski closed subset 
defined over $\RR$, then
\bqn
\lambda\big(\varphi^{-1}(\vV(\RR))\big)=0\,.
\eqn
\end{prop}
This will follow from the following two lemmas.

\begin{lemma}\cite{Kaimanovich_ern}\label{lem:5.3} If $A\subset\partial\DD$ is a set of
positive measure, then there exists a sequence $\{\gamma_n\}_{n=1}^\infty$ in $\Gamma$
such that $\lim_{n\to\infty}\lambda(\gamma_nA)=1$.
\end{lemma}
\begin{lemma}\label{lem:5.4} Let $\vV\subset\gG/\pP$ be a proper 
Zariski closed subset defined over $\RR$ and $\{g_n\}_{n=1}^\infty$
a sequence in $G$.  Then there exists a proper Zariski closed
subset $\wW\subset\gG/\pP$ defined over $\RR$ and a subsequence 
$\{g_{n_k}\}_{k=1}^\infty$
such that for every $\epsilon>0$ there exists $K>0$ such that 
for all $k\geq K$
\bqn
g_{n_k}\vV(\RR)\subset \Nn_\epsilon\big(\wW(\RR)\big)\,,
\eqn
where $\Nn_\epsilon$ denotes the $\epsilon$-neighborhood for some fixed
distance on $G/P$.
\end{lemma}

\begin{proof} Passing to a subsequence $\{g_{n_k}\}$ 
we may assume that, if $V:=\vV(\RR)$, $g_{n_k}V$ converges 
to a compact set $F\subset G/P$ in the Gromov--Hausdorff topology.  
It suffices then to show that $F$ is not Zariski dense in $\gG/\pP$.  
To this end, let $I(\vV)\subset\CC[\gG/\pP]$ be the defining ideal 
and pick $d\geq1$ such that the homogeneous component
$I(\vV)_d(\RR)\subset\CC[\gG/\pP]_d(\RR)$ is nonzero. 
Let $\ell:=\dim I(\vV)_d(\RR)$; then we may assume that 
the sequence $g_{n_k} I(\vV)_d(\RR)$ converges to a point $E$
in the Grassmannian $\Gr_\ell\big(\CC[\gG/\pP]_d(\RR)\big)$.
In particular, given $p\in E$, $p\neq0$, 
there exists $p_n\in g_{n_k}I(\vV)_d(\RR)$ with $p_n\to p$.  
It is then not difficult to show that $p$ vanishes on $F$
and one can take $\wW\subset\gG/\pP$ to be the Zariski closure of $F$.
\end{proof}

\begin{proof}[{\it Proof of Proposition~\ref{prop:5.2}}]
Let $V:=\vV(\RR)$ and $A:=\{x\in\partial\DD:\,\varphi(x)\in V\}$.
Assume that $\lambda(A)>0$ and pick a sequence $\{\gamma_n\}_{n=1}^\infty$ 
in $\Gamma$ as in Lemma~\ref{lem:5.3} such that 
\bqn
\lambda(\gamma_nA)\to1\,.
\eqn
Passing to a subsequence, let $\wW\subset\gG/\pP$ be 
the proper Zariski closed subset given by Lemma~\ref{lem:5.4}.
For every $m\geq1$, let $N(m)$ be such that for all $n\geq N(m)$
\bqn
\rho(\gamma_n)V\subset\Nn_{1/m}(W)\,,
\eqn
where $W:=\wW(\RR)$.  Setting $E_N:=\cup_{n=N}^\infty\gamma_nA$,
we have thus
\bq\label{eq:neigh}
\varphi\big(E_{N(m)}\big)\subset\Nn_{1/m}(W)\,.
\eq
But now $E_{N(m)}\subset\partial\DD$ is a set of full measure
and so is $E:=\cap_{m\geq1}E_{N(m)}$, which, by (\ref{eq:neigh}),
implies now that $\varphi(E)\subset W$. 
This implies that 
\bqn
\essim(\varphi)\subset W\subset\wW\subset\gG/\pP
\eqn
and contradicts the Zariski density of $\rho(\Gamma)$
since $\essim(\varphi)$ is $\rho(\Gamma)$-invariant.
\end{proof}

\subsection{Exploiting maximality}\label{subsec:exploiting} 
Let now $\rho:\Gamma\to G$ be a maximal representation with Zariski dense image 
and $\varphi:\partial\DD\to\cs$ be the $\rho$-equivariant measurable
map given by Theorem~\ref{thm:4.2}.  Having introduced the essential image
$\essim(\varphi)\subset\cs$, we will now study the essential graph
$\essgr(\varphi)\subset\partial\DD\times\cs$ of $\varphi$
defined as the support of the pushforward of the round measure $\lambda$
on $\partial\DD$ under the map 
\bqn
\ba
\partial\DD&\to\,\,\partial\DD\times\cs\\
x\,\,&\mapsto\big(x,\varphi(x)\big)\,.
\ea
\eqn
For this we will use (\ref{eq:bb}) in Corollary~\ref{cor:4.4} in an essential way.
The following properties of the generalized Maslov cocycle $\beta_{\cs}$ follow from
\cite{Clerc_maslov_tube, Clerc_Orsted_2, Clerc_Neeb}:

\begin{lemma}\label{lem:5.6}
\be
\item $\beta_{\cs}:\cs^3\to\{-\frac{\rk_\Xx}{2}\}+\ZZ$ is a $G$-invariant cocycle;
\item $|\beta_{\cs}(x,y,z)|\leq \frac{\rk_\Xx}{2}$;
\item if $\beta_{\cs}(x,y,z)=\frac{\rk_\Xx}{2}$ then $x,y,z$ are pairwise transverse;
\item $\cs^{(3)}=\sqcup_{i=0}^{\rk_\Xx}\Oo_{-\rk_\Xx+2i}$, where $\Oo_{-\rk_\Xx+2i}$
is open in $\cs^3$ and $\beta_{\cs}$ takes on the value $-\frac{\rk_\Xx}{2}+i$
on $\Oo_{-\rk_\Xx+2i}$;
\item if $x,\{x_n\},\{x_n'\}\in\cs$, $\lim_{n\to\infty}x_n=x$ and
$\beta_{\cs}(x,x_n',x_n)=\frac{\rk_\Xx}{2}$, then $\lim_{n\to\infty}x_n'=x$.
\ee
\end{lemma}

Finally, for $x\in\cs=G/Q\subset\gG/\qQ$, let $\vV_x\subset\gG/\qQ$
be the proper Zariski closed $\RR$-subset of all points $y\in\gG/\qQ$
which are not transverse to $x$, so that $V_x:=\vV_x(\RR)$
is the set of points in $\cs$ which are not transverse to $x$.

\begin{lemma}\label{lem:5.7} Let $(x_1,f_1), (x_2,f_2), (x_3,f_3)$
be points in $\essgr(\varphi)$ so that $x_1,x_2,x_3$
are pairwise distinct and $f_1,f_2,f_3$ are pairwise transverse.
Then
\bqn
\beta_{\cs}(f_1,f_2,f_3)=\rk_\Xx\, \beta_{\partial\DD}(x_1,x_2,x_3)\,.
\eqn
\end{lemma}

\begin{proof} Let $I_i$, $i=1,2,3$ be pairwise disjoint open intervals
containing $x_i$ such that  for all $y_i\in I_i$
\bqn
\beta_{\partial\DD}(y_1,y_2,y_3)=\beta_{\partial\DD}(x_1,x_2,x_3)\,.
\eqn
Let $U_i$, $i=1,2,3$ be neighborhoods of $f_i$ such that 
$U_1\times U_2\times U_3\subset\cs^{(3)}$.  Then 
\bqn
A_i=\{x\in I_i:\varphi(x)\in U_i\}
\eqn
is of positive measure and hence it follows from (\ref{eq:bb}) in Corollary~\ref{cor:4.4}  
that for almost every $(y_1,y_2,y_3)\in A_1\times A_2\times A_3$,
\bqn
 \beta_{\cs}\big(\varphi(y_1), \varphi(y_2), \varphi(y_3)\big)
=\rk_\Xx\,\beta_{\partial\DD}(y_1,y_2,y_3)
=\rk_\Xx\,\beta_{\partial\DD}(x_1,x_2,x_3)\,.
\eqn
Thus setting $\varepsilon=2\beta_{\partial\DD}(x_1,x_2,x_3)\in\{\pm1\}$,
we have for almost every $(y_1,y_2,y_3)\in A_1\times A_2\times A_3$, that
\bqn
\big(\varphi(y_1),\varphi(y_2),\varphi(y_3)\big)
\in(U_1\times U_2\times U_3)\cap\Oo_{\varepsilon \rk_\Xx}\,,
\eqn
which implies, since the neighborhood $U_i$ can be chosen
arbitrarily small, that $(f_1,f_2,f_3)\in\overline{\Oo_{\varepsilon \rk_\Xx}}$.
But
\bqn
 \overline{\Oo_{\varepsilon \rk_\Xx}}\cap\cs^{(3)}
=\overline{\Oo_{\varepsilon \rk_\Xx}}\cap\big(\cup_{i=0}^{\rk_\Xx}\Oo_{-\rk_\Xx+2i}\big)
=\Oo_{\varepsilon \rk_\Xx}
\eqn
which, together with the assumption that $(f_1,f_2,f_3)\in\cs^{(3)}$,
implies that $(f_1,f_2,f_3)\in\Oo_{\varepsilon \rk_\Xx}$ 
and hence proves the lemma.
\end{proof}

\begin{lemma}\label{lem:5.8} Let $(x_1,f_2),(x_2,f_2)\in \essgr(\varphi)$
with $x_1\neq x_2$.  Then $f_1$ is transverse to $f_2$.
\end{lemma}
\begin{proof} For $x,y\in\partial\DD$, let 
\bqn
((x,y)):=\left\{z\in\partial\DD:\,\beta_{\partial\DD}(x,z,y)=\frac{1}{2}\right\}\,.
\eqn
We will use the obvious fact that for almost every $x\in\partial\DD$,
$\big(x,\varphi(x)\big)\in \essgr(\varphi)$.  Using Proposition~\ref{prop:5.2},
we can find $a\in((x_1,x_2))$ such that $\big(a,\varphi(a)\big)\in \essgr(\varphi)$
and $\varphi(a)\notin V_{f_1}\cup V_{f_2}$, that is
$\varphi(a)$ is transverse to $f_1$ and $f_2$. 
Then by the same argument, we can find $b\in((x_2,x_1))$ such that 
$\big(b,\varphi(b)\big)\in \essgr(\varphi)$, and $\varphi(b)$ is transverse
to $f_1,f_2$ and $\varphi(a)$.  Applying the cocycle property of $\beta_{\cs}$
and Lemma~\ref{lem:5.7}, we obtain
\bqn
\ba
  0
=&\beta_{\cs}\big(\varphi(a),f_2,\varphi(b)\big)
 -\beta_{\cs}\big(f_1,f_2,\varphi(b)\big)\\
 &\hphantom{XXX}+\beta_{\cs}\big(f_1,\varphi(a),\varphi(b)\big)
 -\beta_{\cs}\big(f_1,\varphi(a),f_2\big)\\
=&\frac{\rk_\Xx}{2}-\beta_{\cs}\big(f_1,f_2,\varphi(b)\big)
 +\frac{\rk_\Xx}{2}-\beta_{\cs}\big(f_1,\varphi(a),f_2\big)\,,
\ea
\eqn
which, together with Lemma~\ref{lem:5.6}(2) implies that 
$\beta_{\cs}\big(f_1,f_2,\varphi(b)\big)=\frac{\rk_\Xx}{2}$; using Lemma~\ref{lem:5.6}(3) 
we conclude that $f_1$ and $f_2$ are transverse.
\end{proof}

For a subset $A\subset\partial\DD$, let 
\bqn
F_A=\big\{f\in\cs:\,\hbox{ there exists }x\in A\hbox{ such that }(x,f)\in \essgr(\varphi)\big\}\,,
\eqn
and set 
\bqn
((x,y]]:=((x,y))\cup\{y\}\,.
\eqn

\begin{lemma}\label{lem:5.9} Let $x\neq y$ in $\partial\DD$.  Then
$\overline{F_{((x,y]]}}\cap F_x$ and $\overline{F_{[[y,x))}}\cap F_x$
consist each of one point.
\end{lemma}
\begin{proof} We start with two observations: first, if $A\cap B=\emptyset$, 
it follows from Lemma~\ref{lem:5.8} that $F_A\cap F_B=\emptyset$;
moreover, if $A$ is closed then $F_A$ is also closed.
We prove now that $\overline{F_{((x,y]]}}\cap F_x$ consists of one point,
the other statement can be proved analogously.

Let $f,f'\in\overline{F_{((x,y]]}}\cap F_x$
and let $(x_n,f_n)\in \essgr(\varphi)$ be a sequence such that 
\bqn
x_n\in((x,y]],\qquad\lim x_n=x, \text{ and} \lim f_n=f\,.
\eqn
Observe now that if $z\in((x,y]]$, writing $F_{((x,y]]}=F_{((x,z))}\cup F_{[[z,y]]}$
and taking into account that $F_{[[z,y]]}$ is closed and disjoint from $F_x$,
we get that
\bqn
\overline{F_{((x,y]]}}\cap F_x=\overline{F_{((x,z))}}\cap F_x\,,
\eqn
and hence also
\bq\label{eq:fxy}
\overline{F_{((x,y]]}}\cap F_x=\overline{F_{((x,z]]}}\cap F_x
\eq
for all $z\in((x,y))$.  
Using (\ref{eq:fxy}), we may find another sequence
$(y_n,f_n')\in \essgr(\varphi)$ such that 
\bqn
y_n\in((x,x_n)),\qquad\lim y_n=x, \text{ and} \lim f_n'=f'\,.
\eqn
Then it follows from Lemma~\ref{lem:5.7} and Lemma~\ref{lem:5.8}
that 
\bqn
\beta_{\cs}(f,f_n',f_n)=\rk_\Xx \beta_{\partial\DD}(x,y_n,x_n)=\frac{\rk_\Xx}{2}\,.
\eqn
Since $\lim f_n=f$, by Lemma~\ref{lem:5.6}(5), this however implies that $\lim f_n'=f$ 
and hence $f=f'$.
\end{proof}

Here is an interesting corollary about the structure of $\essgr(\varphi)$
which spells out precisely to which extent $\essgr(\varphi)$ is in general not 
the graph of a map.

\begin{cor}\label{cor:5.10} For every $x\in\partial\DD$, 
$F_x$ consists of one or two points.
\end{cor}
\begin{proof} Pick $y_-,x,y_+$ positively oriented in $\partial\DD$ and $f\in F_x$.
For every neighborhood $U$ of $f$, one of the sets
\bqn
\ba
&\big\{z\in[[y_-,x)):\,\varphi(z)\in U\big\}\\
&\big\{z\in((x,y_+]]:\,\varphi(z)\in U\big\}
\ea
\eqn
is of positive measure.  This implies that 
\bqn
F_x\subset\overline{F_{[[y_-,x))}}\cup\overline{F_{((x,y_+]]}}
\eqn
and thus 
\bqn
F_x=\big(\overline{F_{[[y_-,x))}}\cap F_x\big)\cup\big(\overline{F_{((x,y_+]]}}\cap F_x\big)\,,
\eqn
which, together with Lemma~\ref{lem:5.9} proves the assertion.
\end{proof}

\begin{proof}[Proof of Theorem~\ref{thm:5.1}] 
We use Lemma~\ref{lem:5.9} in order to define for every $x\in\partial\DD$
\bq\label{eq:f+f-}
\ba
\varphi_-(x)=&\overline{F_{((x,y_+]]}}\cap F_x\\
\varphi_+(x)=&\overline{F_{[[y_-,x))}}\cap F_x\
\ea
\eq
where $y_+\neq x$ and $y_-\neq x$ are arbitrary.  Then $\varphi_+$ and $\varphi_-$
are clearly respectively right and left continuous.
The strict $\rho$-equivariance of $\varphi_+$ and $\varphi_-$ 
follows from the invariance of $\essgr(\varphi)\subset\partial\DD\times\cs$
under the diagonal $\Gamma$-action together with (\ref{eq:f+f-})
and the fact that $\Gamma$ acts in an orientation preserving way on
$\partial\DD$. 

Properties (3) and (4) are immediate consequences of Lemmas~\ref{lem:5.7}
and \ref{lem:5.8}.

Concerning the uniqueness of $\varphi_+$ for example, 
let $\psi$ be a right continuous $\rho$-equivariant map.
Since the $\rho(\Gamma)$-action on $\cs$ is proximal, 
we have $\psi(x)=\varphi_+(x)$ for almost every $x\in\partial\DD$ \cite{Furstenberg_63}.
Fix $x\in\partial\DD$: then pick a sequence $y_n\in((x,y))$ such that
\bqn
\lim y_n=x \text{ and }\psi(y_n)=\varphi_+(y_n)\,.
\eqn
This implies, since $\psi$ and $\varphi_+$ are both right continuous, 
that $\psi(x)=\varphi_+(x)$.
\end{proof}

%%% Local Variables: 
%%% mode: latex
%%% TeX-master: "toledo"
%%% End: 

%% file: structure2.tex
\section{Structure of maximal representations and boundary maps: the general case}\label{sec:structure2}
In this section we present the proofs of Theorem~\ref{thm_intro:thm3} and
Theorem~\ref{thm_intro:thm5} in the introduction;
this relies on the results obtained in \S\S~\ref{sec:structure1} and \ref{sec:regularity}
and on the relation between maximal representations and tight homomorphisms, which were introduced in \cite{Burger_Iozzi_Wienhard_tight}. 

Let $G$ be a group of type (RH).  We briefly recall the structure of $\hcb^2(G,\RR)$
in terms of simple components and the explicit form of the Gromov norm.
Let $\Xx$ be the symmetric space associated to $G$; 
setting $\gx:=\Iso(\Xx)^\circ$ we have a canonical projection
$q:G\to\gx$ which by hypotheses has compact kernel.
Let $\Xx=\Xx_1\times\dots\times\Xx_n$ be the decomposition into irreducible factors
and $p_i:\gx\to\gxi$ the canonical projections.
We have then the following isometric isomorphisms:
\be
\item $q^\ast:\hcb^2(\gx,\RR)\to\hcb^2(G,\RR)$;
\item 
%\bqn
%\ba
%\prod\hcb^2(\gxi,\RR)&\to\hcb^2(\gx,\RR)\\
%(\alpha_i)&\mapsto\sum_{i=1}^n p_i^\ast(\alpha_i)\,.
%\ea
%\eqn
$\prod\hcb^2(\gxi,\RR)\to\hcb^2(\gx,\RR)$, $(\alpha_i)\mapsto\sum_{i=1}^n p_i^\ast(\alpha_i)$.
\ee
Defining for ease of notation $\kxb$ to be the bounded K\"ahler class of $\gx$,
let 
\bqn
\kxib:=p_i^\ast(\kappa_{\Xx_i}^\mathrm{b})\qquad\hbox{ and }\qquad\kgib:=q^\ast(\kxib)\,.
\eqn
Then 
\bqn
\{\kgib:\,1\leq i\leq n\}
\eqn
is a basis of $\hcb^2(G,\RR)$ and the norm of an element
\bqn
\kappa=\sum_{i=1}^n\lambda_i\kgib
\eqn
equals (see \cite[(2.15)]{Burger_Iozzi_Wienhard_tight}
\bq\label{eq:norm}
\|\kappa\|=\sum_{i=1}^n|\lambda_i|\frac{\r_{\Xx_i}}{2}\,.
\eq
In the rest of this section $L$ will always denote a finite connected covering 
of $\PU(1,1)$ and $\Gamma <L$ a lattice.
The following lemma is a routine verification using the Definition~\ref{def:3.4} 
of maximality and the Milnor--Wood type bounds in Corollary~\ref{cor:3.3}.

\begin{lemma}\label{lem:6.1} Let $\rho:\Gamma\to G$ be a homomorpism.
\be
\item If $\Gamma_0<\Gamma$ is a subgroup of finite index, then
$\rho$ is maximal if and only if $\rho|_{\Gamma_0}$ is maximal;
\item $\rho$ is maximal if and only if $q\circ\rho:\Gamma\to\gx$ is maximal;
\item $\rho$ is maximal if and only if $p_i\circ q\circ\rho:\Gamma\to\gxi$ is maximal
for every $i=1,\dots,n$. 
\ee
\end{lemma}

Recall now from \cite[Definition 2.11.]{Burger_Iozzi_Wienhard_tight} that 
if $H$ is a locally compact group, a continuous homomorphism $\rho:H\to G$
is {\it tight} if
\bqn
\|\rho^\ast(\kgb)\|=\|\kgb\|\,.
\eqn

\begin{lemma}\label{lem:6.2} If $\rho:\Gamma\to G$ is maximal, then it is tight.
\end{lemma}

\begin{proof} Using that
\bqn
\Tb\big(\rho^\ast(\kgb)\big)=\tb(\rho)\klb
\eqn
and that, because of maximality,
\bqn\tb(\rho)=\rk_\Xx=2\|\kgb\|
\eqn
we get that 
\bqn
\Tb\big(\rho^\ast(\kgb)\big)=2\|\kgb\|\klb 
\eqn
which, together with
%\bqn
$\big|\Tb\big(\rho^\ast(\kgb)\big)\big|\leq\|\rho^\ast(\kgb)\|$ 
%\eqn
and $\|\klb\|=1/2$ implies that
\bqn
\|\kgb\|\leq\|\rho^\ast(\kgb)\|\,. 
\eqn
Since the reverse inequality holds always true, we have proved the lemma.
\end{proof}

Let $H$ be of type (RH), $\Yy$ the associated symmetric space, $\Yy=\Yy_1\times,\dots,\times\Yy_m$
the decomposition into irreducible factors and $\{\khib:\,1\leq i\leq m\}$ the basis
of $\hcb^2(H,\RR)$  obtained as above.  Given a continuous homomorphisms $\sigma:H\to G$
and writing 
\bqn
\sigma^\ast(\kgb)=\sum_{i=1}^m\lambda_i\khib\,.
\eqn
It follows from (\ref{eq:norm}) and (\ref{eq:gromov-norm}) that $\sigma$ is tight if and only if
\bqn
\sum_{i=1}^m|\lambda_i|\rk_{\Yy_i}=\rk_\Xx\,.
\eqn
Finally we recall that a continuous homomorphism $\sigma:H\to G$ is {\it positive}
if $\lambda_i\geq 0$ for $1\leq i\leq m$.

\begin{lemma}\label{lem:6.3} Let $\rho:\Gamma\to H$ and $\sigma: H\to G$ be homomorphisms,
where $\sigma$ is continuous and $H,G$ are of type (RH).
\be
\item If $\sigma\circ\rho$ is maximal then $\sigma$ is tight;
\item if $\sigma\circ\rho$ is maximal and $\sigma$ is positive then $\rho$ is maximal;
\item if $\rho$ is maximal and $\sigma$ is tight and positive, then $\sigma\circ\rho$ is maximal.
\ee
\end{lemma}

\begin{proof} With the notation introduced above, let
\bqn
\sigma^\ast(\kgb)=\sum_{i=1}^m\lambda_i\khib\,.
\eqn
Thus
\be
\item[(\ref{lem:6.3}.a)] $\tb(\sigma\circ\rho)=\sum_{i=1}^m\lambda_i\tb(\rho,\khib)$;
\item[(\ref{lem:6.3}.b)] $|\tb(\rho,\khib)|\leq\r_{\Yy_i}$;
\item[(\ref{lem:6.3}.c)] $\sum_{i=1}^m|\lambda_i|\rk_{\Yy_i}\leq\rk_\Xx$.
\ee
Thus if $\sigma\circ\rho$ is maximal, the equality (\ref{lem:6.3}.a) combined with 
(\ref{lem:6.3}.b) and (\ref{lem:6.3}.c) implies that we have equality in (\ref{lem:6.3}.c)
and hence $\sigma$ is tight.

If $\sigma$ is positive, that is $\lambda_i\geq 0$ for $1\leq i\leq m$,
and $\sigma\circ\rho$ is maximal, we get from (\ref{lem:6.3}.a), (\ref{lem:6.3}.b) and (\ref{lem:6.3}.c)
that
\bqn
\tb(\rho,\khib)=\rk_{\Yy_i}
\eqn
which, together with Lemma~\ref{lem:6.1}(3) implies that $\rho$ is maximal.

Finally, if $\rho$ is maximal we get from Lemma~\ref{lem:6.1}(3)
that $\tb(\rho,\khib)=\rk_{\Yy_i}$ for $1\leq i\leq m$,
and if $\sigma$ is tight and positive then 
\bqn
 \tb(\sigma\circ\rho)
=\sum_{i=1}^m\lambda_i\tb(\rho,\khib)
=\sum_{i=1}^m\lambda_i\rk_{\Yy_i}
=\rk_\Xx
\eqn
and hence $\sigma\circ\rho$ is maximal.
\end{proof}

\begin{proof}[{\it Proofs of Theorems~\ref{thm_intro:thm3} and ~\ref{thm_intro:thm5}}]
In order to prove these results we place ourselves, as we may,
in the slightly more general context in which $\Gamma$ is a lattice in a finite covering
of $\PU(1,1)$.  Let now $\gG$ and $G=\gG(\RR)^\circ$ be as in the statement of 
Theorem~\ref{thm_intro:thm3}, that is $\gG$ is a connected
semisimple algebraic group defined over $\RR$ such that $G=\gG(\RR)^\circ$ is of Hermitian type, 
and let $\rho:\Gamma\to G$ be a maximal representation.
Set $\hH:=\overline{\rho(\Gamma)}^Z$.  Since $\rho$ is maximal, it is in particular
a tight homomorphism (Lemma~\ref{lem:6.2}) and hence \cite[Theorem 4]{Burger_Iozzi_Wienhard_tight}
applies. In particular, $\hH$ is reductive, $H:=\hH(\RR)^\circ$ has compact centralizer
in $G$ and is of type (RH); furthermore if $\Yy$ denotes the symmetric space associated
to $H$ then there is a unique $H$-invariant complex structure on $\Yy$ such that the inclusion
$i: H\to G$ is tight and positive.

Setting $\Gamma_0=\rho^{-1}(\Gamma\cap H)$ and $\rho_0:=\rho|_{\Gamma_0}:\Gamma_0\to H$,
we have from Lemma~\ref{lem:6.1}(1) that $i\circ\rho_0:\Gamma_0\to H$ is maximal and,
since $i$ is tight and positive, from Lemma~\ref{lem:6.3}(2) that $\rho_0:\Gamma_0\to H$ is maximal as well.
Composing $\rho_0$ with $p_i\circ q:H\to H_\Yy\to H_{\Yy_i}$, where $\Yy_i$, $1\leq i \leq m$ are the irreducible factors of $\Yy$,   
the resulting homomorphisms $\rho_{0,i}:\Gamma_0\to H_{\Yy_i}$ are maximal
with Zariski dense image.  Theorem~\ref{thm:4.1} then implies that $\Yy_i$ 
is of tube type and $\rho_{0,i}$ is injective, modulo the center of $\Gamma_0$, and with discrete image.
This implies that $\rho:\Gamma\to G$ is injective (modulo the center) and with discrete image.
Since $\Yy$ is of tube type and $i:H\to G$ is tight, there is a unique maximal subdomain
$\Tt\subset\Xx$ of tube type with $i(\Yy)\subset\Tt$ (see \cite[Theorem~10(1)]{Burger_Iozzi_Wienhard_tight}); 
it is moreover $H$-invariant and 
hence (by uniqueness) $\hH(\RR)$-invariant and thus $\rho(\G)$-invariant.
This completes the proof of Theorem~\ref{thm_intro:thm3}.

Applying Theorem~\ref{thm:5.1} to every irreducible factor of $\Yy$,
we get, say, a left continuous strictly $\rho_0$-equivariant map $\varphi:\partial\DD\to\cs_\Yy = \cs_{\Yy_1} \times \cdots \times \cs_{\Yy_m}$.
Since $i:H\to G$ is tight, we have also a canonical $i$-equivariant map
$\hat i:\cs_\Yy\to\cs_\Xx$; applying judiciously the uniqueness property in \cite[Theorem~4.1.]{Burger_Iozzi_Wienhard_tight},
we deduce that $\hat i\circ\varphi:\partial\DD\to\cs_\Xx$ is $\rho$-equivariant.
Finally, writing
\bqn
i^\ast(\kxb)=\sum_{i=1}^m\lambda_i\kappa_{\Yy,i}^\mathrm{b}\,,
\eqn
for $\lambda_i\geq 0$, we have (see \cite[Lemma~5.9.]{Burger_Iozzi_Wienhard_tight})
\bqn
 \beta_{\cs_\Xx}\big(\hat i(x),\hat i(y),\hat i(z)\big)
=\sum_{j=1}^m\lambda_j\beta_{\cs_{\Yy_j}}(x_j,y_j,z_j), 
\eqn
where $x = (x_1, \cdots, x_m),\, y = (y_1, \cdots, y_m),\, z=(z_1, \cdots, z_m) \in \cs_{\Yy}$. 
This, together with (\ref{eq:bb}) in Corollary~\ref{cor:4.4}, implies
\bqn
 \beta_{\cs_\Xx}\big(\hat i\varphi(a),\hat i\varphi(b),\hat i\varphi(c)\big)
=\left(\sum_{i=1}^m\lambda_i\rk_{\Yy_i}\right)\beta_{\partial\DD}(a,b,c)
\eqn
and concludes the proof since $i$ is tight and positive.
\end{proof}

%%% Local Variables: 
%%% mode: latex
%%% TeX-master: "toledo"
%%% End: 

%% file: rotation.tex
\section{Rotation numbers and applications to groups of Hermitian type}\label{sec:rotation}
In this section we introduce and study rotation numbers on locally compact groups
and compute them for groups of Hermitian type.  The results are of independent interest. Here they 
are used in an essential way
in the computation of the Toledo invariant in the case of surfaces with boundary.

\subsection{Basic definitions and properties}\label{subsec:7.1}
Let $G$ be a locally compact group and $\kappa\in\hhcb^2(G,\ZZ)$ a bounded integer valued Borel class.
Let $B<G$ be a closed subgroup and consider the first few terms of the long exact sequence
\bq\label{eq:seq}
\quad\quad
\xymatrix@1{
 0\ar[r]
&\homc(B,\RR/\ZZ)\ar[r]^{\quad\delta}
&\hhcb^2(B,\ZZ)\ar[r]
&\hcb^2(B,\RR)\ar[r]
&\dots}
\eq
coming from the coefficient sequence \eqref{eq:coeff}; 
denote by $\kappa_\RR$ the image in $\hcb^2(G,\RR) $ of an element $\kappa\in\hhcb^2(G,\ZZ)$.
Then if $\kappa_\RR|_B=0$, we let $f_B:B\to\RR/\ZZ$
denote the unique continuous homomorphism with $\delta(f_B)=\kappa|_B$.  
In particular this applies to $B=\overline{\langle g\rangle}$ for any $g\in G$ and we define
\begin{defi}\label{defi:rot_number}
The rotation number of $g\in G$  with respect to $\kappa\in\hhcb^2(G,\ZZ)$ is  
\bqn
\rotk(g):=f_{\overline{\langle g\rangle}} (g)\,.
\eqn
\end{defi}
Using that the exact sequence in \eqref{eq:seq} is natural with respect to group homomorphisms,
one verifies easily the following properties

\begin{lemma}\label{lem:7.1} 
\be
\item $\rotk:G\to\RR/\ZZ$ is invariant under conjugation;
\item if $\kappa_\RR|_B=0$, then $\rotk|_B$ is a continuous homomorphism and $\delta(\rotk|_B)=\kappa|_B$;
\item if $\sigma:G_1\to G_2$ is a continuous homomorphism and $\kappa_1=\sigma^\ast(\kappa_2)$,
then 
\bqn
\rot_{\kappa_1}(g_1)=\rot_{\kappa_2}\big(\sigma(g_1)\big)
\eqn
for all $g_1\in G_1$.
\ee
\end{lemma}

In the study of rotation numbers $\rotk$, quasimorphisms play an important role.  
We quickly review the basic definitions.  
If $A=\ZZ,\RR$, a function $f:G\to A$ is a quasimorphism if the function $df:G\to A$
\bqn
df(x,y)=f(xy)-f(x)-f(y)
\eqn
is bounded. When $A=\RR$, a quasimorphism is homogeneous if
\bqn
f(g^n)=nf(g)
\eqn
for $n\in\ZZ$ and $g\in G$.  Any quasimorphism $f:G\to A$ can be made homogeneous
by setting 
\bqn
Hf(x):=\lim_{n\to\infty}\frac{f(x^n)}{n}\in\RR
\eqn
and it is a standard fact that $f-Hf$ is bounded.

\begin{lemma}\label{lem:7.2} Assume that $\kappa$ vanishes when considered as an ordinary class
in $\hhc^2(G,\ZZ)$.  Let $f:G\to\ZZ$ be a Borel map such that $df$ is bounded and represents $\kappa$
seen as a bounded class.  Then $f$ is a quasimorphism and 
\bqn
\rotk(g)\equiv Hf(g)\mod\ZZ
\eqn
where $Hf$ is the homogenization of $f$.
\end{lemma}

\begin{proof} By Lemma~\ref{lem:7.1}(3), it suffices to show the assertion for $G=\ZZ$.
Let $p:\RR\to\RR/\ZZ$ be the canonical projection; then $p\circ Hf:\ZZ\to\RR/\ZZ$
is a homomorphism since $Hf$ is homogeneous, and we claim that $\kappa=\delta(p\circ Hf)$.
Indeed let $\sigma:\RR/\ZZ\to\RR$ be the Borel section of $p$ with values in $[0,1)$;
then we have
\bqn
 \delta(p\circ HF)
=-\d(\sigma\circ p \circ Hf)
=df-d((f-Hf+\sigma\circ p\circ Hf)
\eqn
where we have used in the last equality that $dHf=0$.  
The claim then follows from the fact that $f$ and $Hf-\sigma\circ p\circ Hf$ take integral values
and both $f-Hf$ and $\sigma\circ p\circ Hf$ are bounded.  
From the claim we the get that 
\bqn
\rotk(g)=p\circ Hf(g)
\eqn
for all $g\in \ZZ$, which is the assertion that needed to be proved.
\end{proof}

Next we have:

\begin{lemma}\label{lem:7.3} Let $f:G\to\RR$ be a homogeneous Borel quasimorphism.
Then $f$ is continuous.
\end{lemma}

\begin{proof} We show first that $f$ is locally bounded.  
Let $C:=\sup_{x,y}|df(x,y)|$ and $E_N=\big\{x\in G:\,|f(x)|\leq N\big\}$.
Then $E_N=E_N^{-1}$ and for $N$ large enough is of positive Haar measure. 
For such $N$ we deduce that $E_N\cdot E_N$ is a neighborhood of $e\in G$;
since $f$ is a quasimorphism we have that $E_N\cdot E_N\subset E_{N+2C}$
which implies that $f$ is bounded in a neighborhood of $e$ and hence,
by the quasimorphism property, locally bounded.
Fix now $\varphi:G\to[0,\infty)$ continuous with compact support and 
of total integral one.  Since $f$ is locally bounded,
we have that for every $n\in\NN$
\bqn
F_n(x)=\frac{1}{n}\int_G\big(f(x^ny)-f(y)\big)\varphi(y)dy
\eqn
is defined and continuous.  Since $f$ is homogeneous we have
\bqn
 \big|f(x)-F_n(x)\big|
&=&\left|\frac{1}{n}f(x^n)-F_n(x)\right|\\
&=&\left|\frac{1}{n}\int_G\big(f(x^n)+f(y)-f(x^ny)\big)\varphi(y)dy\right|
\leq\frac{C}{n}
\eqn
which implies that $f$ is the uniform limit of a sequence of continuous functions 
and therefore continuous.
\end{proof}

Now we come to the main goal of this subsection, which is the continuity of $\rotk$;
this is shown by exhibiting a direct relationship with a certain quasimorphism 
on a central extension of $G$.  More precisely let as before $\kappa\in\hhcb^2(G,\ZZ)$;
then $\kappa$ can be seen as a class in $\hc^2(G,\ZZ)$
and hence gives rise by \cite{Mackey_57} to a topological central extension
\bqn
\xymatrix@1{
 0\ar[r]
&\ZZ\ar[r]^i
&G_\kappa\ar[r]^p
&G\ar[r]
&e
}
\eqn
that is $G_\kappa$ is a locally compact group, $i$ and $p$ are continuous, 
$i(\ZZ)$ is a closed central subgroup of $G_\kappa$
and $G_\kappa/i(\ZZ)$ is topologically isomorphic to $G$.

\begin{prop}\label{prop:7.4} There exists a continuous homogeneous quasimorphism $f:G_\kappa\to\RR$
such that 
\be
\item $f\big(i(n)g\big)=n+f(g)$, for $n\in\ZZ$, $g\in G_\kappa$;
\item $\rotk\big(p(g)\big)\equiv f(g)\mod\ZZ$.
\ee
\end{prop}

\begin{cor}\label{cor:7.5} Let $\kappa\in\hhcb^2(G,\ZZ)$.  Then $\rotk:G\to\RR/\ZZ$ is continuous.
\end{cor}

\begin{proof}[{\it Proof of Proposition~\ref{prop:7.4}}] Let $c:G^2\to\ZZ$ be a bounded Borel
cocycle which represents $\kappa\in \hhcb^2(G,\ZZ)$ and which we assume to be normalized.
Then $G_\kappa$ is a Borel group isomorphic to the Borel space $G\times\ZZ$ with multiplication given by 
\bqn
(g_1,n_1)(g_2,n_2)=\big(g_1g_2,n_1+n_2+c(g_1,g_2)\big)\,.
\eqn
Define $f_1:G_\kappa\to\ZZ$ by $f_1(g,m):=m$.  Then $f_1$ is a Borel function
and $df_1$ is a bounded Borel cocycle representing $p^\ast(\kappa)\in\hhcb^2(G_\kappa,\ZZ)$.
let $f:G_\kappa\to\RR$ be the homogenization of $f_1$.  
Then Lemma~\ref{lem:7.2} implies that 
\bqn
\rotk\big(p(g)\big) = \rot_{p^\ast(\kappa)}(g) \equiv f(g) \mod\ZZ
%f(g)\cong\rot_{p^\ast(\kappa)}(g)=\rotk\big(p(g)\big)
\eqn
and Lemma~\ref{lem:7.3} that $f$ is continuous.  
Finally $f$ satisfies (1) because $f_1$ does and $i(\ZZ)$ is central.
\end{proof}

\subsection{Rotation numbers on groups of Hermitian type}\label{subsec:7.2}
We begin first by determining $\hhcb^2(G,\ZZ)$ for $G$ of Hermitian type.
The main points are summarized in the following

\begin{prop}\label{prop:7.6}  Let $G$ be a group of Hermitian type and $K<G$ a maximal compact subgroup.
\be
\item The comparison map $\hhcb^2(G,\ZZ)\to\hhc^2(G,\ZZ)$ is an isomorphism;
\item the map 
\bq
\ba
\hhcb^2(G,\ZZ)&\to\homc(K,\RR/\ZZ)\\
\kappa&\mapsto\rot_{\kappa|_K}
\ea
\eq
is an isomorphism;
\item the change of coefficient map
\bqn
\hhcb^2(G,\ZZ)\to\hcb^2(G,\RR)
\eqn
is injective and its image is a lattice.
\ee
\end{prop}

\begin{proof} (1) follows from the commutativity of the diagram
\bqn
\xymatrix{
 0\ar[r]
&\hhcb^2(G,\ZZ)\ar[r]\ar[d]
&\hhcb^2(G,\RR)\ar[r]\ar[d]
&\hcb^2(G,\RR/\ZZ)\ar@{=}[d]\\
 0\ar[r]
&\hhc^2(G,\ZZ)\ar[r]
&\hhc^2(G,\RR)\ar[r]
&\hc^2(G,\RR/\ZZ)
}
\eqn
and the fact that with real coefficients the comparison map is an isomorphism.

\medskip
\noindent
(2) follows from (1) and the fact that in ordinary Borel cohomology
the restriction 
\bqn
\hhc^2(G,\ZZ)\to\hhc^2(K,\ZZ)
\eqn
is an isomorphism \cite{Wigner_73}.

\medskip
\noindent
(3) follows from (1) and the corresponding statement in ordinary cohomology.
\end{proof}

In view of the preceding proposition we can refer for every $u\in\homc(K,\RR/\ZZ)$ to the class
$\kappa$ associated to $u$, and conversely.  

We turn now  to the explicit computation of the rotation number function.
let $G=KAN$ be an Iwasawa decomposition;  recall the refined Jordan decomposition,
namely that every $g\in G$ is a product $g=g_eg_hg_n$, where $g_e$ is contained
in a compact subgroup, $g_h$ and $g_n$ are conjugated to an element respectively in $A$ and $N$;
moreover the elements $g_e,g_h$ and $g_n$ commute pairwise.  Then we have:

\begin{prop}\label{prop:7.7}  Let $\kappa\in\hhcb^2(G,\ZZ)$ and $u\in\homc(K,\RR/\ZZ)$
be the corresponding homomorphism.
\be
\item $\rotk|_{AN}$ is the trivial homomorphism;
\item for $g\in G$, let $g_e$ be the elliptic component in the refined Jordan decomposition of $g$ and
$k\in C(g_e)\cap K$, where $C(g_e)$ denotes the $G$-conjugacy class of $g_e$.  Then $\rotk(g)=u(k)$.
\ee
\end{prop}

\begin{proof}  (1) Since $B=AN$ is amenable, $\kappa_\RR|_B=0$ and thus
Lemma~\ref{lem:7.1}(2) implies that $\rotk|_{AN}$ is a continuous homomorphism
and hence differentiable.  Since $[B,B]=N$, then $\rotk(N)=0$.
The restriction $\rotk|_A$ is invariant under the Weyl group $\Nn_K(A)/\Zz_K(A)$ and hence its differential 
\bqn
D_e\rotk|_A:\fraka\to\RR
\eqn
is a linear form invariant under the Weyl group; it must therefore vanish since $G$ is semisimple
and thus $\rotk|_A$ is trivial as well.

\medskip
\noindent
(2) Let $g=g_eg_hg_n$ be the refined Jordan decomposition of $g$.  Since the subgroup $C$
generated by $g_e, g_h$ and $g_n$ is Abelian, $\rotk|_C$ is a homomorphism,
hence 
\bqn
\rotk(g)=\rotk(g_e)+\rotk(g_h)+\rotk(g_n)\,.
\eqn
Taking into account (1) and the fact that $g_h$ and $g_n$ are conjugate respectively to  elements 
in $A$ and $N$, we get that
\bqn
\rotk(g)=\rotk(g_e)=\rotk(k)=u(k)
\eqn
which concludes the proof.
\end{proof}

For the next result, if $u\in\homc(K,\RR/\ZZ)$, let us denote by 
$u_\ast:\pi_1(G)\to\ZZ$ the homomorphism induced by $u$ on the level of fundamental 
groups, where we have identified $\pi_1(K)$ with $\pi_1(G)$.

\begin{thm}\label{thm:7.8} Let $u:K\to\RR/\ZZ$ be a continuous homomorphism and $\kappa\in\hhcb^2(G,\ZZ)$
the associated class.  
\be
\item $\rotk$ is continuous;
\item the unique continuous lift $\widetilde\rotk:\widetilde G\to\RR$ such that 
$\widetilde\rotk(e)=0$ is a continuous homogeneous quasimorphism and satisfies
\bqn
\widetilde\rotk(zg)=u_\ast(z)+\widetilde\rotk(g)
\eqn
for all $z\in\pi_1(G)$ and for all $g\in\widetilde G$.
\ee
\end{thm}

As a consequence we obtain a description of the space $\Qq(\widetilde G)_\ZZ$
of continuous homogeneous quasimorphism $f:\widetilde G\to\RR$ such that 
$f\big(\pi_1(G)\big)\subset\ZZ$.

\begin{cor}\label{cor:7.9} The maps
\bqn
\ba
\homc(K,\RR/\ZZ)&\to     \Qq(\widetilde G)_\ZZ  \to     \hom\big(\pi_1(G),\ZZ\big)\\
  u \qquad  &\mapsto\, \, \widetilde\rot_\kappa\, \, \, \mapsto\,\,\widetilde{\rot_\kappa}|_{\pi_1(G)}=u_\ast
\ea
\eqn
are group isomorphisms.
\end{cor}

\begin{proof}[{\it Proof of Theorem~\ref{thm:7.8}}] The first assertion is a special case of Corollary~\ref{cor:7.5}.
Let then $p:G_\kappa\to G$ be the Lie group central extension determined by $\kappa$ and 
$\pi:\widetilde G\to(G_\kappa)^\circ$ the canonical projection.  
If $f:G_\kappa\to\RR$ is the continuous homogeneous quasimorphism given by Proposition~\ref{prop:7.4},
then it follow from Proposition~\ref{prop:7.4}(2) that $f|_{G_\kappa}\circ\pi:\widetilde G\to\RR$
is a continuous lift of $\rotk$ to $\widetilde G$ which moreover vanishes at $e$.  
Hence $\widetilde\rotk=f|_{(G_\kappa)^\circ}\circ\pi$, 
which implies the remaining assertion in the theorem.
\end{proof}

\begin{proof}[{\it Proof of Corollary~\ref{cor:7.9}}] Since the composition of the two arrows 
is the isomorphism 
\bqn
\ba
\homc(K,\RR/\ZZ)&\to      \hom\big(\pi_1(G),\ZZ\big)\\
   u \qquad  &\mapsto\quad\qquad u_\ast
\ea
\eqn
it suffices to show that the second morphism is injective.
If $f_1,f_2:\widetilde G\to\RR$ are homogeneous continuous quasimorphisms
which induce the same homomorphism $h:\pi_1(G)\to\ZZ$
then their difference $f_1-f_2$ is $\pi_1(G)$-invariant and hence descends to a homogeneous
quasimorphism $G\to\RR$ which therefore vanishes since $G$ is connected semisimple
with finite center.  Thus $f_1=f_2$, which completes the proof.
\end{proof}

\begin{rem}\label{rem:rot_numbers}
For special classes 
$\kappa \in \hhcb^2(G,\ZZ)$ the rotation number $\rot_\kappa$ coincides with previously known constructions:
\begin{enumerate}
\item If $G=\Homeo_+(S^1)$ is the group of orientation preserving
homeomorphisms of the circle (viewed as an abstract group) and
$e^\mathrm{b}\in\hb^2(G,\ZZ)$ is the bounded Euler class, 
Ghys \cite{Ghys_87} observed that $\rot_{e^\mathrm{b}}(\varphi)$ 
is the classical rotation number of the homeomorphism $\varphi$. 

\item To obtain the symplectic rotation number defined by Barge and Ghys \cite{Barge_Ghys}
for $G=\mathrm{Sp}(2n,\RR)$, we have to consider the class $\kappa$ which corresponds to the homomorphism
$u:K= \mathrm{U}(n)\to\TT$ defined by $u(k)=(\det k)^2$.  

\item If $\Dd$ is an irreducible symmetric domain of tube type, $G=\aut(\Dd)^\circ$
and $K$ is the stabilizer of $0\in\Dd$, Clerc and Koufany construct
a homomorphism $\chi:K\to\TT$ using the Jordan algebra determinant.
The rotation number function and the quasimorphism constructed in their paper 
\cite[Theorem~10.3 and Proposition~10.4]{Clerc_Koufany},
coincide then respectively with $\rot_\kappa$ and $\widetilde{\rot_\kappa}$, where $\kappa$ is the class corresponding 
to $\chi$. 

We observe moreover that if $u:K\to\TT$ is the complex Jacobian at $0$,
then for every $k\in K$ we have that 
\bqn
\chi(k)^{p_\Xx}=u(k)^2\,,
\eqn
which incidentally shows that $\kgb$ is in the image of 
$\hhcb^2(G,\ZZ)$.
\end{enumerate}
\end{rem}

%%% Local Variables: 
%%% mode: latex
%%% TeX-master: "toledo"
%%% End: 

%% file: formula.tex
\section{Toledo numbers: formula and applications to representation varieties}\label{sec:formula}

\subsection{The formula}\label{subsec:8.1}
Let $\Sigma$ be a connected oriented surface and $G$ a group of Hermitian type.
In this subsection we establish a formula for the Toledo invariant $\tksr$
where $\kappa$ is a bounded integral class, and we concentrate on the case
in which $\bs\neq\emptyset$; we mention at the end the formula in the case in which
$\bs=\emptyset$.

The boundary of $\Sigma$ is the union $\bs=\bigsqcup_{j=1}^n C_j$
of oriented circles and we fix a presentation
\bqn
 \pi_1(\Sigma)
=\bigg\<a_1,b_1,\dots,a_g,b_g,c_1,\dots,c_n:\,
\prod_{i=1}^g[a_i,b_i]\prod_{j=1}^nc_j=e\bigg>
\eqn
where $g$ is the genus of $\Sigma$ and $c_j$ is freely homotopic to $C_j$ 
with positive orientation. Combining now the long exact sequence 
in bounded cohomology associated to the pair of spaces $(\Sigma,\partial\Sigma)$
and the one associated to the usual coefficient sequence in \eqref{eq:coeff},
we obtain
\bqn
\xymatrix{
 0\ar[r]
&\hb^2(\Sigma,\bs,\RR)\ar[r]^-{j_\bs}
&\hb^2(\Sigma,\RR)\ar[r]
&\hb^2(\bs,\RR)=0\\
& 
&\hb^2(\Sigma,\ZZ)\ar[u]\ar[r]
&\hb^2(\bs,\ZZ)\ar[u]\\
&
&
&\h^1(\bs,\RR/\ZZ)\ar[u]^\delta\\
&
& 
&0\ar[u]
}
\eqn
where we have used that $\hb^i(\bs,\RR)=0$, $i\geq 1$.
We have then the following congruence relation:

\begin{lemma}\label{lem:8.1} Let $\alpha\in\hb^2(\Sigma,\ZZ)$ and denote by
$\alpha_\RR$ its image in $\hb^2(\Sigma,\RR)$.
\bq\label{eq:8.1}
 \ba
       \big\< j_{\partial \Sigma}^{-1} (\alpha_\RR), [\Sigma, \partial\Sigma]\big\>
\equiv &-\big\<\delta^{-1}(\alpha|_{\partial \Sigma}), [\partial \Sigma] \big\>\mod\ZZ\\
      =&-\sum_{j=1}^n\big\<\delta^{-1}\alpha|_{C_j},[C_j]\big\>\,\,, 
\ea
\eq
where we view $j_{\partial \Sigma}^{-1} (\alpha_\RR)$ as ordinary relative singular cohomology class.
\end{lemma}

\begin{proof}
Let $c\in\Zz_\mathrm{b}^2(\Sigma,\ZZ)$ be a $\ZZ$-valued bounded cocycle
representing $\alpha\in\hb^2(\Sigma,\ZZ)$, 
and let $c|_{\partial\Sigma}\in\Zz_\mathrm{b}^2(\partial\Sigma,\ZZ)$ 
be its restriction to the boundary $\partial\Sigma$.  
Since the fundamental groups of the components of $\partial\Sigma$ are amenable,
there exists a bounded $\RR$-valued $1$-cochain $c'\in F_\mathrm{b}^1(\partial\Sigma,\RR)$
such that $dc'=c|_{\partial\Sigma}$.  Let $c''\in F_\mathrm{b}^1(\partial\Sigma,\RR/\ZZ)$
be the corresponding $\RR/\ZZ$-valued $1$-cochain on $\partial\Sigma$:  
for any $1$-simplex $t\in S_1(\partial\Sigma)$ 
we have that 
\bq\label{eq:c'c''}
\<c',t\>\equiv\<c'',t\>\mod\ZZ\,
\eq
and moreover, since $c|_{\partial\Sigma}$ is $\ZZ$-valued, 
$c''$ is a $1$-cocycle which represents the class 
$\delta^{-1}(\alpha|_{\partial\Sigma})\in\h^1(\partial\Sigma,\RR/\ZZ)$.

On the other hand, we can extend $c'$ to a $1$-cochain $\widetilde{c'}$ on $\Sigma$
by setting 
\bqn
\widetilde{c'}(\sigma)=
\begin{cases}
c'(\sigma)&\hbox{if } \sigma\in S_1(\partial\Sigma)\\
\quad 0         &\hbox{otherwise}\,,
\end{cases}
\eqn
so that $c-d\widetilde{c'}\in\Zz_\mathrm{b}^2(\Sigma,\partial\Sigma,\RR)$ 
is a cocycle which represents $j_{\partial\Sigma}^{-1}(\alpha_\RR)$.

Let now $s$ be a two-chain which represents the relative fundamental class 
$[\Sigma,\partial\Sigma]$, 
so that $\partial s$ represents the fundamental class $[\partial\Sigma]$. 
From the definition of $\widetilde{c'}$, from (\ref{eq:c'c''}) and
the fact that $\<c,s\>\in\ZZ$, it follows that 
\bqn
\ba
      \<c-d\widetilde{c'},s\>
     =\<c,s\>-\<d\widetilde{c'},s\>
%\equiv-\<d\widetilde{c'},s\>\mod\ZZ
%     =-\<\widetilde{c'},\partial s\>\mod\ZZ\\
%     &=-\<c',\partial s\>\mod\ZZ
\equiv-\<c'',\partial s\>\mod\ZZ\,,
\ea
\eqn
thus completing the proof.
\end{proof}

We apply the above general lemma to the situation at hand and show the following

\begin{lemma}\label{lem:8.2} Let $\kappa\in\hhcb^2(G,\ZZ)$ and
$\rho:\ps\to G$ a homomorphism.  Then
\bqn
\tksr\equiv-\sum_{j=1}^n\rotk\big(\rho(c_j)\big)\mod\ZZ\,.
\eqn
\end{lemma}

\begin{proof} Consider the commutative diagram
\bqn
\xymatrix{
 \hb^2\big(\ps,\ZZ\big)\ar[r]^{g_\Sigma}\ar[d]
&\hb^2(\Sigma,\ZZ)\ar[d]\\
 \hb^2\big(\pi_1(C_j),\ZZ\big)\ar[r]
&\hb^2(C_j,\ZZ)\\
 \h^1\big(\pi_1(C_j),\RR/\ZZ\big)\ar[u]_\delta\ar[r]
&\h^1(C_j,\RR/\ZZ)\ar[u]
}
\eqn
where the horizontal arrows are isomorphisms, the vertical arrows between  
first and second row are restriction maps and those between third and second are connecting homomorphisms.
The commutativity implies the first equality
\bqn
 \big\langle\delta^{-1}\big([g_{\Sigma} \rho^\ast(\kappa)]|_{C_j}\big),[C_j]\big\rangle
=\delta^{-1}\big(\rho^\ast(\kappa)|_{\pi_1(C_j)}\big)(c_j)
=\rotk\big(\rho(c_j)\big)\,,
\eqn
while the second is the definition of $\rotk$ (see Definition~\ref{defi:rot_number}).  This, together with
Lemma~\ref{lem:8.1} applied to $\alpha=g_\Sigma \big(\rho^\ast(\kappa)\big)$
implies the result.
\end{proof}

Now we come to the formula for the Toledo invariant.
Observe first that when $\bs\neq\emptyset$, $\ps$ is a free group
and thus any homomorphism $\rho:\ps\to G$ admits a lift $\widetilde\rho:\ps\to \widetilde G$.

\begin{proof}[{\it Proof of Theorem~\ref{thm_intro:thm8}}] 
Using the equivariance property of $\widetilde\rotk$ in Theorem~\ref{thm:7.8}(2),
one checks that 
\bqn
\mathrm{R}(\rho):=-\sum_{j=1}^n\widetilde\rotk\big(\widetilde\rho(c_j)\big)
\eqn
does not depend on the choice of the lift $\widetilde\rho$.
Thus the map 
\bqn
\ba
\hom\big(\ps,G\big)&\to\RR\\
\rho\qquad&\mapsto\mathrm{R}(\rho)
\ea
\eqn
is well defined and continuous since $\widetilde\rotk$ is continuous
and $\widetilde G\to G$ is a covering.  This implies with Proposition~\ref{prop:3.8}
that the map 
\bq\label{eq:8}
\rho\mapsto\tksr-\mathrm{R}(\rho)
\eq 
is continuous;
on the other hand (from Lemma~\ref{lem:8.2}) we know that this map is $\ZZ$-valued 
and hence, since $\hom\big(\ps,G\big)$ is connected, that \eqref{eq:8}
is constant: evaluation at the trivial homomorphism implies   
that this constant is zero thus showing the theorem.
\end{proof}

Finally, we indicate briefly the formula when $\bs=\emptyset$.  Let
\bqn
[\,\cdot\,,\,\cdot\,]\tilde{\ }:G\times G\to\widetilde G
\eqn
denote the $\widetilde G$-valued commutator map.  Recall that when 
$\bs=\emptyset$ the consideration of ordinary cohomology suffices.

\begin{thm}\label{thm:8.4} Let $\kappa\in\hhc^2(G,\ZZ)$ and $u\in\homc(K,\RR/\ZZ)$
the associated homomorphism.  If $\rho:\ps\to G$ is a representation, we have
\bqn
\tksr=-u_\ast\left(\prod_{j=1}^g\big[\rho(a_i),\rho(b_i)\big]\tilde{\ }\right)
\eqn
where, as usual, $u_\ast:\pi_1(K)=\pi_1(G)\to \ZZ$ is the homomorphism 
induced by $u$.
\end{thm}

\begin{rem}\label{rem:8.5} One may prove the formula in Theorem~\ref{thm:8.4}
by cutting $\Sigma$ along the separating curve $[a_1,b_1]$ and combine
the formula in Theorem~\ref{thm:8.4} applied to each component
together with the additivity property in Proposition~\ref{prop:3.1}(1).
\end{rem}

\begin{rem}
The formula in Theorem~\ref{thm:8.4} generalizes Milnor's classical formula for the Euler number of a representation into $\mathrm{GL}^+(2)$ \cite{Milnor}. 
\end{rem}
\subsection{The bounded fundamental class and generalized $w_1$-classes}\label{subsec:8.2}
The group $\PU(1,1)$ acts effectively on the circle $\partial\DD$
and we have seen that if $\rho_1$ and $\rho_2$ are maximal representations
of $\ps$ into $\PU(1,1)$ then the resulting actions on $\partial\DD$
are semiconjugate in the sense of Ghys \cite{Ghys_87}. 
If $e^\mathrm{b}\in\hb^2\big(\PU(1,1),\ZZ\big)$ denotes the bounded
Euler class, or more precisely the restriction to $\PU(1,1)$ of the bounded
Euler class of the group of orientation preserving homeomorphisms of $\partial\DD$,
then $\rho_1^\ast(e^\mathrm{b})=\rho_2^\ast(e^\mathrm{b})$ \cite{Ghys_87}.
Thus we obtain a canonical class
\bqn
\kszb\in\hb^2\big(\ps,\ZZ\big)
\eqn
associated to the oriented surface $\Sigma$ and which plays
the role of the classical fundamental class when $\bs\neq\emptyset$.
We propose to call it the {\it bounded fundamental class of $\Sigma$};
observe that even when $\bs=\emptyset$, this class contains more information
than the usual fundamental class since
\bqn
\hb^2\big(\ps,\ZZ\big)\to\h^2\big(\ps,\ZZ\big)
\eqn
is never injective.  
Thus we obtain that for a representation $\rho:\ps\to\PU(1,1)$ the following
are equivalent
\be
\item $\rho$ is maximal;
\item $\rho$ comes from a complete hyperbolic structure on $\Sigma^\circ$;
\item $\rho^\ast(e^\mathrm{b})=\kszb$.
\ee
For general groups $G$ of Hermitian type, an analogue
of the equivalence of (1) and (3) holds for real coefficients;
the extent to which it does not hold for integral coefficients 
will lead to nontrivial invariants for maximal representations.

Let now $G$ be of Hermitian type and, as usual, let 
$\{\kgib:\,1\leq i\leq n\}$ be the basis of $\hcb^2(G,\RR)$
determined by the decomposition of the associated symmetric space
$\Xx=\Xx_1\times\dots\Xx_n$ into irreducible factors.  
We define a linear form $\lambda_G$ on $\hcb^2(G,\RR)$ by
\bqn
\lambda_G(\kgib)=\rk_{\Xx_i}\,.
\eqn
Let $\ksb\in\hb^2\big(\ps,\RR\big)$ denote the real class which is the image of $\kszb$
by change of coefficients.
Then it follows from our results obtained so far that:

\begin{cor}\label{cor:8.6} For a homomorphism $\rho:\ps\to G$ the following
are equivalent:
\be
\item $\rho$ is maximal:
\item $\rho^\ast(\kappa)=\lambda_G(\kappa)\ksb$ for all $\kappa\in\hcb^2(G,\RR)$.
\ee
\end{cor}

Let now $\Dd$ be the bounded symmetric domain associated to $G$, 
$\cs$ its Shilov boundary and $Q$ the stabilizer of some point in $\cs$;
let $e_G$ be the exponent of the finite group $Q/Q^\circ$.
We will furthermore denote by $\hmaxpg$ the set of maximal representations
of $\pi_1(\Sigma)$ into $G$.

Let $\kappa\in\hhcb^2(G,\ZZ)$ and $\rho_0:\ps \to G$ a maximal representation then Theorem~\ref{thm_intro:thm9} states that for every maximal representation $\rho:\ps\to G$ the map 
\bqn 
\ba
\mathrm{R}^{\rho_0}_\kappa(\rho):\ps&\to\quad\RR/\ZZ\\
\gamma\quad&\mapsto\rotk\big(\rho(\gamma)\big) - \rotk\big(\rho_0(\gamma)\big)
\ea
\eqn
is a homomorphism, which takes values in $e_G^{-1}\ZZ/\ZZ$ if $\Dd$ is of tube type. 

Before we turn now to the proof of Theorem~\ref{thm_intro:thm9} we give an example and state a few preliminary lemmas.
\begin{exo}\label{exo:8.8}
Let $V$ be a real symplectic vector space of dimension $2n$. 
Let $K=\mathrm{U}(V,J)$ be the maximal compact subgroup of $\Sp(V)$ given by 
the choice of a compatible complex structure $J$ on $V$,
$\det:K\to\TT$ the complex determinant and let $\kappa\in\hhcb^2\big(\Sp(V),\ZZ\big)$
be the bounded integer class associated to the homomorphism $u:K\to\RR/\ZZ$,
where $e^{2\pi i u}=\det$.

When $n=1$, the associated rotation number $\rot_\kappa(\rho):\ps\to\ZZ/2\ZZ$, $\gamma \mapsto\rot_\kappa\big(\rho(\gamma)\big)$ associates to an element $\gamma \in \ps$ the sign of the 
eigenvalue of $\rho(\gamma) \in \Sp(V) \cong \SL(2,\RR)$. In particular, $\rot_\kappa(\rho)$ itself is not a homomorphism. 

On the other hand, when $n$ is even, let $\Delta:\SL(2,\RR)\to \Sp(V)$ be the homomorphism corresponding 
to a diagonal disk $\DD\to\Xx$, and choose a hyperbolization $h:\ps\to \SL(2,\RR)$. 
Setting $\rho_0=\Delta\circ h$, we have that $\rotk(\rho_0)=0$
and hence $\rotk(\rho) = \mathrm{R}^{\rho_0}_\kappa(\rho)$ is therefore a homomorphism.
This homomorphism is related to the first Stiefel-Whitney class of the following bundle. 
%\end{exo}
%
%\begin{exo}\label{exo:8.8} Let $V$ be a real symplectic vector space,
Let $\Ll(V)$ the Grassmannian of Lagrangian subspaces in $V$ and
$\Ll\to\Ll(V)$ the tautological bundle.  In \cite{Burger_Iozzi_Labourie_Wienhard}
we have shown that if $S=\Gamma\backslash\DD$ is a closed hyperbolic surface and
$\rho:\Gamma\to\Sp(V)$ is a maximal representation, the equivariant map
\bqn
\varphi:\partial\DD\to\Ll(V)
\eqn
in Theorem~\ref{thm_intro:thm5} is continuous.
Composing $\varphi$ with the visual map $T_1\DD\to\partial\DD$ and pulling back
the bundle $\Ll$, we get a vector bundle $\Ll_\rho\to T_1S$ with base 
the unit tangent bundle of $S$.  
%Let $K=\mathrm{U}(V,J)$ be the maximal compact subgroup of $\Sp(V)$ given by 
%the choice of a compatible complex structure $J$ on $V$,
%$\det:K\to\TT$ the complex determinant and let $\kappa\in\hhcb^2\big(\Sp(V),\ZZ\big)$
%be the bounded integer class associated to the homomorphism $u:K\to\RR/\ZZ$,
%where $e^{2\pi i u}=\det$.  
Then the first Stiefel--Whitney class 
$w_1(\Ll_\rho)\in\h^1(T_1S,\ZZ/2\ZZ)$ is given by the composition 
of the projection $\pi_1(T_1S)\to\Gamma$ and $\rot_\kappa(\rho):\Gamma\to\ZZ/2\ZZ$.
\end{exo}

This example shows that in some cases already $\rot_\kappa(\rho)$ is a homomorphism, whereas in general only the difference 
$\mathrm{R}^{\rho_0}_\kappa(\rho)$ for a fixed maximal representation $\rho_0$ is a homomorphism, which generalizes 
the first Stiefel--Whitney class; in fact when $\bs\neq\emptyset$, 
the map $\varphi$ is in general not continuous as the case of $\PU(1,1)$ 
already shows, and there is no (continuous) bundle in sight.

\begin{lemma}\label{lem:8.8} Let $\rho:\ps\to G$ be a maximal representation. 
Then for every $\gamma\in\ps$, $\rho(\gamma)$ has at least one fixed point in $\cs$.
\end{lemma}

\begin{proof} This follows at once from Theorem~\ref{thm_intro:thm5}, 
more specifically from the strict equivariance of the left continuous map
$\varphi:\partial\DD\to\cs$.
\end{proof}

\begin{lemma}\label{lem:8.9} The restriction map $\hcb^2(G,\RR)\to\hcb^2(Q,\RR)$
is the zero map.
\end{lemma}

\begin{proof} 
Let $\Xx = \Xx_1 \times \cdots \times \Xx_n$ be a decomposition of the symmetric space associated to $G$ into irreducible factors. Then $\cs = \cs_1 \times \cdots \times \cs_n$, where $\cs_i$ is the Shilov boundary of $\Xx_i$. 
Let $p_i: \cs \to \cs_i$ be the projection onto the $i$-th factor and set $\beta_{\cs, i} = {p_i}^\ast \beta_{\cs_i}$, where 
$\beta_{\cs_i}$ is the generalized Maslov cocycle of $\cs_i$.

Let $\kappa = \sum_{i=1}^n \lambda_i \kgib$, where $\{\kgib:\,1\leq i\leq n\}$ is the basis of $\hcb^2(G,\RR)$. 
Applying \cite[Corollary~2.3]{Burger_Iozzi_app},
we have for any $Q$-invariant Borel set $Z\subset\cs$ a commutative diagram
\bqn
\xymatrix{
 \h^2\big(\balt(\cs^\bullet)^Q\big)\ar[r]\ar[d]
&\hcb^2(Q,\RR)\\
 \h^2\big(\balt(Z^\bullet)^Q\ar[ur]\big)
&
}
\eqn
where the class $[\b:= \sum_{i=1}^n \lambda_i \b_{\cs,i}]\in\h^2\big(\balt(\cs^\bullet)^Q\big)$
goes to $\kappa|_Q$ (see \S~\ref{subsubsec:2.1.3});
taking now $Z$ to be the $Q$-fixed point in $\cs$ and observing that 
$\b|_{Z^3}=0$, we get that $\kappa|_Q=0$.
\end{proof}

\begin{lemma}\label{lem:8.10} $\rotk|_Q:Q\to\RR/\ZZ$ is a homomorphism and 
if $\Dd$ is of tube type $\rotk$ is trivial on $Q^\circ$, 
and hence $\rotk(Q)\subset e_G^{-1}\ZZ/\ZZ$.
\end{lemma}

\begin{proof} The first assertion follows from the fact that $\kappa_\RR|_Q=0$
(see Lemma~\ref{lem:8.9}) and from Lemma~\ref{lem:7.1}(2).

Let $Q=MA_QN_Q$ be the Langlands decomposition of $Q$; 
then $\rotk$ is trivial on $A_QN_Q$ (see Proposition~\ref{prop:7.7}(1)).
Now $M^\circ$ is reductive with compact center and 
we may assume that $\Zz(M^\circ)\subset K\cap Q$.
If then $\Dd$ is of tube type, the Lie algebra of $K\cap Q$ is contained
in the Lie algebra of $[K,K]$ \cite[Theorem~4.11]{Koranyi_Wolf_65_annals}
and since $\rotk|_K$ is a homomorphism, it is therefore trivial on $[K,K]$ 
and hence on $\Zz(M^\circ)^\circ$.  Since $\rotk$ is also trivial 
on every connected almost simple factor of $M^\circ$, 
we obtain finally that it is trivial on $Q^\circ$.
\end{proof}

\begin{proof}[{\it Proof of Theorem~\ref{thm_intro:thm9}}] Let $\kappa\in\hhcb^2(G,\ZZ)$ 
and $\rho,\rho_0:\ps\to G$ be maximal representations.
Corollary~\ref{cor:8.6} implies that the real class in $\hb^2\big(\ps,\RR\big)$
which corresponds to $\rho^\ast(\kappa)-\rho(\kappa)\in\hb^2\big(\ps,\ZZ\big)$ vanishes and hence
\bq\label{eq:8.7}
\rho^\ast(\kappa)-\rho_0^\ast(\kappa)=\delta(h)
\eq
for a unique homomorphism $h:\ps\to\RR/\ZZ$. 
Restricting the equality \eqref{eq:8.7} to cyclic subgroups,
we get 
\bqn
\mathrm{R}^{\rho_0}_\kappa(\rho) = \rotk\big(\rho(\gamma)\big)-\rotk\big(\rho_0(\gamma)\big)=h(\gamma)
\eqn
for all $\gamma\in\ps$.  

%Let now $\Delta:L\to G$ be the homomorphism corresponding to a diagonal disk $\DD\to\Xx$,
%where $L$ is some finite covering of $\PU(1,1)$,
%and choose a hyperbolization $h:\ps\to L$ such that $\rot_\alpha\big(h(\gamma)\big)=0$
%for all $\alpha\in\hhcb^2(L,\ZZ)$ and $\gamma\in\ps$.
%Setting $\rho_0=\Delta\circ h$, we have that $\rotk(\rho_0)=0$
%and hence $\rotk(\rho)=h$ is therefore a homomorphism. 

Now from Lemma~\ref{lem:8.8} we know that every $\rho(\gamma)$ and $\rho_0(\gamma)$ is conjugate to an element of $Q$
and thus if $\Dd$ is of tube type we have from Lemma~\ref{lem:8.10}
that $\mathrm{R}^{\rho_0}_\kappa(\rho)\in e_G^{-1}\ZZ/\ZZ$.
The last assertion follows then from the fact that $\rotk$ is continuous and
$\hom\big(\ps,e_G^{-1}\ZZ/\ZZ\big)$ is finite.
\end{proof}

\subsection{Applications to representation varieties}\label{subsec:8.3}
Let $G$ be a group of Hermitian type. 
If $\bs=\emptyset$, it is well known that 
$\hmaxpg$ is a union of components of $\hpg$ and hence if $G$ is real algebraic,
the set of maximal representations is a real semialgebraic set.  In the case in which
$\bs\neq\emptyset$, $\hpg$ is connected;
it is then necessary to study certain naturally defined subsets of the representation variety.

We assume $\bs\neq\emptyset$ and use the presentation in \eqref{eq:ps}.
Let $\Cc=(\Cc_1,\dots,\Cc_n)$ be a set of conjugacy classes in $G$. Then
\bqn
 \hom^\Cc\big(\pi_1(\Sigma),G\big)
:=\big\{\rho\in\hom\big(\pi_1(\Sigma),G\big):\,
\rho(c_i)\in\Cc_i,\,1\leq i\leq n\big\}
\eqn
is a real semialgebraic set and

\begin{cor}\label{cor:8.11} For any $\kappa\in\hcb^2(G,\RR)$ the map
$\rho\mapsto\tksr$ is constant on connected components of $\hom^\Cc\big(\pi_1(\Sigma),G\big)$.
\end{cor}

\begin{proof} Since $\hcb^2(G,\RR)$ is spanned by integral classes 
(see Proposition~\ref{prop:7.6}(3)), we may assume that $\kappa$ is integral
in which case $\tksr$ is congruent $\mod\ZZ$ to $-\sum_{i=1}^n\rotk\big(\rho(\c_i)\big)$;
the latter is then constant for $\rho\in\hom^\Cc\big(\pi_1(\Sigma),G\big)$.
Thus, since $\rho\mapsto\tksr$ is continuous, it is locally constant
which proves the corollary.
\end{proof}

In view of Lemma~\ref{lem:8.8} a particularly suitable space of representations
in relation with the study of maximal representations is
\bqn
 \hom^{\cs}\big(\pi_1(\Sigma),G\big)
:=\big\{\rho\in\hom\big(\pi_1(\Sigma),G\big):\,
  \rho(c_i)\hbox{ has at least}\\ \hbox{ one fixed point in }\cs,\,1\leq i\leq n\big\}\,.
\eqn
Indeed we have:
\bqn
\hmaxpg\subset \hom^{\cs}\big(\pi_1(\Sigma),G\big)\,.
\eqn

\begin{cor}\label{cor:8.12} Let $\kappa\in\hhcb^2(G,\ZZ)$ and assume that $\Dd$ is of tube type.
Then we have that
\bqn
\tksr\in e_G^{-1}\ZZ
\eqn
for every $\rho\in\hom^{\cs}\big(\pi_1(\Sigma),G\big)$ and $\hmaxpg$ is a union 
of connected components of $\hom^{\cs}\big(\pi_1(\Sigma),G\big)$.
In particular, if $G$ is a real algebraic group we conclude that the set of maximal
representations of $\ps$ into $G$ is a real semialgebraic set.
\end{cor}

\begin{proof} If $\Dd$ is of tube type and $\rho(c_i)$ fixes a point in $\cs$,
we have by Lemma~\ref{lem:8.10} that $\rotk\big(\rho(c_i)\big)\in e_G^{-1}\ZZ$
and hence, by Lemma~\ref{lem:8.2}, $\tksr\in e_G^{-1}\ZZ$.
Since $\hhcb^2(G,\ZZ)$ spans $\hcb^2(G,\RR)$, we get that for every $\kappa\in\hcb^2(G,\RR)$,
the map $\rho\mapsto\tksr$ is locally constant on $\hom^{\cs}\big(\pi_1(\Sigma),G\big)$
which implies the assertion.
\end{proof}

%%% Local Variables: 
%%% mode: latex
%%% TeX-master: "toledo"
%%% End: 

%% file: examples.tex
\section{Examples}\label{sec:examples}
The aim of this section is to prove Theorem~\ref{thm_intro:thm4} in the introduction.
In this case $G=\gG(\RR)^\circ$, 
where $\gG$ is a connected algebraic group defined over $\RR$ and $G$ is of Hermitian type.  
As usual, $t:\DD^\rk\to\Xx$ is a maximal polydisk (where $\rk=\rk_\Xx$), 
and $d:\DD\to\Xx$ is the composition of the diagonal embedding $\DD\to\DD^\rk$ with $t$.
Accordingly, we have homomorphisms
\bqn
\tau:\SU(1,1)^\rk\to G\qquad\hbox{ and }\qquad\Delta:\SU(1,1)\to G
\eqn
which satisfy 
\bqn 
\tau^\ast(\kgb)=\kappa_{\SU(1,1)^\rk}^\mathrm{b}
\qquad\hbox{ and }\qquad
\Delta^\ast(\kgb)=\rk\,\kappa_{\SU(1,1)}^\mathrm{b}\,.
\eqn
As a consequence, if $h, h_1,\dots,h_\rk:\ps\to\SU(1,1)$ are hyperbolizations,
the composition of
\bqn
\ba
\ps&\to\SU(1,1)^\rk\\
\gamma&\mapsto\big(h_1(\gamma),\dots,h_\rk(\gamma)\big)
\ea
\eqn
with $\tau$, $h_r$ and $\Delta\circ h$ define maximal representations.
We will need the following

\begin{lemma}\label{lem:9.1} If $\Xx$ is of tube type there exists 
$u\in\Zz_G(\Image\Delta)$ such that $G$ is generated by 
$\Image\tau\cup u(\Image\tau)u^{-1}$.
\end{lemma}

\begin{proof}%[{\it Sketch of the Proof.}]
For every $u\in\Zz_G(\Image\Delta)$ let
$H_u$ denote the subgroup of $G$ generated by $\Image \tau$ and 
$u(\Image \tau)u^{-1}$, and let $\frakh_u$ denote its Lie algebra.
Then $H_u$ is of Hermitian type and of the same rank as $G$, and 
the embedding of the symmetric space $\Yy_u$ associated to $H_u$ into the
symmetric space $\Xx$ associated to $G$ is holomorphic. 
Moreover, if $Z_0$ is the generator of the center of the maximal compact subgroup 
of $\Image\Delta$ which gives the complex structure on the disk $d(\DD)$, 
then $Z_0 \in\frakh_u \subset \frakg$ gives the complex structure on $\Yy_u$ 
and on $\Xx$. In particular, the embedding of Lie algebras 
$\frakh_u \hookrightarrow \frakg$ is an
$(\h_2)$-homomorphism, \cite{Satake_book}. 

Fixing a base point in the image of the tight holomorphic disk $d(\DD)$ 
in $\Xx$,  we may assume that the Cartan decompositions of $\frakh_u$ and
$\frakg$ are compatible.
Let $\frakk_0\subset\frakg$ denote the Lie algebra of $\Zz_G(\Image\Delta)$ 
and ${\mathfrak t}\subset \frakg$ the Lie algebra of $\Image\tau$,
and let ${\mathfrak t} = \frakl \oplus \frakr$ be the Cartan decomposition of 
$\frakt$. With a case by case analysis using the Satake--Ihara classification
of $(\h_2)$-homomorphisms, \cite{Satake_hol_65, Ihara_65}, 
one can determine elements $v\in \frakk_0$ such that 
$\Ad(\exp v)\frakr\cup\frakr$ will not be contained in
any noncompact Lie subalgebra $\frakh \subset \frakg$ given by an
$(\h_2)$-homomorphism. 
\end{proof}

Now let $h:\ps\to\SU(1,1)$ be a hyperbolization as in the statement of the theorem
and choose a simple closed geodesic $C\subset\Sigma^\circ$ separating
$\Sigma$ into two components $\Sigma_1,\Sigma_2$.
With the hypotheses at hand, we can find simple closed geodesics
$C_i\subset\Sigma_i$ not intersecting $C$.
Let $h_t^{(i)}$ denote the hyperbolizations of $\Sigma^\circ$ 
obtained by multiplying the length of $C_i$ by a factor $(1+t)$, $t\geq0$,
while keeping $h_t^{(i)}$ constant on $\pi_1(C)\hookrightarrow\pi_1(\Sigma_i)$.
Fix $0<\epsilon_1<\dots<\epsilon_\rk$; then the composition $\rho_t^{(i)}$ with 
\bqn
\ba
\pi_i(\Sigma_i)&\to\SU(1,1)^\rk\\
\gamma&\mapsto\big(h_{\epsilon_1t}^{(i)}(\gamma),\dots,h_{\epsilon_\rk t}^{(i)}(\gamma)\big)
\ea
\eqn
is maximal and its Zariski closure coincides with $\Image\tau$.
Choose now $u\in\Zz_G(\Image\Delta)$ as in Lemma~\ref{lem:9.1} 
and define the representation $\rho_t:\ps\to G$ by 
\bqn
\ba
\rho_t(\gamma):=
\begin{cases}
\rho_t^{(1)}&\hbox{if}\gamma\in\pi_1(\Sigma_1)\\
u\rho_t^{(2)}u^{-1}&\hbox{if}\gamma\in\pi_1(\Sigma_2)\,.
\end{cases}
\ea
\eqn
Then $\rho_t$ is maximal by the additivity property (see Proposition~\ref{prop:3.1})
and from Lemma~\ref{lem:9.1} we deduce that $\rho_t$ has Zariski dense image for $t>0$.

\vfill
%%% Local Variables: 
%%% mode: latex
%%% TeX-master: "toledo"
%%% End: 

%% file: notation.tex
\appendix

\section{Index of Notation}
\begin{footnotesize}
\begin{tabular}{ll}
$G^\circ$ & connected component of the identity in $G$\\

$\cs$ & Shilov boundary of a bounded symmetric domain \\

$\r_\Xx$ & rank of the symmetric space $\Xx$\\

$G_\Xx$ & connected component of the group of isometries of $\Xx$\\

$\Delta(x,y,z)$ & smooth triangle with geodesic sides and vertices in $x,y,z$\\

$\Delta:L\to G$ & homomorphism associated to a diagonal disk\\

 $\tau:L^{\rk_\Xx}\to G$ & homomorphism associated to a maximal polydisk\\

$\Omega^\bullet(\Xx)^G$ & complex of $G$-invariant differential forms on $\Xx$\\

$\hc^\bullet(G,\RR)$ & continuous cohomology of $G$ with  $\RR$ coefficients\\

$\hcb^\bullet(G,\RR)$& bounded continuous cohomology of $G$ with $\RR$- coefficients\\

$\hhcb^\bullet(G,A)$ & Borel cohomology of $G$ with  $A = \RR,\, \ZZ $ or $\RR/\ZZ$ coefficients\\

$\hhc^\bullet(G,A)$& bounded Borel cohomology of $G$ with  $A = \RR,\, \ZZ $ or $\RR/\ZZ$ coefficients\\

$\h^\bullet(X,Y,A)$ & relative singular cohomology with $A = \RR,\, \ZZ $ or $\RR/\ZZ$ coefficients\\

 $\hb^\bullet(X,Y,A)$& relative bounded singular cohomology with $A = \RR,\, \ZZ $ or $\RR/\ZZ$ coefficients\\

$\h^\bullet(X,A)$ &singular cohomology with $A = \RR,\, \ZZ $ or $\RR/\ZZ$ coefficients\\

 $\hb^\bullet(X,A)$ & bounded singular cohomology with $A = \RR,\, \ZZ $ or $\RR/\ZZ$ coefficients\\

$\h^\bullet\big(\pi_1(X),A\big)$ &group cohomology with $A = \RR,\, \ZZ $ or $\RR/\ZZ$ coefficients\\

 $\hb^\bullet\big(\pi_1(X),A\big)$& bounded group cohomology with $A = \RR,\, \ZZ $ or $\RR/\ZZ$ coefficients\\

$\big(\balt(\cs^\bullet)\big)$ & complex of bounded alternating Borel cocycles on $\cs$\\

$\Zz\linftya\big((\partial\DD)^\bullet,\RR\big)^\Gamma$ & cocycles in the complex of $\G$-invariant alternating\\
& bounded measurable functions on $\partial \DD$ \\

$S_m(Y)$& set of singular $m$-simplices in $Y$ \\

$F_\mathrm{b}(Y,\RR)$ & space of bounded $m$-cochains\\

$\kg$& K\"ahler class in $\hc^2(G,\RR)$\\

$\kgb$& bounded K\"ahler class in $\hcb^2(G,\RR)$\\

$\tksr$ & Toledo number with respect to $\kappa \in \hcb^2(G,\RR)$\\

$\tsr$ & Toledo number with respect to bounded K\"ahler class $\kgb \in \hcb^2(G,\RR)$\\

$\rot_\kappa$ & rotation number associated to  $\kappa \in \hhcb^2(G,\ZZ)$\\

$\widetilde{\rot_\kappa}$&  homogeneous quasimorphism which lifts $\rot_\kappa$\\

$\Nn_K(A)$ & normalizer of $A$ in $K$\\

$\Zz_K(A)$ & centralizer of $A$ in $K$\\

$\Zz(\Gamma)$ & center of $\Gamma$\\

%$\kgib$ & basis elements of $\hcb^2(G,\RR)$ coming from a decomposition of $G$ into almost simple factors

\end{tabular}
\end{footnotesize}

%% file: BurgerIozziWienhard.bbl
\def\cprime{$'$} \def\cprime{$'$} \def\cprime{$'$} \def\cprime{$'$}
  \def\cprime{$'$}
\providecommand{\bysame}{\leavevmode\hbox to3em{\hrulefill}\thinspace}
\providecommand{\MR}{\relax\ifhmode\unskip\space\fi MR }
\providecommand{\MRhref}[2]{%
  \href{http://www.ams.org/mathscinet-getitem?mr=#1}{#2}
}
\providecommand{\href}[2]{#2}
\begin{thebibliography}{10}

\bibitem{Barge_Ghys}
J.~Barge and {\'E}.~Ghys, \emph{Cocycles d'{E}uler et de {M}aslov}, Math.~Ann.
  \textbf{294} (1992), no.~2, 235--265.

\bibitem{Blanc}
Ph. Blanc, \emph{Sur la cohomologie continue des groupes localement compacts},
  Ann.~Sci.~\'Ecole Norm.~Sup.~(4) \textbf{12} (1979), no.~2, 137--168.

\bibitem{Borel_comm}
A.~Borel, \emph{Class functions, conjugacy classes and commutators in
  semisimple {L}ie groups}, Algebraic groups and {L}ie groups,
  Austral.~Math.~Soc.~Lect.~Ser., vol.~9, Cambridge Univ.~Press, Cambridge,
  1997, pp.~1--19.

\bibitem{Bradlow_GarciaPrada_Gothen}
S.~B. Bradlow, O.~Garc{\'{\i}}a-Prada, and P.~B. Gothen, \emph{Surface group
  representations in {${\rm PU}(p,q)$} and {H}iggs bundles}, J.~Diff.~Geom.
  \textbf{64} (2003), no.~1, 111--170.

\bibitem{Bradlow_GarciaPrada_Gothen_cacanerveuse}
\bysame, \emph{Maximal surface group representations in isometry groups of
  classical {H}ermitian symmetric spaces}, Geom.~Dedicata \textbf{122} (2006),
  185--213.

\bibitem{Brooks}
R.~Brooks, \emph{Some remarks on bounded cohomology}, Riemann surfaces and
  related topics: Proceedings of the 1978 Stony Brook Conference (State
  Univ.~New York, Stony Brook, N.Y., 1978) (Princeton, N.J.), Ann.~of
  Math.~Stud., vol.~97, Princeton Univ.~Press, 1981, pp.~53--63.

\bibitem{Burger_Iozzi_Labourie_Wienhard}
M.~Burger, F.~Labourie A.~Iozzi, and A.~Wienhard, \emph{Maximal representations
  of surface groups: {S}ymplectic {A}nosov structures}, Quaterly Journal of
  Pure and Applied Mathematics \textbf{1} (2005), no.~3, 555--601, Special
  Issue: In Memory of Armand Borel, Part 2 of 3.

\bibitem{Burger_Iozzi_formula}
M.~Burger and A.~Iozzi, \emph{A useful formula in bounded cohomology}, to
  appear in "S\'eminaires et Congr\`es", nr.~18,
  http://www.math.ethz.ch/$\sim$iozzi/grenoble.ps.

\bibitem{Burger_Iozzi_app}
\bysame, \emph{Boundary maps in bounded cohomology}, Geom.~Funct.~Anal.
  \textbf{12} (2002), 281--292.

\bibitem{Burger_Iozzi_supq}
\bysame, \emph{Bounded {K}\"ahler class rigidity of actions on {H}ermitian
  symmetric spaces}, Ann.~Sci.~\'Ecole Norm.~Sup.~(4) \textbf{37} (2004),
  no.~1, 77--103.

\bibitem{Burger_Iozzi_difff}
\bysame, \emph{Bounded differential forms, generalized {M}ilnor--{W}ood
  inequality and an application to deformation rigidity}, Geom. Dedicata
  \textbf{125} (2007), no.~1, 1--23.

\bibitem{Burger_Iozzi_Wienhard_tight}
M.~Burger, A.~Iozzi, and A.~Wienhard, \emph{Tight homomorphisms and hermitian
  symmetric spaces}, to appear in Geom.~Funct.~Anal. {\sf
  http://www.arXiv.org/math.DG/0710.5641}.

\bibitem{Burger_Iozzi_Wienhard_ann}
\bysame, \emph{Surface group representations with maximal {T}oledo invariant},
  C.~R.~Acad.~Sci.~Paris, S\'er.~I \textbf{336} (2003), 387--390.

\bibitem{Burger_Iozzi_Wienhard_kahler}
\bysame, \emph{Hermitian symmetric spaces and {K}\"ahler rigidity},
  Transform.~Groups \textbf{12} (2007), no.~1, 5--32.

\bibitem{Burger_Monod_GAFA}
M.~Burger and N.~Monod, \emph{Continuous bounded cohomology and applications to
  rigidity theory}, Geom.~Funct.~Anal. \textbf{12} (2002), 219--280.

\bibitem{Clerc_maslov_tube}
J.~L. Clerc, \emph{L'indice de {M}aslov g\'en\'eralis\'e}, J.~Math.~Pures
  Appl.~(9) \textbf{83} (2004), no.~1, 99--114.

\bibitem{Clerc_maslov}
\bysame, \emph{An invariant for triples in the {S}hilov boundary of a bounded
  symmetric domain}, Comm. Anal. Geom. \textbf{15} (2007), no.~1, 147--173.

\bibitem{Clerc_Koufany}
J.~L. Clerc and K.~Koufany, \emph{Primitive du cocycle de {M}aslov
  g\'en\'eralis\'e}, Math. Ann. \textbf{337} (2007), 91--138.

\bibitem{Clerc_Neeb}
J.~L. Clerc and K.-H. Neeb, \emph{Orbits of triples in the {S}hilov boundary of
  a bounded symmetric domain}, Transform.~Groups \textbf{11} (2006), no.~3,
  387--426.

\bibitem{Clerc_Orsted_TG}
J.~L. Clerc and B.~{\O}rsted, \emph{The {M}aslov index revisited},
  Transform.~Groups \textbf{6} (2001), no.~4, 303--320.

\bibitem{Clerc_Orsted_2}
\bysame, \emph{The {G}romov norm of the {K}aehler class and the {M}aslov
  index}, Asian J.~Math. \textbf{7} (2003), no.~2, 269--295.

\bibitem{Fock_Goncharov_local}
V.~Fock and A.~Goncharov, \emph{Moduli spaces of local systems and higher
  {T}eichm\"uller theory}, Publ.~Math.~Inst.~Hautes \'Etudes Sci. (2006),
  no.~103, 1--211.

\bibitem{Fock_Goncharov_convex}
\bysame, \emph{Moduli spaces of convex projective structures on surfaces},
  Adv.~Math. \textbf{208} (2007), no.~1, 249--273.

\bibitem{Furstenberg_63}
H.~Furstenberg, \emph{A {P}oisson formula for semisimple {L}ie groups}, Ann.~of
  Math. \textbf{77} (1963), 335--383.

\bibitem{Ghys_87}
{\'E}.~Ghys, \emph{Groupes d'hom\'eomorphismes du cercle et cohomologie
  born\'ee}, The {L}efschetz centennial conference, Part III, (Mexico City
  1984), Contemp.~Math., vol.~58, American Mathematical Society, RI, 1987,
  pp.~81--106.

\bibitem{Goldman_thesis}
W.~M. Goldman, \emph{Discontinuous groups and the {E}uler class}, Thesis,
  University of California at Berkeley, 1980.

\bibitem{Goldman_88}
\bysame, \emph{Topological components of spaces of representations},
  Invent.~Math. \textbf{93} (1988), no.~3, 557--607.

\bibitem{Goldman_torus}
\bysame, \emph{The modular group action on real {${\rm SL}(2)$}-characters of a
  one-holed torus}, Geom.~Topol. \textbf{7} (2003), 443--486 (electronic).

\bibitem{Gothen}
P.~B. Gothen, \emph{Components of spaces of representations and stable
  triples}, Topology \textbf{40} (2001), no.~4, 823--850.

\bibitem{Gromov_82}
M.~Gromov, \emph{Volume and bounded cohomology}, Inst.~Hautes \'Etudes
  Sci.~Publ.~Math. \textbf{56} (1982), 5--99.

\bibitem{Guichard_hyperconvex}
O.~Guichard, \emph{Sur les repr\'esentations des groupes de surface}, preprint
  2004.

\bibitem{Hernandez}
L.~{H}ern\'andez Lamoneda, \emph{Maximal representations of surface groups in
  bounded symmetric domains}, Trans.~Amer.~Math.~Soc. \textbf{324} (1991),
  405--420.

\bibitem{Hitchin}
N.~J. Hitchin, \emph{Lie groups and {T}eichm\"uller space}, Topology
  \textbf{31} (1992), no.~3, 449--473.

\bibitem{Ihara_65}
S.~Ihara, \emph{Holomorphic imbeddings of symmetric domains},
  J.~Math.~Soc.~Japan \textbf{19} (1967), 261--302. \MR{35 \#5656}

\bibitem{Ivanov}
N.~V. Ivanov, \emph{Foundations of the theory of bounded cohomology}, J.~of
  Soviet Mathematics \textbf{37} (1987), no.~1, 1090--1115.

\bibitem{Kaimanovich_ern}
V.~A. Kaimanovich, \emph{S{AT} actions and ergodic properties of the horosphere
  foliation}, Rigidity in dynamics and geometry (Cambridge, 2000), Springer,
  Berlin, 2002, pp.~261--282.

\bibitem{Koranyi_Wolf_65_annals}
A.~Kor{\'a}nyi and J.~A. Wolf, \emph{Realization of hermitian symmetric spaces
  as generalized half-planes}, Ann.~of Math.~(2) \textbf{81} (1965), 265--288.

\bibitem{Koziarz_Maubon}
V.~Koziarz and J.~Maubon, \emph{Harmonic maps and representations of
  non-uniform lattices of {${\rm PU}(m,1)$}}, preprint, {\sf
  arXiv:math.DG/0309193}.

\bibitem{Labourie_energy}
F.~Labourie, \emph{Cross {R}atios, {A}nosov representations and the energy
  functional on {T}eichm\"uller space}, Annales Scientifiques de l'ENS, to
  appear, {\sf http://www.arXiv.org/math.DG/0512070}.

\bibitem{Labourie_anosov}
\bysame, \emph{{A}nosov flows, surface groups and curves in projective space},
  Invent. Math. \textbf{165} (2006), no.~1, 51--114.

\bibitem{Strohm}
C.~Loeh-Strohm, \emph{The proportionality principle of simplicial volume},
  preprint, {\sf http://www.arXiv.org/math.AT/0504106}, April 2005.

\bibitem{Lusztig_positivity}
G.~Lusztig, \emph{Total positivity in reductive groups}, Lie theory and
  geometry, Progr.~Math., vol. 123, Birkh\"auser Boston, Boston, MA, 1994,
  pp.~531--568.

\bibitem{Lusztig}
\bysame, \emph{Total positivity in partial flag manifolds}, Represent.~Theory
  \textbf{2} (1998), 70--78 (electronic).

\bibitem{Mackey_57}
G.~W. Mackey, \emph{Les ensembles bor\'eliens et les extensions des groupes},
  J.~Math.~Pures Appl.~(9) \textbf{36} (1957), 171--178.

\bibitem{Milnor}
J.~Milnor, \emph{On the existence of a connection with curvature zero},
  Comment.~Math.~Helv. \textbf{32} (1958), 215--223.

\bibitem{Monod_book}
N.~Monod, \emph{Continuous bounded cohomology of locally compact groups},
  Lecture Notes in Math., no. 1758, Springer-Verlag, 2001.

\bibitem{Satake_hol_65}
I.~Satake, \emph{Holomorphic imbeddings of symmetric domains into a {S}iegel
  space}, Amer.~J.~Math. \textbf{87} (1965), 425--461. \MR{33 \#4326}

\bibitem{Satake_book}
\bysame, \emph{Algebraic structures of symmetric domains}, Kan\^o Memorial
  Lectures, vol.~4, Iwanami Shoten, Tokyo, 1980.

\bibitem{Toledo_89}
D.~Toledo, \emph{Representations of surface groups in complex hyperbolic
  space}, J.~Diff.~Geom. \textbf{29} (1989), no.~1, 125--133.

\bibitem{Van_Est}
W.~T. van Est, \emph{Group cohomology and {L}ie algebra cohomology in {L}ie
  groups, {I}, {II}}, Nederl.~Akad.~Wetensch.~Proc.~Series
  A.~\{56\}=Indag.~Math. \textbf{15} (1953), 484--504.

\bibitem{Wienhard_mapping}
A.~Wienhard, \emph{The action of the mapping class group on maximal
  representations}, Geom.~Dedicata \textbf{120} (2006), 179--191.

\bibitem{Wigner_73}
D.~Wigner, \emph{Algebraic cohomology of topological groups},
  Trans.~Amer.~Math.~Soc. \textbf{178} (1973), 83--93.

\end{thebibliography}
